\documentclass{amsart}

\usepackage[utf8]{inputenc} 
\usepackage[T1]{fontenc}    
\usepackage{hyperref}       
\usepackage{url}            
\usepackage{booktabs}       
\usepackage{amsfonts}       
\usepackage{nicefrac}       
\usepackage{microtype}      
\usepackage{xcolor}         

\usepackage{amssymb}   
\usepackage{algorithm}
\usepackage{algpseudocode}  
\usepackage{amsthm}
\usepackage{graphicx}
\usepackage{stmaryrd}
\usepackage{nccmath}
\usepackage{geometry}
\usepackage{amsmath}
\usepackage{enumitem}
\usepackage{amsaddr}
\newcommand{\R}{\mathbb{R}}
\newcommand{\C}{\mathbb{C}}
\newcommand{\N}{\mathbb{N}}
\newcommand{\E}{\mathbb{E}}

\makeatletter
\def\namedlabel#1#2{\begingroup
    #2%
    \def\@currentlabel{#2}%
    \phantomsection\label{#1}\endgroup
}
\makeatother

\usepackage{physics}

\DeclareMathOperator*{\Diag}{Diag}

\DeclareMathOperator*{\Id}{Id}

\theoremstyle{plain}
\newtheorem{thm}{Theorem}[section]
\newtheorem{lem}{Lemma}
\newtheorem{crl}{Corollary}
\newtheorem{ppt}{Proposition}
\theoremstyle{definition}
\newtheorem{dft}{Definition}
\newtheorem{as}{Assumption}
\newtheorem{ex}{Example}
\newtheorem{rmk}{Remark}
\newtheorem{notat}{Notation}

\begin{document}

\title[Asymptotic non-linear shrinkage for weighted sample covariance]{Asymptotic non-linear shrinkage and eigenvector overlap for weighted sample covariance}

\author[]{Benoît Oriol}
\address{CEREMADE, Université Paris-Dauphine, PSL, Paris, France}
\address{Société Générale Corporate and Investment Banking, Puteaux, France}

\date{March 17, 2025}

\begin{abstract}We compute asymptotic non-linear shrinkage formulas for covariance and precision matrix estimators for weighted sample covariances, and the joint sample-population eigenvector overlap distribution, in the spirit of Ledoit and Péché. We detail explicitly the formulas for exponentially-weighted sample covariances. We propose an algorithm to numerically compute those formulas. Experimentally, we show the performance of the asymptotic non-linear shrinkage estimators. Finally, we test the robustness of the theory to a heavy-tailed distributions.\end{abstract}



\maketitle

\section*{Introduction}
Covariance estimation is a central topic in multivariate analysis. In high dimension and proportional number of samples, known as Kolmogorov asymptotics, the sample covariance behaves badly due to an inherent spectrum deformation. This phenomenon is explained by Random Matrix Theory (RMT), originally by the Marcenko-Pastur theorem \cite{Marcenko1967}.

One method to cope with this phenomenon and unbias the sample covariance spectrum is called shrinkage \cite{Ledoit2020a, Clifford2020}: the idea is to multiply (to "shrink" in the literature) each sample covariance eigenvalue by a specific factor to counter the deformation induced by high dimension. The idea goes back to Stein \cite{Stein1956}.

Linear shrinkage (where we shrink each eigenvalue by the same constant) was deeply studied, \cite{Ledoit2004, Chen2010, Ikeda2016, Yang2019, Oriol2025} to cite some. Then, RMT gave the mathematical tools to understand asymptotically non-linear shrinkage. Ledoit and Péché \cite{Ledoit2009} found asymptotic equations for optimal non-linear shrinkage of the sample covariance in the class of rotation-invariant covariance and precision matrix estimators. In the case of covariance estimation, the formulas were extended for more general noise models, additive and multiplicative, in \cite{Bun2016}.

This work opened a way to several methods to estimate those complex non-linear optimal formulas which uses the asymptotic Cauchy-Stieltjes transform of the sample spectrum: estimation through cross-validation \cite{Abadir2011}, methods of moments \cite{Rao2008, Bai2010b, Li2013, Kong2017}, discrete Marcenko-Pastur inversion with QuEST \cite{Karoui2006, Ledoit2015, Ledoit2016}, or kernel density estimation \cite{Ledoit2019, Ledoit2020b} thanks to the theoretical tools on kernel sample covariance spectrum density estimation studied by Jing in 2010 \cite{Jing2010}. The applications of those estimators are diverse, and an extensive list was made by Ledoit and Wolf \cite{Ledoit2016} in their literature review. We can include new works in it, as climatology \cite{Popov2022}, neuroscience \cite{Honnorat2019}, or sensor monitoring \cite{Steland2018}. More recently, a CLT on non-linear shrinkage formulas was used to design a statistical test on the covariance structure \cite{Bodnar2024}, and a general study on the asymptotic deviation of a larger set of functionals needed for non-linear shrinkage is studied in \cite{Louart2023}.

These methods focus on the standard sample covariance. However, applications in time series such as neurosciences \cite{Honnorat2022}, finance \cite{Bun2016, Ding2022, Bongiorno2023}, suggest that the sample covariance, even shrunk, suffers from the non-stationarity of the data. The weighted sample covariance appears naturally in finance when the volatility is time-dependent \cite{Bun2016}. Moreover, weighting schemes, such as the exponential weighted moving average (EWMA), are a model-free approach, and represent a transparent candidate for complex and hard-to-model dynamics. The weighted sample covariance was used in a shrinkage setup under Gaussian setting and EWMA weights in \cite{Pafka2004} for portfolio management. This was further studied for covariance filtering \cite{Daly2010}. The spectrum of the EWMA sample covariance in an uncorrelated setting was studied in \cite{Svensson2007}, and Tan et al. in 2023 \cite{Tan2023} developed a shrinkage algorithm of EWMA sample covariance, \textit{à la} NERCOME, with applications in dynamic brain connectivity.

Theoretically, the asymptotic spectrum of the weighted sample covariance has been studied, particularly in the wireless communication field, see \cite{Tulino2004, Bai2010} for two different reviews. The convergence of the spectrum of the weighted sample covariance in high dimension was studied in the model of MIMO systems for example \cite{Couillet2011}. The result can be also found in the literature under the name of separable covariance matrices or random matrices with correlated inputs \cite{Monvel1996, Burda2005, Paul2009}. The version of the spectral convergence theorem with the weakest set of hypotheses is stated in \cite{Zhang2007}. 

The support of the asymptotic spectrum and the presence or not of eigenvalues outside the asymptotic support in finite samples was studied in the case of the unweighted sample covariance \cite{Bai1998}, and used in practice in the QuEST algorithm \cite{Ledoit2015, Ledoit2016}. For the weighted sample covariance and more general types of product random matrices, the subject was addressed more recently \cite{Couillet2015, Bai2023, Ding2023}.


However, to the best of our knowledge, no extension to formulas for asymptotic non-linear shrinkage of the weighted sample covariance have been derived yet, as it has been done for the unweighted sample covariance. A rotation-invariant covariance estimator has been derived in \cite{Bun2016}, but the result is difficultly tractable, without direct way to implement it in practice. This work addresses this gap in the literature: we generalize the asymptotic equation of the spectrum in the spirit of Ledoit and Péché \cite{Ledoit2009} and use it to derive asymptotic non-linear shrinkage formulas. Those formulas give tractable asymptotic optimal rotation-invariant covariance and precision matrix estimators. 
We compute analytical formulas in function of the asymptotic Cauchy-Stieltjes transform in the case of the exponentially-weighted sample covariance.

Similarly to the unweighted scenario, the computation of those formulas is a challenge, as it depends on the asymptotic Cauchy-Stieltjes transform of the sample spectrum. Here, we propose an algorithm to compute the non-linear shrinkage formulas using samples and weights only, making \textit{bona fide} estimators of the covariance and precision matrices. It uses auto-differentiation to make the implementation simple and fast.

The experimental part numerically highlights the performance of the exact non-linear shrinkage. Moreover, we look at how they behave when the theoretical assumptions are not met, specifically when the underlying noise distribution has heavy tails.

\section{Notation, definitions and hypotheses}
Notation is not constant across major works on the spectrum of sample covariances and their use for shrinkage. In our work, we follow mostly Silverstein \cite{Silverstein1995a}.

\begin{notat}[The data matrix]
	There are $N$ samples of dimension $n$. We have:
	 \begin{itemize}
	 	\item $c_n = \frac{n}{N}$ the concentration ratio,
	 	\item $Z_n$ is the noise $n \times N$ matrix composed of independent centered complex entries of unit variance, 
	 	\item $T_n$ is the true covariance, a non-negative definite Hermitian matrix of size $n \times n$, 
	 	\item $W_n$ is the weight matrix, an $N \times N$ diagonal non-negative real matrix, 
		\item $Y_n = T_n^{1/2} Z_n $ is the observed data matrix.  
	\end{itemize}
\end{notat}

Remark that the mean of the noise distribution is assumed to be zero. In practice, if the mean is unknown, we can center the samples using the empirical mean. As noted in \cite{Ledoit2020b}, Remark 3, demeaning the samples does not change the asymptotic spectrum properties of the (weighted) sample covariance, as it is an operation of rank one. So for the theoretical asymptotic part, without loss of generality, we can assume that the mean is zero.

A practical heuristic to take in account this statistical operation in finite samples is given in Section 7 \cite{Ledoit2020a}: the authors propose to replace the number of samples $N$ by $N-1$ in the implementation if the samples were empirically centered.

The object of interest is the \textbf{weighted} sample covariance $B_n$ and its spectrum.
\begin{notat}[Weighted sample covariance]
	For $n \in \N^*$, the weighted sample covariance is defined by:
	\begin{equation*}\label{}
	\begin{aligned}
		B_n := \frac{1}{N} Y_n W_n Y_n^*.
	\end{aligned}
	\end{equation*} 
	Moreover, we denote by $(\tau_1^{(n)},...,\tau_n^{(n)})$ the eigenvalues of $T_n$ in decreasing order, and $(\lambda_1^{(n)},...,\lambda_n^{(n)})$ the eigenvalues of $B_n$ in decreasing order.
\end{notat}

\begin{ex}[Standard weighting]
	$W_n = I_{N(n)}$ is the standard weighting. In this situation, $B_n$ is the standard sample covariance, and its asymptotic spectrum is described by Marcenko and Pastur \cite{Marcenko1967}.
\end{ex}

\begin{ex}[Exponentially weighted scheme]\label{exEWMA}
	In time series analysis, a common choice is the exponentially weighted (EWMA) scheme. Parametrized by some $\alpha \in \R_+^*$, we define the weights as:
	\begin{equation*}\label{}
	\begin{aligned}
		\forall i \in \llbracket 1,N \rrbracket, (W_n)_{ii} = \beta e^{-\alpha i/N}, \\
		\beta = e^{-\alpha /N} \frac{1 - e^{-\alpha /N}}{1 - e^{-\alpha}}.
	\end{aligned}
	\end{equation*}
\end{ex}

The spectrum of $B_n$ is studied through its empirical spectrum distribution (e.s.d.).
\begin{notat}[Empirical spectrum distribution]
	We consider a Hermitian matrix $A$ of size $n \times n$ with real eigenvalues $(\mu_1,...,\mu_n)$. We define the empirical spectrum distribution of $A$, denoted $F^A$, as:
	\begin{equation*}\label{}
	\begin{aligned}
		F^A := \frac{1}{n} \sum_{i=1}^n \mathrm{1}_{[\mu_i, +\infty[}.
	\end{aligned}
	\end{equation*}
\end{notat}

We describe now several assumptions, mostly the same used in Ledoit and Péché \cite{Ledoit2009}. These assumptions fix the framework of what we call "high dimensional setting": the dimension and number samples grow linearly together, and the empirical spectrum distribution converges in this setting. 

The assumption on bounded $12^{th}$ moments to be finite may appear quite strong for practitioner, however our experimental work lead to think that only bounded $4^{th}$ moments are enough for the results to be practically applicable.
\begin{as}
	We assume the following hypotheses.
	\begin{itemize}
		\item[\namedlabel{H1}{H1}:] $Z_n$ is a $(n,N)$ matrix of real or complex independent random variables with zero mean, unit variance and bounded $12^{th}$ moments by a constant $B \in \R$ independent of $n$.
		\item[\namedlabel{H2}{H2}:] $T_n$ is a random Hermitian positive definite matrix, $W_n$ is a diagonal random positive definite matrix, and $Z_n$, $W_n$ and $T_n$ are mutually independent.
		\item[\namedlabel{H3}{H3}:] $c_n = \frac{n}{N} \rightarrow c \in \R_+^*$ as $n \rightarrow \infty$.
		\item[\namedlabel{H4}{H4}:] $F^{T_n} \underset{n \rightarrow \infty}\implies H$ almost surely, where $\implies$ denotes weak convergence. $H$ defines a probability distribution function (p.d.f.), whose support $S_{F^{T_n}}$ is included in the compact interval $[h_1, h_2]$ with $0 < h_1 \leq h_2 < \infty$.
		\item[\namedlabel{H5}{H5}:] $F^{W_n} \underset{n \rightarrow \infty}\implies D$ a.s.. $D$ defines a probability distribution function, whose support $S_{F^{W_n}}$ is included in the compact interval $[d_1, d_2]$ with $0 < d_1 \leq d_2 < \infty$. 
	\end{itemize}
\end{as}

\begin{rmk}
	Remark that Assumption \ref{H1} requires only the independence of the $(Z_n)_{i,j}$, not that they are identically distributed. This is the main difference from the Assumptions made in \cite{Ledoit2009} for the unweighted case.
\end{rmk}

\section{$\Theta^g$ characterization, non-linear shrinkage oracle formulas and overlap distribution}
In this section, we find asymptotic optimal formulas to correct the spectrum of weighted sample covariance and precision matrices in high dimension, and the asymptotic density of the sample-population eigenvector overlap. A central result in this study is the extension of Marcenko-Pastur theorem to weighted covariance matrices. This theorem states that, in the high-dimensional setting, the empirical spectrum $F^{B_n}$ of the weighted sample covariance almost surely converges weakly to a deterministic spectrum $F$. This result can be found in different works under different sets of hypotheses \cite{Monvel1996, Burda2005, Zhang2007, Paul2009}, the weakest assumptions being found in Theorem 1.2.1, \cite{Zhang2007}. With our notation and assumptions, which are not optimal for this result, the theorem can be stated as following:

\begin{thm}[$F^{B_n}$ convergence, from Theorem 1.2.1 \cite{Zhang2007}]\label{wMP}
	Assume the conditions \ref{H1}-\ref{H5}. Then, almost surely, $F^{B_n} \implies F$, where $F$ is a deterministic distribution depending only on $c, H$ and $D$ and is characterized by its Cauchy-Stieltjes transform, denoted $m(\cdot)$. For all $z \in \C_+ := \{ z \in \C|\Im[z] > 0\}$:
	\begin{equation}\label{}
	\begin{aligned}
		 &m(z) = -\frac{1}{z}\int \frac{1}{\tau X(z) + 1}dH(\tau),
	\end{aligned}
	\end{equation}
	where for all $z \in \C_+$, $X(z)$ is the unique solution in $\C_+$ of the following equation:
	\begin{equation}\label{}
	\begin{aligned}
		X(z) &=-\int \frac{\delta}{z -  \delta c \int \frac{\tau}{\tau X(z) + 1}dH(\tau)}dD(\delta).
	\end{aligned}
	\end{equation}
\end{thm}

However, when it comes to minimize the loss of a covariance or precision matrix estimator, the behavior of the spectrum $F^{B_n}$ is not sufficient anymore as the eigenvectors plays also a role in the loss. 
\begin{notat}[Eigenvectors of $B_n$ and $T_n$]
	We denote $(u_1^{(n)},...,u_n^{(n)})$ a set of eigenvectors of $B_n$ associated to the eigenvalues $(\lambda_1^{(n)},...,\lambda_n^{(n)})$, and $(v_1^{(n)},...,v_n^{(n)})$ a set of eigenvectors of $T_n$ associated to the eigenvalues $(\tau_1^{(n)},...,\tau_n^{(n)})$.
\end{notat}

In this section, we look at functionals $\Theta_n^{(g)}$ introduced in a work of Ledoit and Péché \cite{Ledoit2009} that carry the information needed to shrink properly the weighted sample eigenvalues for covariance and precision matrix estimation, taking in account the non-ideal projection of $u$ onto $v$. For that, we need a notation to apply real functions to $\R$-diagonalizable matrices.
\begin{notat}[Extension of real functions to $\R$-diagonalizable matrices]
	With $g: \R \rightarrow \R$ and $M = P\Lambda P^{-1}$, $M$ a $\R$-diagonalizable $p \times p$ matrix, $p \in \N^*$, $\Lambda$ a real diagonalization of $M$, we denote:
	\begin{equation*}\label{}
	\begin{aligned}
		 g(M) := P\tilde \Lambda P^{-1},
	\end{aligned}
	\end{equation*}
	where $\forall i \neq j \in \llbracket 1,p \rrbracket, \tilde \Lambda_{ii} = g(\Lambda_{ii}) \text{ and } \tilde \Lambda_{ij} = 0$.
\end{notat}

We can now define the functionals $\Theta_n^{(g)}$ which are the objects of interest of this section.
\begin{dft}[$\Theta^g_n$]\label{thetadef}
For $g: \R \rightarrow \R$ a bounded function with a finite number of discontinuity, we define:
\begin{equation}\label{}
\begin{aligned}
	 \forall z \in \C_+, \Theta_n^g(z) &= \frac{1}{n} \sum_{i=1}^n \frac{1}{\lambda_i - z} \sum_{j=1}^n |u_i^* v_j|^2 g(\tau_j) \\
	 &= \frac{1}{n} \tr\left( (B_n - zI)^{-1}g(T_n)\right).
\end{aligned}
\end{equation}
For $k \in \mathbb{Z}$, we denote $\Theta^{(k)}_n = \Theta^g_n$ for $g: t \mapsto t^k$.
\end{dft}

\subsection{Asymptotic equation on $\Theta^g$}
In the spirit of Theorem 2 by Ledoit and Péché \cite{Ledoit2009} for equally weighted sample covariance, the following result characterizes the limit of the functionals $\Theta^g_n$ defined above. 
\begin{thm}\label{FEg}
	Assume the conditions \ref{H1}-\ref{H5}. For $g: \R \rightarrow \R$ a bounded function with a finite number of discontinuity, for all $z \in \C_+$, $\Theta^{g}_n(z) \rightarrow \Theta^{g}(z)$ almost surely and:
	\begin{equation}\label{}
	\begin{aligned}
		 &\Theta^{g}(z) = -\frac{1}{z}\int \frac{g(\tau)}{\tau X(z) + 1}dH(\tau),
	\end{aligned}
	\end{equation}
	where for all $z \in \C_+$, $X(z)$ is the unique solution in $\C_+$ of the following equation:
	\begin{equation}\label{}
	\begin{aligned}
		X(z) &=-\int \frac{\delta}{z -  \delta c \int \frac{\tau}{\tau X(z) + 1}dH(\tau)}dD(\delta).
	\end{aligned}
	\end{equation}
\end{thm}
The proof of this result is given in the Appendix \ref{apxB}. As remarked in Ledoit and Péché \cite{Ledoit2009}, for all $z \in \C_+$, the kernel of integration $\kappa_z(\tau) = \frac{1}{\tau X(z) + 1}$ is conserved disregarding the choice of $g$.

The next result controls the behavior of $X$ and $\Theta^g$ near the real line, and is essential to derive oracle shrinkage formulas for rotation-invariant precision and covariance matrix estimators. The continuity on $\R^*$ is known for $m$ and implies that $F$ has a density on $\R^*$ \cite{Couillet2015}, we generalize the result for the functionals $\Theta^g$.
\begin{thm}[Continuity on the real line]\label{c0}
	Assume the conditions \ref{H1}-\ref{H5}. Then, for $\lambda \in \R$ if $c > 1$ or $\lambda \in \R^*$ otherwise, we have that $\check{X}(\lambda) := \lim_{z \in \C_+ \rightarrow \lambda} X(z)$ exists. For $g: [h_1,h_2] \rightarrow \R$ bounded with a finite number of discontinuity points, for $\lambda \in \R$ if $c < 1$ or $\lambda \in \R^*$ otherwise, $\check \Theta^{g}(\lambda) := \lim_{z \in \C_+ \rightarrow \lambda} \Theta^{g}(z)$ exists and for $\lambda \in \R^*$:
	\begin{equation}\label{thetag}
	\begin{aligned}
		&\check \Theta^{g}(\lambda) = -\frac{1}{\lambda} \int \frac{g(\tau)}{\tau \check X(\lambda) + 1}dH(\tau).
	\end{aligned}
	\end{equation}
	If $c < 1$, we have furthermore:
	\begin{equation}
	\begin{aligned}
		&\check \Theta^{g}(0) = \frac{\int \check m(0)g(\tau)\tau^{-1}dH(\tau)}{\int \tau^{-1} dH(\tau)} \text{ with } \check m(\cdot) := \check \Theta^{(0)}(\cdot).
	\end{aligned}
	\end{equation}
\end{thm}

The proof of this result is given in the Appendix \ref{apxC}. The equation \eqref{thetag} of Theorem \ref{c0} is of utter importance in the derivation and computations of the asymptotic non-linear shrinkage formulas: it reduces the complexity of computing the $\check \Theta^g$, central objects in non-linear shrinkage formulas as it will be stated in Section \ref{formsection}, to only two objects ($\check X$ and $H$), which are both retrievable from the samples, as we will see in the Section \ref{numsection}.

The next part aims at using this result in the setting of rotation-invariant covariance and precision matrix estimation. By properly choosing $g$, with $x \mapsto x$ and $x \mapsto 1/x$, we will derive oracle shrinkage formulas for weighted estimators.

\subsection{Oracle shrinkage formulas for rotation-invariant covariance matrix estimators}\label{formsection}
An important application of the Theorem \ref{FEg} is the possibility to derive formulas for the optimal shrinkage of weighted sample eigenvalues in the class of rotation-invariant estimators, as it was done in the non-weighted case in Theorem 4 and 5 by Ledoit and Péché \cite{Ledoit2009}. The minimization is asymptotically equivalent under Frobenius, Inverse Stein or Minimum Variance loss \cite{Ledoit2020b}. 
\begin{dft}[Class of rotation-invariant estimator]
	Given a weighted sample covariance $S_n$ of size $(n,n), n \in \N^*$, of the diagonalized form $U_n \Lambda_n U_n^*$ with $\Lambda$ diagonal, real, we define the class of rotation-invariant estimators as:
	\begin{equation*}\label{}
	\begin{aligned}
		 \mathcal{C}_n = \left\{U_nD_nU_n^* | D_n = \Diag(d_1,...,d_n), (d_i)_{i=1}^n \in \R^n \right\}.
	\end{aligned}
	\end{equation*}
\end{dft}
This class of estimators allows modifying (shrinking) the eigenvalues of the weighted sample covariance estimator, but not the eigenvectors. The performance of an estimator in $\mathcal{C}_n$ is thus directly linked to its capacity to correct the weighted sample covariance spectrum. Remark that the class of rotation invariant estimators suits to covariance matrix estimation as well as precision matrix estimation.

\subsubsection{Asymptotic shrinkage for covariance matrix estimation}
We first focus on the problem covariance estimation. It is very similar regarding the method to precision matrix estimation in our framework, which will be exposed after this part.

Under Frobenius norm, we can minimize the norm to the true covariance $\Sigma_n$ and find the optimal covariance estimator in $\mathcal{C}_n$. This is an "oracle" estimator because the optimal eigenvalues $\tilde{d}_i$ still depends on $\Sigma_n$.
\begin{ppt}[Oracle rotation-invariant covariance estimator, \cite{Ledoit2009} p.9]
	With $\lVert \cdot \rVert$ the Frobenius loss, the minimization problem $\min_{D_n} \left \lVert U_n D_n U_n^* - \Sigma_n \right \rVert$ has the following minimizer:
	\begin{equation*}\label{}
	\begin{aligned}
		 \tilde{D}_n = \Diag(\tilde{d}_1,...,\tilde{d}_n),
	\end{aligned}
	\end{equation*}
	where $\forall i \in \llbracket 1,n \rrbracket, \tilde{d}_i = u_i^* \Sigma_n u_i$. $(\tilde{d}_1,...,\tilde{d}_n)$ are denoted as "the oracle shrunk covariance eigenvalues".
\end{ppt}

The key object to study in order to derive asymptotic shrinkage formulas, introduced originally by Ledoit and Péché \cite{Ledoit2009}, is the following $\Delta_n$. 

The object $\Delta_n$ contains the necessary information to find each $\tilde d_i$, which is what we need to shrink our sample eigenvalues. Using this object has a direct interest: it can be expressed using $\Theta_n^{(1)}$, which asymptotic behavior is characterized by Theorem \ref{FEg}. 

To sum up, through $\Delta_n$, we can describe the asymptotic behavior of the $\tilde d_i$, and thus find asymptotic shrinkage formulas.
\begin{dft}[$\Delta_n$, from \cite{Ledoit2009}]\label{}
We define for all $x \in \R$:
\begin{equation*}
\begin{aligned}
	 &\Delta_n(x) = \frac{1}{n} \sum_{i=1}^n \tilde d_i 1_{[\lambda_i, +\infty[}(x).
\end{aligned}
\end{equation*}
When the sample eigenvalues are all distinct, we can retrieve $\tilde d_i$ through $\Delta_n$:
\begin{equation*}\label{}
\begin{aligned}
	 &\forall i \in \llbracket 1,n\rrbracket, \tilde d_i = \lim_{\varepsilon \rightarrow 0^+} \frac{\Delta_n(\lambda_i + \varepsilon) - \Delta_n(\lambda_i - \varepsilon)}{F_n(\lambda_i + \varepsilon) - F_n(\lambda_i - \varepsilon)}.
\end{aligned}
\end{equation*}
Finally, for all $x \in \R$ such that $\Delta_n$ is continuous at $x$:
\begin{equation*}\label{}
\begin{aligned}
	 \Delta_n(x) = \lim_{\eta \rightarrow 0^+} \frac{1}{\pi} \int_{-\infty}^x \Im[\Theta^{(1)}_n(\xi + i\eta)]d\xi.
\end{aligned}
\end{equation*}
\end{dft}

We can now describe the asymptotic behavior of $\Delta_n$, which naturally leads to asymptotic shrinkage formulas.
\begin{thm}[Oracle covariance shrinkage formula]\label{csh}
Assume the conditions \ref{H1}-\ref{H5}.
There exists a nonrandom function $\Delta: \R \rightarrow \R$ such that a.s, $\forall x \in \R_+, \Delta_n(x) \rightarrow \Delta(x)$.
Moreover, we have by Stieltjes-Perron formula: 
\begin{equation*}\label{}
\begin{aligned}
	 &\Delta(x) = \lim_{\eta \rightarrow 0^+} \frac{1}{\pi} \int_{-\infty}^{x} \Im\left[\Theta^{(1)}(\lambda + i\eta) \right] d\lambda.
\end{aligned}
\end{equation*}
Suppose $c \neq 1$. Then:
	\begin{equation*}\label{}
	\begin{aligned}
		 \Delta(x) = \int_{-\infty}^x h(\lambda) dF(\lambda),
	\end{aligned}
	\end{equation*}
	where for $\lambda \in \R \backslash \{x \in S_F, F'(x) = 0\}$:
	\begin{equation*}\label{}
	\begin{aligned}
		&h(\lambda) = \frac{\int \frac{\tau^2}{|\tau \check X(\lambda) + 1|^2}dH(\tau)}{\int \frac{\tau}{|\tau \check X(\lambda) + 1|^2}dH(\tau)} \text{, if } \lambda \neq 0,\\
		&h(\lambda) = \frac{1}{(c-1) \check X(0)}\text{, if } \lambda = 0 \text{ and } c > 1,\\
		&h(\lambda) = 0 \text{, otherwise.}
	\end{aligned}
	\end{equation*}
\end{thm}

The proof of this result is given in the Appendix \ref{apxD}. Practically, this theorem means that replacing each weighted sample covariance eigenvalue $\lambda_i$ by $h(\lambda)$, provided we know $\check{X}(\lambda)$ and $H$, will asymptotically minimize the Frobenius loss, in the class of rotation-invariant estimators. While Bun \textit{et al.} \cite{Bun2016} mentioned the problem of tractability of their formula for covariance shrinkage in Section III.B of their work, our different approach from the functionals $\Theta^g$ has two main advantages: it makes the formulas computationally tractable (the estimation of $H$ and $\check X$ is addressed in Section \ref{numsection}), and it provides also formulas for precision matrix estimation, detailed in Section \ref{precform}.

\begin{figure}[]
\centering
\includegraphics[width=0.63\linewidth]{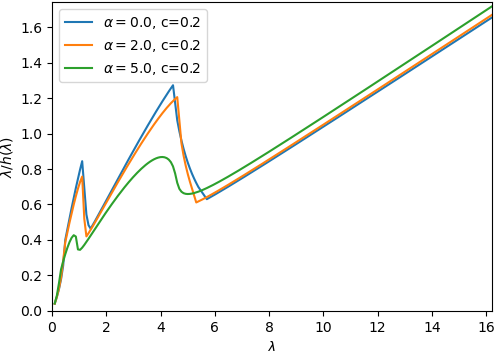}
\caption{Shrinkage intensities $\lambda/h(\lambda)$ in function of $\lambda$, for $c=0.2$, $H = \frac{1}{5} \mathbf{1}_{[1,\infty[} + \frac{2}{5} \mathbf{1}_{[3,\infty[} + \frac{2}{5} \mathbf{1}_{[10,\infty[}$ and $D$ exponentially-weighted with parameter $\alpha \in \{0,2,5\}$ as introduced in Example \ref{exEWMA}.}
\label{fig:s_c4}
\end{figure}

For $\lambda \in \R^*\backslash \{x \in S_F, F'(x) = 0\}$, we denote $\frac{\lambda}{h(\lambda)}$ as the \textit{shrinkage intensity associated to $\lambda$}. An example of shrinkage intensities in function of $\lambda$ for $c=0.2$ is given in figure \ref{fig:s_c4}. The true eigenvalue distribution used in the figure is $H = 0.2 \times 1_{[1,\infty[} + 0.4 \times 1_{[3,\infty[} + 0.4 \times 1_{[10,\infty[}$ as introduced by Bai and Silverstein \cite{Bai1998}. The weight distribution is of the form $D(x) = \alpha \log(x) + \text{cst}$ and mimics an exponentially-weighted moving average (EWMA). The proper definition of this weight distribution is given in Section \ref{ewma}.

We can follow a similar methodology to find oracle asymptotic shrinkage formulas for the precision matrix.

\subsubsection{Asymptotic shrinkage for precision matrix estimation}\label{precform}
For the precision matrix $\Gamma_n = \Sigma_n^{-1}$, we can use a very similar approach and extract an oracle estimator and asymptotic shrinkage formulas for the eigenvalues. We use the same method as with the covariance: 
\begin{itemize}
	\item find the oracle rotation-invariant estimator, 
	\item define an object containing the useful information and that we can study asymptotically, 
	\item and describe its asymptotic behavior.
\end{itemize}

\begin{ppt}[Oracle rotation-invariant precision estimator, \cite{Ledoit2009} p.9]
	With $\lVert \cdot \rVert$ the Frobenius loss, the minimization problem $\min_{D_n} \left \lVert U_n D_n U_n^* - \Gamma_n \right \rVert$ has the following minimizer:
	\begin{equation*}\label{}
	\begin{aligned}
		 \tilde{D}_n = \Diag(\tilde{\gamma}_1,...,\tilde{\gamma}_n),
	\end{aligned}
	\end{equation*}
	where $\forall i \in \llbracket 1,n \rrbracket, \tilde{\gamma}_i = u_i^* \Gamma_n u_i$. $(\tilde{\gamma}_1,...,\tilde{\gamma}_n)$ are denoted as "the oracle shrunk precision eigenvalues".
\end{ppt}

The key object $\Psi_n$ is very similar to $\Delta_n$ in its conception.
\begin{dft}[$\Psi_n$ from \cite{Ledoit2009}]\label{}
We define for all $x \in \R$:
\begin{equation*}
\begin{aligned}
	 &\Psi_n(x) = \frac{1}{n} \sum_{i=1}^n \tilde \gamma_i 1_{[\lambda_i, +\infty[}(x).
\end{aligned}
\end{equation*}
$\Psi_n$ contains the necessary information to retrieve each $\tilde \gamma_i$ when the sample eigenvalues are all distinct:
\begin{equation*}\label{}
\begin{aligned}
	 &\forall i \in \llbracket 1,n\rrbracket, \tilde \gamma_i = \lim_{\varepsilon \rightarrow 0^+} \frac{\Psi_n(\lambda_i + \varepsilon) - \Psi_n(\lambda_i - \varepsilon)}{F_n(\lambda_i + \varepsilon) - F_n(\lambda_i - \varepsilon)}.
\end{aligned}
\end{equation*}
Finally, for all $x \in \R$ such that $\Psi_n$ is continuous at $x$:
\begin{equation*}\label{}
\begin{aligned}
	 \Psi_n(x) = \lim_{\eta \rightarrow 0^+} \frac{1}{\pi} \int_{-\infty}^x \Im[\Theta^{(-1)}_n(\xi + i\eta)]d\xi.
\end{aligned}
\end{equation*}
\end{dft}

We can describe the asymptotic behavior of $\Psi_n$ and deduce asymptotic shrinkage formulas.
\begin{thm}[Oracle precision shrinkage formula]\label{psh}
	Assume the conditions \ref{H1}-\ref{H5}.
	There exists a nonrandom function $\Psi: \R \rightarrow \R$ such that a.s, $\forall x \in \R_+, \Psi_n(x) \rightarrow \Psi(x)$.
	Moreover, we have by Stieltjes-Perron formula:
	\begin{equation*}\label{}
	\begin{aligned}
		\Psi(x) = \lim_{\eta \rightarrow 0^+} \frac{1}{\pi} \int_{-\infty}^{x} \Im\left[\Theta^{(-1)}(\lambda + i\eta) \right] d\lambda.
	\end{aligned}
	\end{equation*}
	Suppose $c \neq 1$. Then:
	\begin{equation*}\label{}
	\begin{aligned}
		&\Psi(x) = \int_{-\infty}^{x} t(\lambda) dF(\lambda) 
	\end{aligned}
	\end{equation*}
	where for $\lambda \in \R \backslash \{x \in S_F, F'(x) = 0\}$:
	\begin{equation*}\label{}
	\begin{aligned}
		& t(\lambda) = \frac{\int \frac{1}{|\tau \check X(\lambda) + 1|^2}dH(\tau)}{\int \frac{\tau}{|\tau \check X(\lambda) + 1|^2}dH(\tau)} \text{, if } \lambda \neq 0,\\
		& t(\lambda) = - \check X(0) + \frac{c}{c-1}\int \frac{1}{\tau}dH(\tau) \text{, if } \lambda = 0 \text{ and } c > 1,\\
		& t(\lambda) = 0 \text{, otherwise.}
	\end{aligned}
	\end{equation*}
\end{thm}
The proof of this result is given in the Appendix \ref{apxD}. Practically, similarly to Theorem \ref{csh}, this theorem means that replacing each weighted sample precision eigenvalue $\frac{1}{\lambda_i}$ by $t(\lambda_i)$, provided you know $H$, will asymptotically minimize the Frobenius loss, in the class of rotation-invariant estimators.

\section{Sample and Population eigenvectors overlap}\label{phisection}
The theoretical study of the sample versus the population eigenvectors overlap can be asymptotically studied thanks to the functionals $\Theta^g$. Indeed, we can study the asymptotic behavior of the following bivariate cumulative distribution function, initially introduced in \cite{Ledoit2009} for unweighted sample covariance. We generalize the asymptotic result here to weighted sample covariances.

\begin{dft}
	For $\lambda, \tau \in \R$ and $n \in \N$, we define:
	\begin{equation}\label{}
	\begin{aligned}
		&\Phi_n(\lambda, \tau) = \frac{1}{n} \sum_{i=1}^n \sum_{j=1}^n |u_i^*v_j|^2 \mathrm{1}_{[\lambda_i,+\infty[}(\lambda) \times \mathrm{1}_{[\tau_j,+\infty[}(\tau).
	\end{aligned}
	\end{equation}
\end{dft}

This object contains key information about the overlap $|u_i^*v_j|^2$. In fact, the average of $|u_i^*v_j|^2$ for $i$ and $j$ such that $\lambda_i \in [\underline \lambda, \overline \lambda]$ and $\tau_j \in [\underline \tau, \overline \tau]$ is equal to:
\begin{equation}\label{}
\begin{aligned}
	&\frac{\sum_{i=1}^n \sum_{j=1}^n n|u_i^*v_j|^2 \mathrm{1}_{[\underline \lambda,\overline \lambda]}(\lambda_i) \times \mathrm{1}_{[\underline \tau,\overline \tau[}(\tau_j)}{\sum_{i=1}^n \sum_{j=1}^n\mathrm{1}_{[\underline \lambda,\overline \lambda]}(\lambda_i) \times \mathrm{1}_{[\underline \tau,\overline \tau[}(\tau_j)} = \frac{\Phi_n(\overline \lambda, \overline \tau) - \Phi_n(\overline \lambda, \underline \tau) - \Phi_n(\underline \lambda, \overline \tau) + \Phi_n(\underline \lambda, \underline \tau)}{\left[F_n(\overline \lambda) - F_n(\underline \lambda)\right] \times \left[H_n(\overline \tau) - H_n(\underline \tau)\right]}
\end{aligned}
\end{equation}
whenever the denominator is positive. The point is explained in detail in Section 1.2 \cite{Ledoit2009}.

We can deduce the asymptotic behavior of this bivariate c.d.f. $\Phi_n$, which depicts non-intuitive phenomenons concerning PCA in high dimension for instance \cite{Ledoit2009}.
\begin{thm}\label{phithm}
	Assume the conditions \ref{H1}-\ref{H5}. Then there exists $\Phi$ such that almost surely, $\Phi_n(\lambda,\tau) \rightarrow \Phi_(\lambda,\tau)$ at all points of continuity of $\Phi$. 

	If $c \neq 1$, we have $\forall (\lambda, \tau) \in \R^2, \Phi(\lambda, \tau) = \int_{-\infty}^\lambda \int_{-\infty}^\tau \phi(\ell ,t)dH(t)dF(\ell )$, where $\forall (\lambda, \tau) \in \R^2$:
	\begin{equation}\label{}
	\begin{aligned}
		& \phi(\ell ,t) = \frac{t}{\int \frac{u|t \check X(\ell ) + 1|^2}{|u \check X(\ell ) + 1|^2}dH(u)} \text{, if } \ell >0,\\
		& \phi(\ell ,t) = \frac{c}{(c-1)(t \check X(0) + 1)} \text{, if } \ell =0 \text{ and } c > 1,\\
		& \phi(\ell ,t) = 0 \text{, otherwise.}
	\end{aligned}
	\end{equation}
\end{thm}

\begin{figure}[]
\centering
\includegraphics[width=0.7\linewidth]{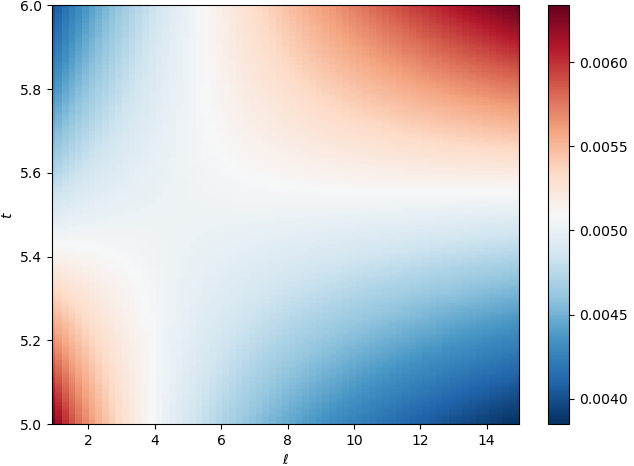}
\caption{Heatmap of the bivariate overlap density $(\ell,t) \mapsto \phi(\ell,t)$ for $H$ uniform in $[5,6]$, concentration ratio $c=0.2$ and $D$ exponentially-weighted with parameter $\alpha=5$.}
\label{fig:phi}
\end{figure}

The proof of this result is given in the Appendix \ref{apxD}. An example of density $(\ell,t) \mapsto \phi(\ell,t)$ is given in Figure \ref{fig:phi}.

A direct application of Theorem \ref{phithm} concerns the choice of threshold $\lambda_{PCA}$ in PCA. We want that the set of eigenvectors $U_{PCA} = \{u_i, \lambda_i \geq \lambda_{PCA}\}$ verifies: $\sum_{u \in U_{PCA}} u^* \Sigma u = \kappa \tr(\Sigma)$, for a ratio of variance $\kappa \in [0,1]$ chosen by the practitioner (often chosen as $\kappa = 0.95$). In other words, we want to find $\lambda_{PCA}$ so that the explained variance of the subspace generated by the eigenvectors $U_{PCA}$ is equal to $\kappa$ times the total variance

Rewriting this key equation for PCA, we have:
\begin{equation}\label{}
\begin{aligned}
	& \frac{1}{n} \sum_{i=1}^n \sum_{j=1}^n \tau_j |u_i^*v_j|^2 \mathrm{1}_{[\lambda_i,+\infty[}(\lambda_{PCA}) = (1-\kappa) \int \tau dH_n(\tau).
\end{aligned}
\end{equation}
Remark that asymptotically, in high dimension or not, choosing the right side to be $ (1-\kappa) \int \tau dH_n(\tau)$ or  $(1-\kappa) \int \ell dF_n(\ell)$, as it is often done in practice, is equivalent. However, the variance explained by $U_{PCA}$ is sensitive to high-dimensional effects. Indeed, asymptotically, from Theorem \ref{phithm}, the PCA equation becomes:
\begin{equation}\label{PCAtrue}
\begin{aligned}
	& \int_{-\infty}^{\lambda_{PCA}} \int_{-\infty}^{+\infty} t\phi(\ell ,t)dH(t)dF(\ell )= (1-\kappa) \int \tau dH(\tau).
\end{aligned}
\end{equation}

Using the PCA defined in a classical setting of low-dimension, \textit{i.e.}, neglecting the effects of high dimensionality, we would find the threshold $\lambda_{naive}$ so that $\sum_{u \in U_{naive}} u^* S u = (1-\kappa) \tr(S)$, with $U_{naive} = \{u_i, \lambda_i \geq \lambda_{naive}\}$. While the right side is asymptotically equivalent to $(1-\kappa) \tr(\Sigma)$ in high dimension, the left one is not. Asymptotically, it would lead to solve:
\begin{equation}\label{PCAnaive}
\begin{aligned}
	& \int_{-\infty}^{\lambda_{naive}} \ell dF(\ell )= \int_{-\infty}^{\lambda_{naive}}\int_{-\infty}^{+\infty} \ell \phi(\ell ,t)dH(t)dF(\ell ) = (1-\kappa) \int \tau dH(\tau).
\end{aligned}
\end{equation}

The core difference between Equation \eqref{PCAtrue} and Equation \eqref{PCAnaive} is the integrand, that goes from $t\phi(\ell ,t)$ to $\ell \phi(\ell ,t)$, \textit{i.e.}, taking in account the true variance of the projected eigenvector $t$, or its in-sample estimate $\ell $.

\begin{figure}[]
\centering
\includegraphics[width=0.48\linewidth]{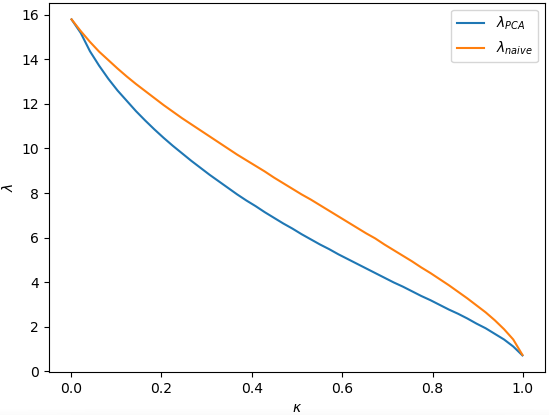}
\includegraphics[width=0.48\linewidth]{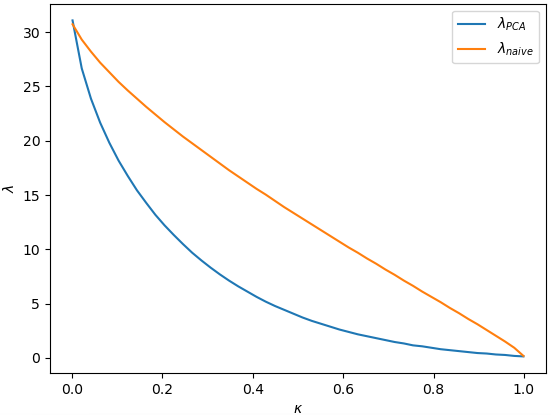}\\
\includegraphics[width=0.48\linewidth]{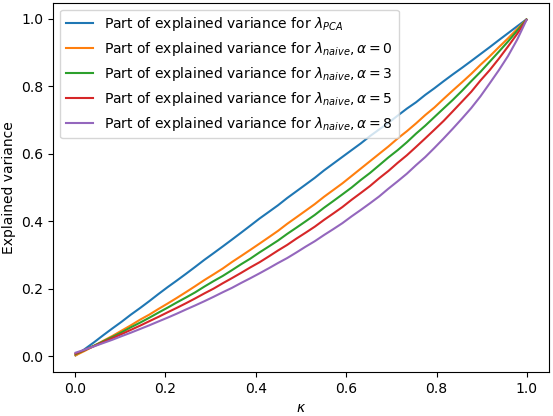}
\includegraphics[width=0.48\linewidth]{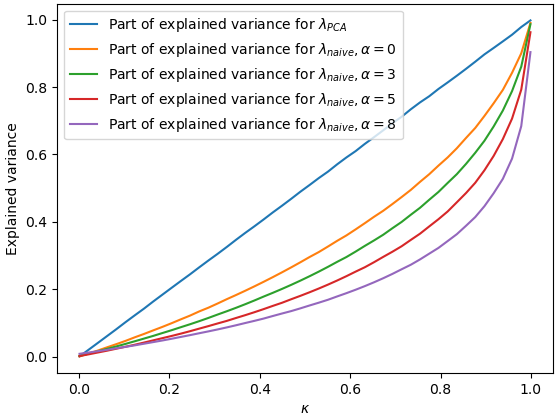}
\caption{(Top) $\lambda_{PCA}$ and $\lambda_{naive}$ in function of $\kappa$, for $H$ uniform in $[1,10]$, concentration ratio $c=0.1$ (left) and $c=0.5$ (right) and $D$ exponentially-weighted with parameter $\alpha=5$. (Down) The proportion of variance explained by the PCA in function of $\kappa$ for $\lambda_{PCA}$ (equal to $\kappa$ by construction) and $\lambda_{naive}$, for $H$ uniform in $[1,10]$, concentration ratio $c=0.1$ (left) and $c=0.5$ (right) and $D$ exponentially-weighted with parameter $\alpha \in \{0,3,5,8\}$.}
\label{fig:PCA}
\end{figure}

In practice, the difference between the corrected method $\lambda_{PCA}$ and the low dimensional one $\lambda_{naive}$ is of utter importance. Some computations are plotted in Figure \ref{fig:PCA} to highlight this phenomenon. As the concentration ratio $c$ increases, the effect of high-dimensionality intensifies. In particular, the part of explained variance that we aim at being $\kappa$ is highly underestimated if we use the low-dimensional approach $\lambda_{naive}$.

\section{Analytical formulas for exponential weight distribution}\label{ewma}
In this section, we want to give an explicit equation for $\check \Theta^{(1)}$ and $\check \Theta^{(-1)}$ in function of $\check m$ with weight distribution that mimics the exponentially weighted moving average (called "exponentially-weighted" or EWMA from this point). We begin by describing mathematically the weight distribution.
\begin{dft}[$\alpha$-exponential weight law]
We fix $\alpha \in \R_+^*$, and we define this law such as its c.d.f. $D$ follows: $D(\beta e^{-\alpha t}) = 1 - t$ for $t \in [0,1]$. Moreover, we impose that $\int \delta dD(\delta) = 1$. We finally have: $\forall x \in [\beta e^{-\alpha}, \beta], D(x) = 1 + \frac{1}{\alpha} \log\left(\frac{x}{\beta}\right)$, with $\beta = \frac{\alpha}{1 - e^{-\alpha}}$.
\end{dft}
We fix the parameter $\alpha \in \R_+^*$ for the remaining of the section. We can now give a more explicit equation to compute $\Theta^{(1)}$ in function of $\check m$.

\begin{ppt}[$\check \Theta^{(1)}$, $\check \Theta^{(-1)}$ and $\check X(\lambda)$ for exponential weight law]\label{exp}
With an $\alpha$-exponential weight law, and $\beta = \frac{\alpha}{1 - e^{-\alpha}}$, we have:
	\begin{equation*}\label{}
	\begin{aligned}
	\forall \lambda \in \R^*, &\check{\Theta}^{(1)}(\lambda) = \frac{1}{\beta c}\frac{e^{\alpha c (1+ \lambda\check{m}(\lambda))} - 1}{1 - e^{-\alpha +\alpha c (1+ \lambda\check{m}(\lambda))}},\\
	\forall \lambda \in \R^*, &\check{\Theta}^{(-1)}(\lambda) = \frac{\check{m}(\lambda)(1+\lambda \check{m}(\lambda))}{\lambda \check{\Theta}^{(1)}(\lambda)} - \frac{1}{\lambda} \int \frac{1}{\tau}dH(\tau),\\
	\forall \lambda \in \R^*, &\check{X}(\lambda) = - \frac{1 + \lambda \check m(\lambda)}{\lambda \check{\Theta}^{(1)}(\lambda)}.
	\end{aligned}
	\end{equation*}
\end{ppt}

The proof of this result is given in the Appendix \ref{apxD}. As for all $\lambda \in (F')^{-1}(\R_+^*)$, $h(\lambda) = \Im[\check{\Theta}^{(1)}(\lambda)]/\Im[\check{m}(\lambda)]$, $t(\lambda) = \Im[\check{\Theta}^{(-1)}(\lambda)]/\Im[\check{m}(\lambda)]$ and $\forall \ell \in \R^*, \phi(\ell ,t) = t/\int \frac{u|t \check X(\ell ) + 1|^2}{|u \check X(\ell ) + 1|^2}dH(u)$, this theorem makes all of this object be computable from $\check m(\lambda)$. In fact, it states that the statistical challenge for exponentially-weighted shrinkage is the same as it is for equally weighted sample covariance: estimate $\check{m}$. 

This estimation can be done through kernel estimation in the unweighted setting, for example \cite{Jing2010,Ledoit2020b}, but the literature is still lacking theoretical support in the weighted setting to the best of our knowledge. Remark that $\frac{1}{\lambda} \int \frac{1}{\tau}dH(\tau)$ is real and thus does not appear in the computation of $t(\lambda)$.

Another way of computing those objects all together, $\check m$, $\check \Theta^g$, and $\phi$, can be done through the estimation of $H$ directly and use of Theorem \ref{c0} then. We detail this approach in the next section.

\section{Numerical implementation}\label{numsection}
We break the computation of the asymptotic non-linear shrinkage formulas and $\phi$ from the samples $X$ into two parts:
\begin{itemize}
	\item retrieving $H$ from the samples and the weights,
	\item computing $\check X(\lambda_i)$ for each sample eigenvalue $\lambda_i$.
\end{itemize} 

\subsection{$H$ retrieval from the sample spectrum}
We propose an application of the universal convergence of the sample spectrum to retrieve the population spectrum distribution $H$ as a mixture of Dirac distributions from observed weighted sample eigenvalues $(\lambda_i)_{i=1}^n$.

The procedure is the following:
\begin{itemize}
	\item[1-] As input, we take the observed sample spectrum distribution $F_n = \frac{1}{n} \sum_{i=1}^n 1_{[\lambda_i,\infty[}$ and the weight matrix $W$.
	\item[2-] Find the estimated true spectrum $\hat H(\tau) = \frac{1}{n}\sum_{i=1}^n 1_{[\tau_i,+\infty[}$ where $\tau$ solves:
		\begin{equation}\label{autoquest}
			\begin{aligned}
				\min_{\tau \in \R^n}\E_{Z}\left[\left\lVert \tilde F_n(Z) - F_n \right\rVert_{\mathcal{W},2}^2 \right]
			\end{aligned}
			\end{equation} 
	where $\lVert \cdot \rVert_{\mathcal{W},2}$ is the $2$-Wasserstein norm and $\tilde F_n(Z) := F^{\frac{1}{N}\sqrt{T}ZWZ^*\sqrt{T}}$ with:
	\begin{itemize}
		\item $T = \Diag\left((\tau_i)_{i=1}^n\right)$, 
		\item $Z$ of size $(n,N)$ with i.i.d. $Z_{ij} \sim \mathcal{N}(0,1)$.
	\end{itemize}
	We use automatic differentiation to solve it.
\end{itemize}

\begin{rmk}[Noise sampling]
	We note that the expectation in \eqref{autoquest} can be computed under any centered and standardized distribution, irrespectively of the noise of the observed phenomenon due to the universality of the asymptotic spectrum.
\end{rmk}

\begin{figure}[]
\centering
\includegraphics[width=0.57\linewidth]{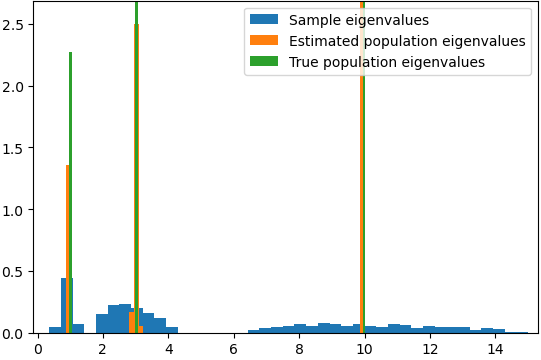}
\includegraphics[width=0.57\linewidth]{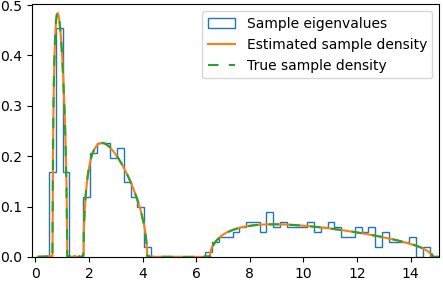}
\caption{(Top) Histograms of sample eigenvalues, estimated population eigenvalues $\hat H$, and true population eigenvalues $H$. (Down) Estimated and true sample density computed on $S_F$ and sample eigenvalues' histogram.}
\label{fig:autoquest}
\end{figure}

An experimental result is shown in Figure \ref{fig:autoquest} to illustrate how the algorithm works. We used $H = \frac{1}{5} \mathbf{1}_{[1,\infty[} + \frac{2}{5} \mathbf{1}_{[3,\infty[} + \frac{2}{5} \mathbf{1}_{[10,\infty[}$, $D$ exponentially weighted with $\alpha = 1$, $c=0.1$, and $Z_{ij} \sim  t_{12}(0,1)$. More experiments are available in the Appendix \ref{apxA}.

\begin{rmk}[Complexity in high dimension]
Regarding the complexity of the algorithm, the main point is the eigenvalue decomposition of the sample covariance at each step of the minimization, which is in $O(n^3)$. When the dimension $n$ too high, and the spectral decomposition is not affordable at each step, we can draw noise matrices $Z$ of size $(n',N')$ instead of $(n,N)$, so that $n'/N' \approx n/N$. We choose $n' < n$ so that the spectrum computation becomes affordable, for example $n' = 200$. Doing that, we sample a spectrum $\tilde F_n$ made of a mixture of $n'$ Dirac distributions, that we can still compare to $F_n$ using Wasserstein distance: it does not require both spectrums to have the same number of Dirac distributions. Obviously, this subsampling of the spectrum will cost accuracy, and a complexity/precision trade-off appears.
\end{rmk}

\subsection{Numerical computation of $\check X$, $\check m$ and $\Im[\check \Theta^g]$}
In this section, we adress the practical problem of numerically computing $\check X$ at a point $\lambda \in \R^*$ while knowing $H$, $D$ and $c$. We saw in Theorem \ref{c0} that computing $\check X(\lambda)$ is enough to compute $\check m(\lambda)$, $\check \Theta^g(\lambda)$ and $\phi(\lambda, \tau)$. 

We define, for $z \in \C$:
\begin{equation*}
\begin{aligned}
	f_z: X \mapsto X + \int \frac{\delta}{z - \delta c \int \frac{\tau}{\tau X + 1}dH(\tau)}dD(\delta)
\end{aligned}
\end{equation*}
We assume that $f_z$ is easy to evaluate at any $X \in \C$ in its domain of definition. This is the practical case where $H$ and $D$ are finite mixtures of Dirac distributions for example, or when they can be efficiently sampled for Monte-Carlo evaluation.

In this situation, $X(z)$ for $z \in \C_+$ is the unique solution to $f_z(X) = 0$ for $X \in \C \backslash \C_-$ and can be solved as a classical minimization problem. In our experiments, a first-order minimization algorithm worked efficiently for this task.

However, computing the limit $\check X(\lambda) = \lim_{z \in \C_+ \rightarrow \lambda} X(z)$ for $\lambda \in \R^*$ is \textit{a priori} a more difficult task. Fortunately, we prove that $\check X(\lambda)$ is a solution to the equation $f_{\lambda}(X) = 0$, which is a first important step. 

But, we cannot apply the same method as with $z \in \C_+$ and directly solve $f_{\lambda}(X) = 0$ for $X \in \C \backslash \C_-$ because there can be many solutions to this equation. Indeed, considering $D = 1_{[1,+\infty]}$ and $H$ a mixture of $N \in \N^*$ Dirac distributions, $f_{\lambda}(X) = 0$ as up to $N$ different solutions in $\C \backslash \C_-$.

We can split the problem into two different scenarios:
\begin{itemize}
	\item if $\check X(\lambda) \in \C_+$, then $\check X(\lambda)$ is the \textbf{unique} solution in $\C_+$ (and not $\C \backslash \C_-$ this time) of the equation $f_{\lambda}(X) = 0$,
	\item otherwise, \textit{i.e.} $\check X(\lambda) \in \R$, the equation (on $X$) $f_{\lambda}(X) = 0$ has only \textbf{real} solutions and $\check X(\lambda)$ is one of them.
\end{itemize}

This is formally stated in the following theorem.
\begin{thm}[$\check X$ computing]\label{checkX}
	Assume the conditions \ref{H1}-\ref{H5}. Then, for $\lambda \in \R^*$, we denote $f_\lambda: X \in \C_+ \cup R \mapsto X + \int \frac{\delta}{\lambda - \delta c \int \frac{\tau}{\tau X + 1}dH(\tau)}dD(\delta)$, where $R$ is its domain of definition included in $\R$. Then, we have that $f_\lambda(\check X(\lambda)) =0$. Moreover,
	\begin{equation*}
	\begin{aligned}
		\check X(\lambda) \in \C_+ &\iff f_\lambda(X) = 0 \text{ has at least a solution in } X \in \C_+, \\
		&\iff \check X(\lambda)  \text{ is the unique solution to } f_\lambda(X) = 0 \text{ in } \C_+.
	\end{aligned}	
	\end{equation*}
\end{thm}

The proof of this result is given in the Appendix \ref{apxE}. Based on this, we suggest the following procedure to compute $\check X(\lambda)$.
\begin{ppt}[Procedure for computing $\check X(\lambda)$]\label{checkXprocedure}
	Let $\lambda \in \R^*$. To compute $\check X(\lambda)$:
	\begin{itemize}
		\item try to solve $f_{\lambda}(X) = 0$ on $\C_+$ (or $\C_{+\varepsilon} := \{z \in \C, \Im[z] \geq \varepsilon)\}$ for some $\varepsilon > 0$), if it succeeds with a solution $X_0 \in \C_+$ (or in $\C_{+\varepsilon}$ for numerical reasons) then $\check X(\lambda) = X_0$,
		\item otherwise, $\check X(\lambda) \in \R$, solve $f_{\lambda + i\eta}(X)$ in $\C_+$, for some $\eta > 0$ small, there is a unique solution $X(\lambda + i \eta)$ and use $\check X(\lambda) = \Re[X(\lambda + i \eta)]$.
	\end{itemize}
\end{ppt}

Moreover, for $n$ sufficiently large, the probability that we find $\lambda_i$ outside the asymptotic support $S_F$ is zero, under the conditions detailed in \cite{Paul2009}. All in all, we do not expect to find any $\check X(\lambda_i) \in \R$ for $n$ sufficiently large. So $f_{\lambda_i}(X) = 0$ will have a unique solution in $\C_+$ with probability $1$ in this setting, which makes the computation $\check X(\lambda_i)$ more robust to numerical approximation of low imaginary parts.

The result of this section is an easy-to-implement algorithm that makes possible to compute the shrinkage coefficients $h(\lambda_i)$ and $t(\lambda_i)$ for each sample eigenvalue $\lambda_i$, and thus define a \textit{bona fide} covariance and precision matrix estimator from the samples $X$.

\section{Experimental results}

\subsection{Performance of the asymptotic oracle shrinkage formulas}
To show the performance of the oracle asymptotic shrinkage formulas in finite sample, we consider the experiment used by Ledoit and Péché \cite{Ledoit2009}. We are looking at the PRIAL of the rotation-invariant covariance estimator using the corrected eigenvalues with the asymptotic oracle formulas of Theorem \ref{csh}. We approximate the shrinkage intensities by numerically solving equations from the asymptotic spectrum on a grid made of the sample eigenvalues $(\lambda_i)_i$. The reference distribution of true eigenvalues is $H = \frac{1}{5} \mathbf{1}_{[1,\infty[} + \frac{2}{5} \mathbf{1}_{[3,\infty[} + \frac{2}{5} \mathbf{1}_{[10,\infty[}$, introduced by Bai and Silverstein \cite{Bai1998} and used by Ledoit and Péché \cite{Ledoit2009}. 

We remind that the PRIAL of an estimator $\hat \Sigma$ of the true covariance $\Sigma$ compared to the sample covariance $S$ is defined as:
\begin{equation*}\label{}
\begin{aligned}
\text{PRIAL}(\hat \Sigma) = \frac{\E[\lVert S - \Sigma \rVert^2] - \E[\lVert \hat \Sigma - \Sigma \rVert^2]}{\E[\lVert S - \Sigma \rVert^2]}
\end{aligned}
\end{equation*}

\begin{figure}[]
\centering
\includegraphics[width=0.63\linewidth]{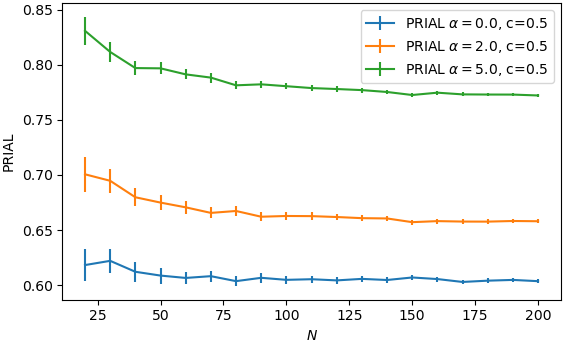}
\caption{PRIAL of the oracle weighted shrunk estimator in function of the number of samples $N$, for $c=0.5$, with $H = \frac{1}{5} \mathbf{1}_{[1,\infty[} + \frac{2}{5} \mathbf{1}_{[3,\infty[} + \frac{2}{5} \mathbf{1}_{[10,\infty[}$ and exponentially-weighted distribution of parameter $\alpha$, Monte-Carlo of $n_{MC} = 50$ draws for each point.}
\label{fig:prial_exp}
\end{figure}

We approximate the expectations by Monte-Carlo. The result of the experiment is shown in Figure \ref{fig:prial_exp}

\subsection{Convergence with heavy tails}
This last experimental part aims at empirically relaxing the strong hypothesis we made in the theoretical part: we assumed bounded $12^{th}$ moments for the underlying distribution of $X$. Indeed, our experiments lead to think that only bounded $4^{th}$ moments are really necessary to see the convergence to the asymptotic density we described in Theorem \ref{FEg}. Even for heavier tails, the estimator seems to still show strong performance.

To visualize this point, we are plotting the histogram of $\Delta_n$ along with the asymptotic theoretical density $\frac{1}{\pi} \Im[\check \Theta^{(1)}(\cdot)]$. 
We draw the $(X_{ij})$ from a normalized t-distribution with varying degree of freedom $\nu \in \{12,3,2.5\}$. At $\nu = 12$, we are at the edge of the theoretical requirement, but we see visible issues with the convergence only with $\nu < 3$, as a non-negligible part of values are outside the theoretical support.

\begin{figure}[]
\centering
\includegraphics[width=0.63\linewidth]{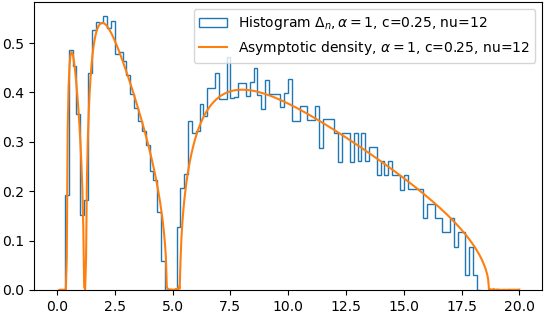}
\includegraphics[width=0.63\linewidth]{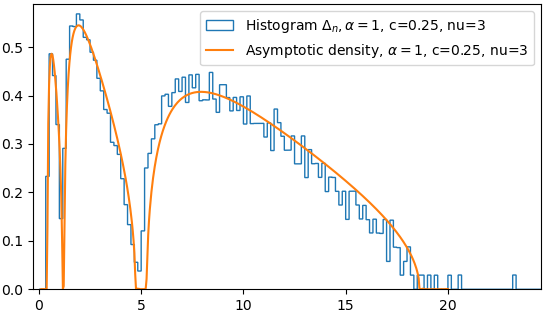}
\includegraphics[width=0.63\linewidth]{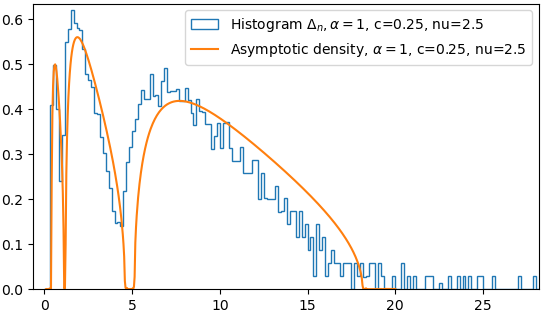}
\caption{Histogram of $\Delta_n$ and its asymptotic theoretical density with underlying noise using t-distribution with $\nu \in \{12, 3, 2.5\}$ top to down, for $n = 2000$, $c=0.25$, $H = \frac{1}{5} \mathbf{1}_{[1,\infty[} + \frac{2}{5} \mathbf{1}_{[3,\infty[} + \frac{2}{5} \mathbf{1}_{[10,\infty[}$ and exponentially-weighted distribution of parameter $\alpha = 1$.}
\label{fig:hist_exp}
\end{figure}

As previously, we chose $H = \frac{1}{5} \mathbf{1}_{[1,\infty[} + \frac{2}{5} \mathbf{1}_{[3,\infty[} + \frac{2}{5} \mathbf{1}_{[10,\infty[}$ and an exponentially-weighted distribution of parameter $\alpha = 1$ for the weights. The concentration ratio is set at $c=0.25$ and the dimension at $n = 2000$. The convergence is almost sure, so Monte-Carlo is not necessary: we only drew one weighted sample covariance per plot. The results are shown in Figure \ref{fig:hist_exp}.

\begin{figure}[]
\centering
\includegraphics[width=0.63\linewidth]{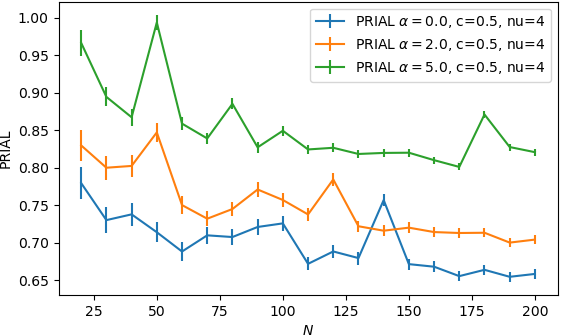}
\includegraphics[width=0.63\linewidth]{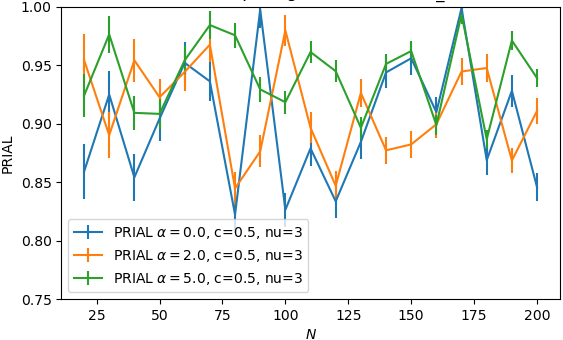}
\includegraphics[width=0.63\linewidth]{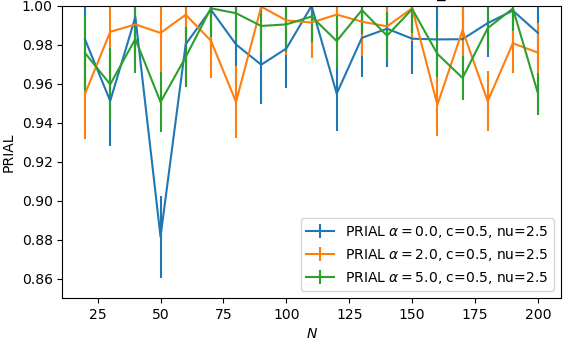}
\caption{PRIAL of the oracle weighted shrunk estimator from a t-distribution with $\nu \in \{ 4, 3, 2.5 \}$ from top to bottom, in function of $N$, for $c=0.5$, with $H = \frac{1}{5} \mathbf{1}_{[1,\infty[} + \frac{2}{5} \mathbf{1}_{[3,\infty[} + \frac{2}{5} \mathbf{1}_{[10,\infty[}$ and $D$ the exponentially-weighted distribution of parameter $\alpha \in \{0,2,5\}$, Monte-Carlo of $n_{MC} = 400$ draws.}
\label{fig:prial_nu}
\end{figure}

Regarding the PRIAL of the asymptotic oracle estimator, the experiment shows a higher robustness of the method, as even with very heavy-tails, such as $\nu = 3$, no particular artifact is appearing in our setting. We used the same parameters as in the previous experimental part concerning the PRIAL, just replacing the underlying Gaussian distribution with a t-distribution of various degrees of freedom $\nu \in \{ 4, 3, 2.5 \}$. The result is shown in figure \ref{fig:prial_nu}, and the Gaussian setting is in figure \ref{fig:prial_exp}.

To conclude this experimental part, according to the previous experiments, the assumption requiring a bounded $12^{th}$ moments of the underlying distribution seems more technical theoretical requirement rather than a real limitation. For $\Theta^g$ convergence, bounded $4^{th}$ moments seem sufficient, as it was noted by Ledoit and Wolf \cite{Ledoit2012} in the equally-weighted setting. 

Regarding the performance of the asymptotic estimator, the convergence of the PRIAL is observed at least until bounded $4^{th}$ moments. For heavier tails, the PRIAL of the non-linear shrinkage estimator has more variance but takes higher values than expected: the non-linear shrinkage estimator is less impacted by heavy tails than the weighted sample covariance.

\section{Conclusion}
In this work, we compute asymptotic optimal formulas for non-linear shrinkage of weighted sample covariance in the class of rotation-invariant covariance and precision matrix estimators. 

Experimentally, we demonstrate the performance in PRIAL of the optimal shrinkage formulas. Finally, we confirm the robustness of the theory when the underlying distribution has only bounded $3^{rd}$ moments instead of $12$ required in the theoretical part.

The main limitation of the use of these formulas concerns their implementation in high dimension. When the dimension and the number of samples are too high, a thousand each for example, to correctly use the algorithm we propose in reasonable time, we currently need to undersample the spectrum at each epoch and lose precision. This point needs further work or extensions. As it was done recently with different approaches \cite{Ledoit2019}, \cite{Ledoit2020b}, analytical estimators of the oracle shrinkage intensities could be derived using kernel density estimation of the weighted empirical spectrum for example. Moreover, the structure of the asymptotic support could be used too to design a weighted version of the QuEST algorithm \cite{Ledoit2016} that can handle very high dimension.

\section{Appendix: Additional experiments on covariance and precision shrinkage}\label{apxA}
In this section, we illustrate how non-linear shrinkage formulas affect the spectrum of the weighted sample covariance and precision matrices. 

\begin{figure}[]
\centering
\includegraphics[width=0.7\linewidth]{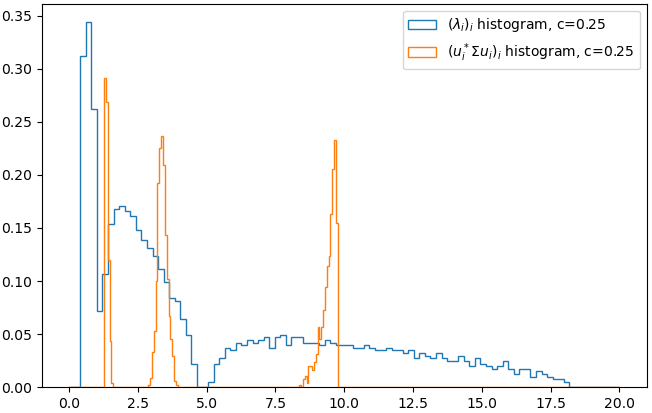}
\includegraphics[width=0.7\linewidth]{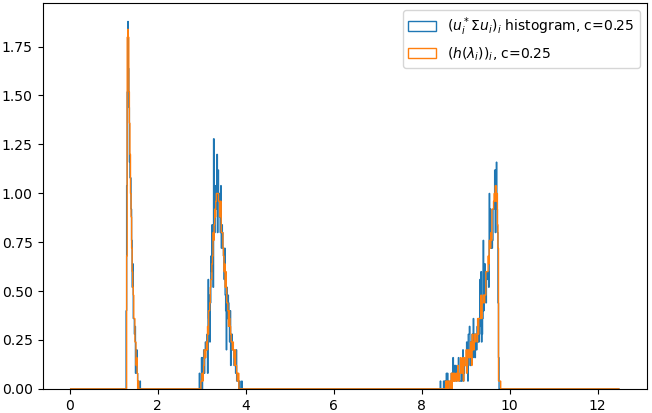}
\caption{Histograms of the weighted sample eigenvalues $\lambda_i$, the optimal shrinkage $u_i^*\Sigma u_i$ (rescaled in top figure for readability) and the asymptotic formula $h(\lambda_i)$, weight distribution $D$ $\alpha$-exponential with $\alpha=1$, $H = \frac{1}{5} \mathbf{1}_{[1,\infty[} + \frac{2}{5} \mathbf{1}_{[3,\infty[} + \frac{2}{5} \mathbf{1}_{[10,\infty[}$, $Z_{ij}$ Gaussian, $c=0.25$, $n=2000$.}
\label{fig:cov_exp1}
\end{figure}

\begin{figure}[]
\centering
\includegraphics[width=0.7\linewidth]{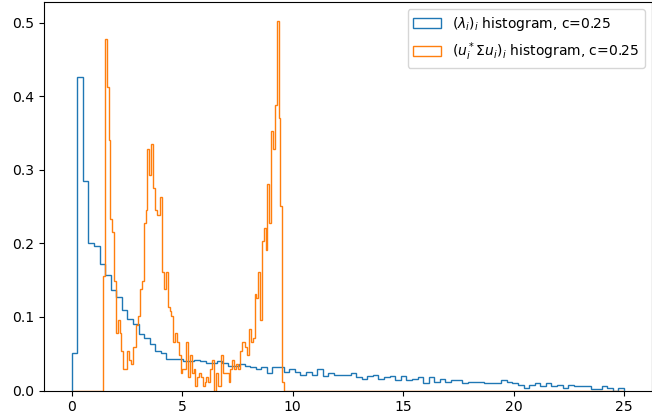}
\includegraphics[width=0.7\linewidth]{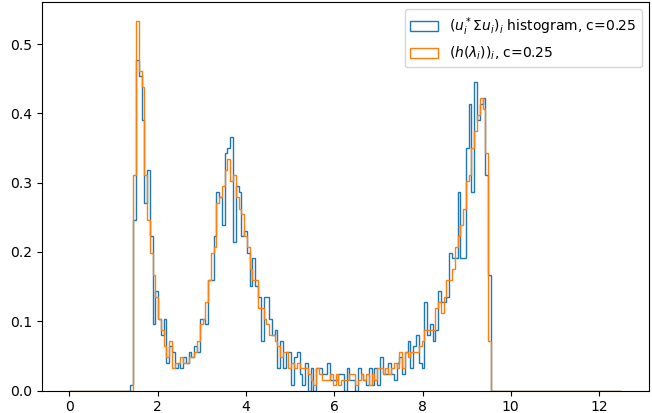}
\caption{Histograms of the weighted sample eigenvalues $\lambda_i$, the optimal shrinkage $u_i^*\Sigma u_i$ and the asymptotic formula $h(\lambda_i)$, weight distribution $D$ $\alpha$-exponential with $\alpha=5$, $H = \frac{1}{5} \mathbf{1}_{[1,\infty[} + \frac{2}{5} \mathbf{1}_{[3,\infty[} + \frac{2}{5} \mathbf{1}_{[10,\infty[}$, $Z_{ij}$ Gaussian, $c=0.25$, $n=2000$.}
\label{fig:cov_exp5}
\end{figure}

\begin{figure}[]
\centering
\includegraphics[width=0.7\linewidth]{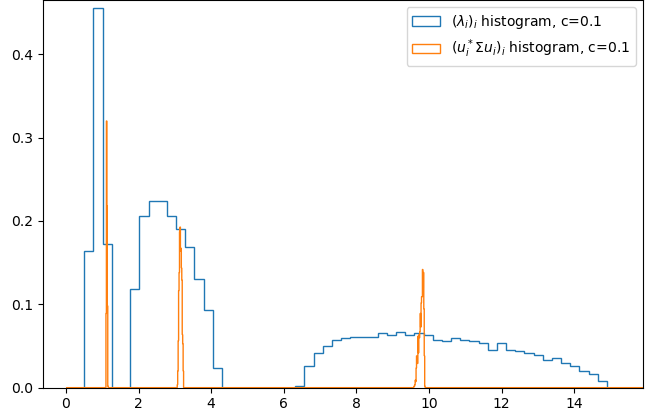}
\includegraphics[width=0.7\linewidth]{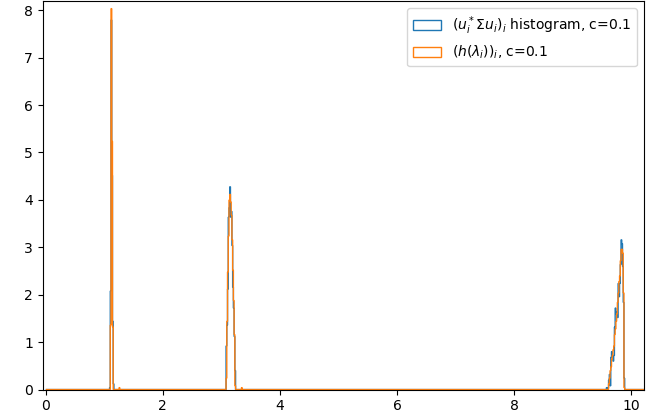}
\caption{Histograms of the weighted sample eigenvalues $\lambda_i$, the optimal shrinkage $u_i^*\Sigma u_i$ (rescaled in top figure for readability) and the asymptotic formula $h(\lambda_i)$, weight distribution $D$ $\alpha$-exponential with $\alpha=1$, $H = \frac{1}{5} \mathbf{1}_{[1,\infty[} + \frac{2}{5} \mathbf{1}_{[3,\infty[} + \frac{2}{5} \mathbf{1}_{[10,\infty[}$, $Z_{ij}$ Gaussian, $c=0.1$, $n=2000$.}
\label{fig:cov_exp1_c10}
\end{figure}

We investigate firstly the case where the weight distribution follows an $\alpha$-exponential weight law, for different $\alpha$ and $c$, for the same distribution $H$ we used for illustration in the corpus: $H = \frac{1}{5} \mathbf{1}_{[1,\infty[} + \frac{2}{5} \mathbf{1}_{[3,\infty[} + \frac{2}{5} \mathbf{1}_{[10,\infty[}$, $Z_{ij}$. In Figures \ref{fig:cov_exp1}, \ref{fig:cov_exp5} and \ref{fig:cov_exp1_c10}, we show the histograms of the weighted sample eigenvalues $(\lambda_i)$, the oracle non-linear shrinkage $u_i^*\Sigma u_i$ and the asymptotic non-linear shrinkage we derived $(h(\lambda_i))$. Through this experiment, we empirically assess the convergence of $\Delta_n(x)$ to $\Delta(x)$ for $x \in \R^*$.

\begin{figure}[]
\centering
\includegraphics[width=0.7\linewidth]{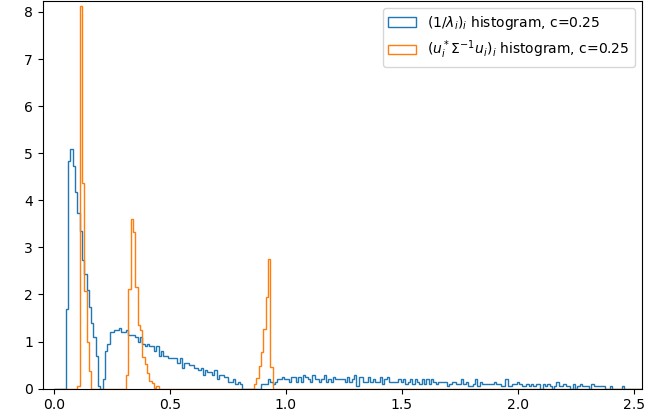}
\includegraphics[width=0.7\linewidth]{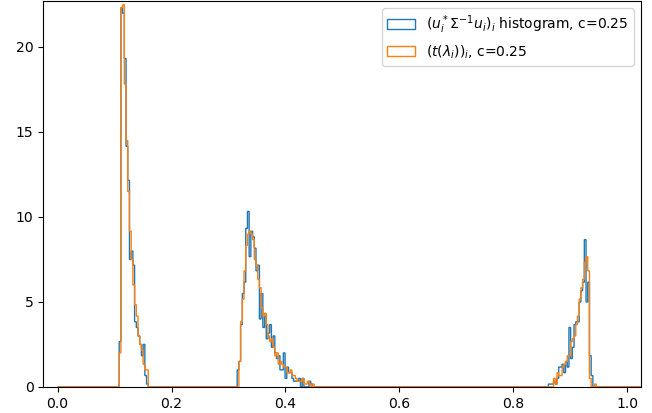}
\caption{Histograms of the weighted sample precision eigenvalues $1/\lambda_i$, the optimal shrinkage $u_i^*\Sigma^{-1} u_i$ (rescaled in top figure for readability) and the asymptotic formula $t(\lambda_i)$, weight distribution $D$ $\alpha$-exponential with $\alpha=1$, $H = \frac{1}{5} \mathbf{1}_{[1,\infty[} + \frac{2}{5} \mathbf{1}_{[3,\infty[} + \frac{2}{5} \mathbf{1}_{[10,\infty[}$, $Z_{ij}$ Gaussian, $c=0.25$, $n=2000$.}
\label{fig:prec_exp1}
\end{figure}

The three last figures give the idea of what does non-linear shrinkage on the weighted sample spectrum estimation. We give one more figure, Figure \ref{fig:prec_exp1} to show an example of \textbf{precision} matrix non-linear shrinkage. Through this experiment, we empirically assess the convergence of $\Psi_n(x)$ to $\Psi(x)$ for $x \in \R^*$. The behavior is very similar to the covariance shrinkage, so we will not dive more into details.

\begin{figure}[]
\centering
\includegraphics[width=0.7\linewidth]{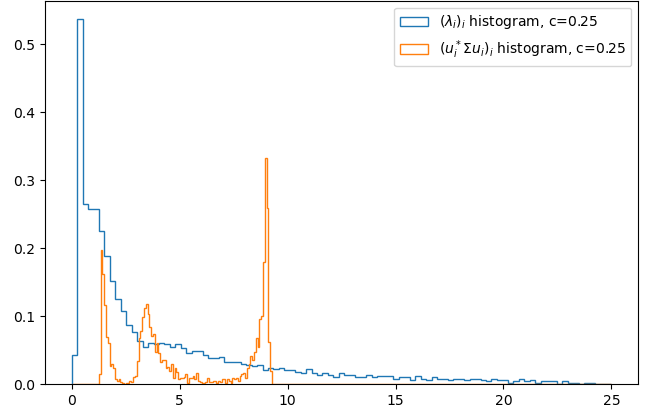}
\includegraphics[width=0.7\linewidth]{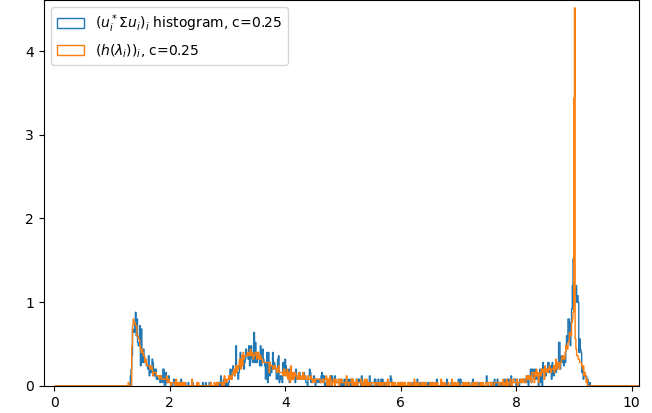}
\caption{Histograms of the weighted sample covariance eigenvalues $\lambda_i$, the optimal shrinkage $u_i^*\Sigma u_i$ (rescaled in top figure for readability) and the asymptotic formula $h(\lambda_i)$, weight distribution $D$ mixture of $5$ Dirac distributions, $H = \frac{1}{5} \mathbf{1}_{[1,\infty[} + \frac{2}{5} \mathbf{1}_{[3,\infty[} + \frac{2}{5} \mathbf{1}_{[10,\infty[}$, $Z_{ij}$ Gaussian, $c=0.25$, $n=2000$.}
\label{fig:cov_5dir}
\end{figure}

\begin{figure}[]
\centering
\includegraphics[width=0.7\linewidth]{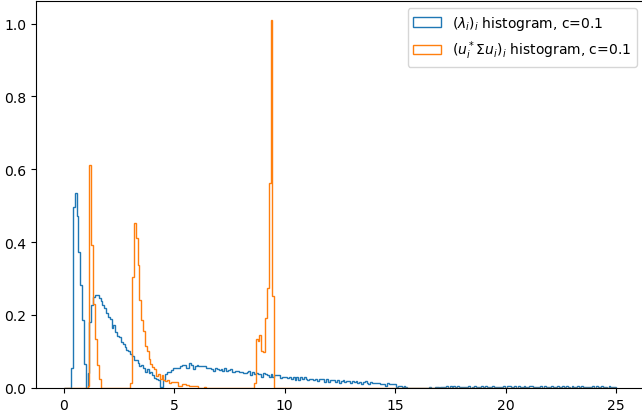}
\includegraphics[width=0.7\linewidth]{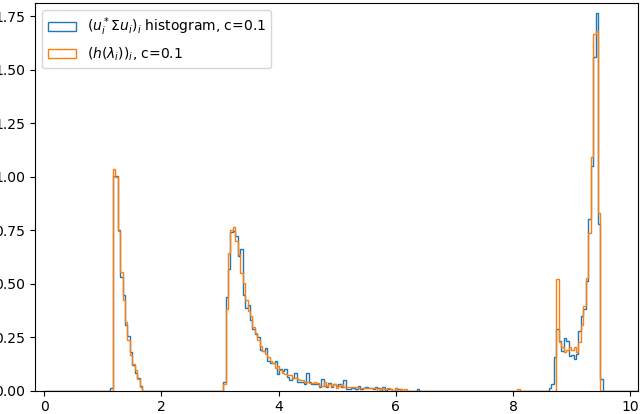}
\caption{Histograms of the weighted sample covariance eigenvalues $\lambda_i$, the optimal shrinkage $u_i^*\Sigma u_i$ (rescaled in top figure for readability) and the asymptotic formula $h(\lambda_i)$, weight distribution $D$ mixture of $5$ Dirac distributions, $H = \frac{1}{5} \mathbf{1}_{[1,\infty[} + \frac{2}{5} \mathbf{1}_{[3,\infty[} + \frac{2}{5} \mathbf{1}_{[10,\infty[}$, $Z_{ij}$ Gaussian, $c=0.1$, $n=2000$.}
\label{fig:cov_5dir_c10}
\end{figure}

We finish this section with a more exotic weight distributions: a mixture of $5$ Dirac distributions, and create an additional gap in the sample spectrum for low $c$: $$D = 0.593 \times \mathbf{1}_{[0.337,+\infty[} +0.297 \times \mathbf{1}_{[0.674,+\infty[} + 0.074 \times \mathbf{1}_{[2.696,+\infty[} + 0.030 \times \mathbf{1}_{[6.74,+\infty[} + 0.006 \times \mathbf{1}_{[33.7,+\infty[}.$$ Results are shown in Figures \ref{fig:cov_5dir} and \ref{fig:cov_5dir_c10} for two different $c$.

\section{Appendix: Proof of Theorem \ref{FEg}}\label{apxB}
For this proof of Theorem \ref{FEg}, we adapt the proof of Ledoit and Péché for Theorem 1.2 \cite{Ledoit2009}. We keep the core structure of the proof, iteratively proving the theorem for more and more complex functions $g$, but we adapt it to our more general setting. In particular, the convergence proofs are different and the technical lemmas too. The idea of the proof mostly relies on the following lemma.
\begin{lem}[Lemma 2.1 \cite{Ledoit2009}]\label{conc}\label{lemma_rmt}
	Let $Y=(y_1,...,y_n)$ be a random vector with independent entries satisfying:
	\begin{equation}\label{}
	\begin{aligned}
			&\E[y_i] = 0, \E[|y_i|^2] = 1, \E[|y_i|^{12}] \leq B,
	\end{aligned}
	\end{equation}
	where the constant $B$ does not depend on $n$. Let also $A$ be a given $n \times n$ matrix. Then there exists a constant $K>0$ independent of $n$, $A$ and $Y$ such that:
	\begin{equation}\label{}
	\begin{aligned}
			&\E\left[\left| YAY^* - \Tr(A) \right|^6\right] \leq K \lVert A \rVert^6n^3.
	\end{aligned}
	\end{equation}
\end{lem}

\begin{rmk}
	Remark that in the original formulation of Lemma 2.1 \cite{Ledoit2009}, which adapts Lemma 3.1 \cite{Silverstein1995b}, the assumption about $Y$ is that the components are i.i.d, however in the proof the "identically distributed" assumption is in fact not used.
\end{rmk}

In order to prove Theorem \ref{FEg}, we will firstly prove it for $g = 1$, for $g = Id$, then for $g$ polynomial by induction, and for $g$ continuous by density, and finally for $g$ bounded with a finite number of discontinuities also by induction.

\subsection{Case $g=1$}
The first part of the proof focuses on $m_n(\cdot)$ and its pointwise a.s. convergence to $m(\cdot)$ in $\C_+$. The result is already known under weaker assumptions - in \cite{Zhang2007} for the weaker set -, but we propose here a proof and its technical lemmas that will be the backbone of the complete proof of Theorem \ref{FEg}. Thus, the development of this proof is of core importance in the whole process, and most of the technical tools and work is introduced, and will make the end of the proof straightforward.

\begin{lem}\label{g1}
	Assume \ref{H1}-\ref{H5}. For all $z \in \C_+$, $m_n(z) := \Theta^{(0)}_n(z) \rightarrow \Theta^{(0)}(z)$ almost surely and:
	\begin{equation*}\label{}
	\begin{aligned}
		 &m(z) := \Theta^{(0)}(z) = -\frac{1}{z}\int \frac{1}{\tau X(z) + 1}dH(\tau),
	\end{aligned}
	\end{equation*}
	where for all $z \in \C_+$, $X(z)$ is the unique solution in $\C_+$ of the following equation:
	\begin{equation*}\label{}
	\begin{aligned}
		X(z) &=-\int \frac{\delta}{z -  \delta c \int \frac{\tau}{\tau X(z) + 1}dH(\tau)}dD(\delta).
	\end{aligned}
	\end{equation*}
\end{lem}

We remind that we denote $B_n = \frac{1}{N} T^{1/2} Z^* W Z T^{1/2}$ and $\underline{B}_n = \frac{1}{N} W^{1/2} Z T Z^* W^{1/2}$. We introduce new objects of interest in the analysis, namely $q_j, r_j$, and $B_{(j)}$.
\begin{notat}
	For $j \in \llbracket 1,N \rrbracket$, we denote:
	\begin{itemize}
		\item $q_j = \frac{1}{\sqrt{n}}Z_{\cdot j}$, 
		\item $r_j =  \frac{1}{\sqrt{N}} T^{1/2}Z_{\cdot j} W_{jj}^{1/2}$, 
		\item $B_{(j)} = B_n - r_j r_j^*$, 
	\end{itemize}
\end{notat}

\subsubsection{Preliminary derivations}
Let's introduce a technical result.
\begin{lem}[Eq (2.1) \cite{Silverstein1995c}]\label{tech0}
	For $B$ a $(n,n)$ matrix, $q \in  \C^n$ for which $B$ and $B + qq^*$ is invertible, we have:
	\begin{equation}\label{}
	\begin{aligned}
		q^* (B + qq^*)^{-1} = \frac{1}{1 + q^*B^{-1}q}q^* B^{-1}.
	\end{aligned}
	\end{equation}
\end{lem}

We have: 
\begin{equation}\label{}
\begin{aligned}
	(B_n -zI) +zI &= \sum_{j=1}^N r_j r_j^*  \\
	I + z(B_n - zI)^{-1} &= \sum_{j=1}^N r_j r_j^* (B_n -zI)^{-1} \\
	I + z(B_n - zI)^{-1} &= \sum_{j=1}^N \frac{1}{1 + r_j^* (B_{(j)} -zI)^{-1}r_j} r_j r_j^* (B_{(j)} -zI)^{-1}  \\
	c_n + z c_n m_n(z) &= 1 - \frac{1}{N} \sum_{j=1}^N \frac{1}{1 + r_j^* (B_{(j)} -zI)^{-1}r_j} \\
	\underline{m}_n(z) &= - \frac{1}{N} \sum_{j=1}^N \frac{1}{z\left(1 + r_j^* (B_{(j)} -zI)^{-1}r_j\right)}.
\end{aligned}
\end{equation} 

Here is the crucial difference between the proof in \cite{Silverstein1995a} in the evenly weighted case and the weighted case. We denote:
\begin{equation}\label{} 
\begin{aligned}
	\alpha_n = - \frac{1}{z\underline{m}_n(z)} \frac{1}{N} \sum_{j=1}^N  W_{jj}\left(r_j^* (B_n -zI)^{-1}r_j - 1 \right),
\end{aligned}
\end{equation} 
while in \cite{Silverstein1995c}, $\alpha_n = 1$ is used.

We have equivalently: 
\begin{equation}\label{}
\begin{aligned}
	\alpha_n &= - \frac{1}{z\underline{m}_n(z)} \frac{1}{N} \sum_{j=1}^N  \frac{W_{jj}}{1 + r_j^* (B_{(j)} -zI)^{-1}r_j} \\
	\text{and } \alpha_n &= \frac{\frac{1}{N} \sum_{j=1}^N \frac{W_{jj}}{1 + r_j^* (B_{(j)} -zI)^{-1}r_j}}{\frac{1}{N} \sum_{j=1}^N \frac{1}{1 + r_j^* (B_{(j)} -zI)^{-1}r_j}},
\end{aligned}
\end{equation} 

Then:
\begin{equation}\label{}
\begin{aligned}
	&(-z \alpha_n \underline{m}_n(z) T_n - zI)^{-1} - (B_n - zI)^{-1} = \\
	&(-z \alpha_n \underline{m}_n(z) T_n - zI)^{-1}\left(z \alpha_n \underline{m}_n(z) T_n + \sum_{j=1}^N r_j r_j^*  \right) (B_n - zI)^{-1} =  \\
	&(-z \alpha_n \underline{m}_n(z) T_n - zI)^{-1}\sum_{j=1}^N \frac{r_j r_j^* (B_{(j)} -zI)^{-1} - \frac{1}{N}\alpha_n T_n (B_n -zI)^{-1}}{1 + r_j^* (B_{(j)} -zI)^{-1}r_j}.
\end{aligned}
\end{equation} 
Applying the trace and dividing by $n$, we obtain:
\begin{equation}\label{}
\begin{aligned}
	&\frac{1}{n} \tr\left((-z \alpha_n \underline{m}_n(z) T_n - zI)^{-1}\right) - m_n(z) = &&\\
	&-\frac{1}{N}\sum_{j=1}^N \frac{1}{z(1 + r_j^* (B_{(j)} -zI)^{-1}r_j)} \Big(W_{jj}q_j^* T^{1/2} (B_{(j)} - zI)^{-1} (\alpha_n \underline{m}_n(z) T_n + I)^{-1}T^{1/2} q_j \\
	& \qquad - \frac{1}{n} \tr \left( (\alpha_n \underline{m}_n(z) T_n + I)^{-1} \alpha_n T_n (B_n - zI)^{-1}\right)\Big)\\
\end{aligned}
\end{equation} 
So,
\begin{equation}\label{eq_cvg}
\begin{aligned}
	&\frac{1}{n} \tr\left((-z \alpha_n \underline{m}_n(z) T_n - zI)^{-1}\right) - m_n(z) = &&\\
	&-\frac{1}{zN}\sum_{j=1}^N q_j^* \left(\frac{W_{jj}T_n^{1/2} (B_{(j)} - zI)^{-1} (\alpha_n \underline{m}_n(z) T_n + I)^{-1}T_n^{1/2}}{1 + r_j^* (B_{(j)} -zI)^{-1}r_j} \right)q_j \\
	& + \frac{1}{zNn} \tr \left( (\alpha_n \underline{m}_n(z) T_n + I)^{-1} \sum_{j=1}^N \frac{\alpha_n T_n}{1 + r_j^* (B_{(j)} -zI)^{-1}r_j} (B_n - zI)^{-1}\right) = \\
	&-\frac{1}{N} \sum_{j=1}^N \frac{1}{z(1 + r_j^* (B_{(j)} -zI)^{-1}r_j)} d_j,
\end{aligned}
\end{equation} 
where:
\begin{equation}\label{}
\begin{aligned}
	&d_j = &&W_{jj}q_j^* T_n^{1/2} (B_{(j)} - zI)^{-1} (\alpha_n \underline{m}_n(z) T_n + I)^{-1}T_n^{1/2}q_j \\
	& &&- \frac{1}{n}\tr \left(W_{jj} (\alpha_n \underline{m}_n(z) T_n + I)^{-1}  T_n (B_n - zI)^{-1}\right).
\end{aligned}
\end{equation}

The strategy of the proof is the following:
\begin{itemize}
	\item prove that $\max_{j \leq N} d_j \longrightarrow 0$ a.s., and that $\frac{1}{n} \tr\left((-z \alpha_{n} \underline{m}_n(z) T_n - zI)^{-1}\right) - m_n(z) \longrightarrow 0$ a.s.,
	\item show that a.s. it exists $\tilde m(z) \in \C \backslash \C_+$ and a subsequence ${n_i}$ so that $\tilde m_{n_i}(z) := -z\alpha_{n_i}\underline{m}_{n_i}(z) \longrightarrow \tilde m(z)$ a.s.,
	\item prove that $\tilde m(z)$ is the unique solution in $\C_-$ of a functional equation, and deduce $\tilde m_{n}(z) \longrightarrow \tilde m(z)$ a.s.,
	\item deduce that a.s. it exists $ m(z) \in \C_+$, uniquely defined in function of $\tilde m(z)$ so that $m_n(z) \longrightarrow m(z)$ a.s.,
	\item similarly deduce that a.s. it exists $ \Theta^{(1)}(z) \in \C_+$, uniquely defined in function of $\tilde m(z)$ so that $\Theta^{(1)}_n(z) \longrightarrow \Theta^{(1)}(z)$ a.s..
\end{itemize}

\subsubsection{Decomposition of $d_j$}

In order to use the Lemma \ref{conc} on $d_j$, we decompose it into negligible terms and a term of the form $q_j^* C q_j - \frac{1}{n}tr(C)$ where $C$ is independent of $q_j$. 

For that, we denote: 
\begin{equation}\label{}
\begin{aligned}
	\tilde m_{n}(z) &:= -z\alpha_n\underline{m}_{n}(z) \\
	\text{and } \tilde m_{(j)}(z) &:= \frac{1}{N} \sum_{i \neq j}  W_{ii}\left(r_i^* (B_{(j)} -zI)^{-1}r_i - 1 \right).
\end{aligned}
\end{equation}

We have the following decomposition of $d_j$: 
\begin{equation}\label{}
\begin{aligned}
	d_j &= d_j^{(1)} + d_j^{(2)} + d_j^{(3)} + d_j^{(4)},
\end{aligned}
\end{equation}
with:
\begin{equation}\label{}
\begin{aligned}
	d_j^{(1)} &= W_{jj}q_j^* T_n^{1/2} (B_{(j)} - zI)^{-1} \left[\left(-\frac{\tilde{m}_n(z)}{z}  T_n + I\right)^{-1} - \left(-\frac{\tilde{m}_{(j)}(z)}{z}  T_n + I\right)^{-1} \right]T_n^{1/2}q_j, \\
	d_j^{(2)} &= W_{jj}q_j^* T_n^{1/2} (B_{(j)} - zI)^{-1}\left(-\frac{\tilde{m}_{(j)}(z)}{z}  T_n + I\right)^{-1}T_n^{1/2}q_j \\
	& \qquad\qquad\qquad\qquad\qquad\qquad - \frac{W_{jj}}{n}\tr \left( \left(-\frac{ \tilde{m}_{(j)}(z)}{z} T_n + I\right)^{-1} T_n (B_{(j)} - zI)^{-1}\right), \\
	d_j^{(3)} &= \frac{W_{jj}}{n}\tr \left(  \left[\left(-\frac{\tilde{m}_{(j)}(z)}{z}  T_n + I\right)^{-1} - \left(-\frac{\tilde{m}_{n}(z)}{z}  T_n + I\right)^{-1} \right] T_n (B_{(j)} - zI)^{-1}\right), \\
	d_j^{(4)} &= \frac{W_{jj}}{n}\tr \left(\left(-\frac{\tilde{m}_{n}(z)}{z}  T_n + I\right)^{-1} T_n \left[(B_{(j)} - zI)^{-1} - (B_{n} - zI)^{-1}\right]\right).
\end{aligned}
\end{equation}

In order to prove that for each $k \in \llbracket 1, 4 \rrbracket, \max_{j \leq N} d_j^{(k)} \longrightarrow 0$ a.s., we need some technical lemmas. They essentially provide the necessary inequalities to prove that $d_j^{(1)}, d_j^{(3)}$ and $d_j^{(4)}$ are indeed negligible, to finally use Lemma \ref{lemma_rmt} on $d_j^{(2)}$ and prove that  $\max_{j \leq N} |d_j| \longrightarrow 0$ a.s.

\subsubsection{Technical lemmas}

\begin{lem}\label{ineq}
	We have the following inequalities, for $j \in \llbracket 1,N \rrbracket$:
	\begin{equation}\label{}
	\begin{aligned}
		&\lVert (B_n - zI)^{-1} \rVert \leq \frac{1}{v}, \\
		&\lVert (B_{(j)} - zI)^{-1} \rVert \leq \frac{1}{v}, \\
		&\frac{1}{|z(1 + r_j^* (B_{(j)} -zI)^{-1}r_j)|} \leq \frac{1}{v}.
	\end{aligned}
	\end{equation} 
\end{lem}

\begin{proof}
	Let $j \in \llbracket 1, N\rrbracket$. The two first inequalities comes from the fact that $B_n$ and $B_{(j)}$ are Hermitian, so for any eigenvalue $\lambda$ of $B_n - zI$ or $B_{(j)} - zI$, we have that $|\lambda| \geq |\Im[\lambda]| = v$. 
	
	For the third inequality, we remark that:
	\begin{equation}\label{}
	\begin{aligned}
		\Im r_j^* (B_{(j)}/z -I)^{-1} r_j & = \frac{1}{2i} r_j^* \left[(B_{(j)}/z -I)^{-1} - (B_{(j)}/z^* -I)^{-1} \right] r_j \\
		& = \frac{v}{|z|^2} r_j^* (B_{(j)}/z -I)^{-1}B_{(j)} (B_{(j)}/ z^* -I)^{-1} r_j  \\
		\Im r_j^* (B_{(j)}/z -I)^{-1} r_j & \geq 0.
	\end{aligned}
	\end{equation} 
	So, we deduce that:
	\begin{equation}\label{}
	\begin{aligned}
		\frac{1}{|z(1 + r_j^* (B_{(j)} -zI)^{-1}r_j)|} \leq \frac{1}{v}.
	\end{aligned}
	\end{equation} 
\end{proof}

\begin{lem}\label{tech1}
	We denote $\bar W =  \frac{1}{N} \sum_{i=1}^N W_{ii}$. For $z = u +iv$, $v >0$ and $j \in \llbracket 1, N\rrbracket$, we have for any non-negative Hermitian matrix $A$:
	\begin{equation}\label{}
	\begin{aligned}
		&\left \lVert \left( - \frac{\tilde m_n(z)}{z} A + I \right)^{-1} \right \rVert \leq f(z,\lVert A \rVert), \text{ and }
		\left \lVert \left( - \frac{\tilde m_{(j)}(z)}{z} A + I \right)^{-1} \right \rVert \leq f(z,\lVert A \rVert),
	\end{aligned}
	\end{equation}
	where:
	\begin{equation}\label{}
	\begin{aligned}
		f(z,\lVert A \rVert) = \begin{cases}
			 \max \left(2, \frac{4}{v} \bar W\lVert A \rVert\right), &\text{ if } u= 0, \\
			 16 \left( \frac{|z|^2}{4v^2|u|}\bar W \lVert A \rVert +1 \right) \times \max \left(\frac{1}{3}, \frac{|u|}{v} \right), &\text{ otherwise.}
		\end{cases}
	\end{aligned}
	\end{equation}
\end{lem}

\begin{proof}
	Let z = u + iv$, v> 0, u \in \R$. For $j \in \llbracket 1,N \rrbracket$, we denote by $(u_{ij})_{i=1}^n$ a set of eigenvectors of $B_{(j)}$ with associated eigenvalues $(\lambda_{ij})_{i=1}^n$. We then derive the following formulation:
		\begin{equation}\label{}
		\begin{aligned}
			R := \Re\left[- \frac{\tilde m_n(z)}{z} \right] &= - \frac{1}{N} \sum_{j=1}^N \frac{W_{jj}}{\left|z + r_j^* \left(\frac{B_{(j)}}{z} - I \right)^{-1}r_j\right|^2} \left(u + \sum_{i=1}^n \frac{|r_j^*u_{ij}|^2}{|\lambda_{ij} - z|^2} (\lambda_{ij}u - |z|^2) \right),\\
			I :=\Im\left[- \frac{\tilde m_n(z)}{z} \right] &= \frac{1}{N} \sum_{j=1}^N \frac{W_{jj}}{\left|z + r_j^* \left(\frac{B_{(j)}}{z} - I \right)^{-1}r_j\right|^2} \left(v + \sum_{i=1}^n \frac{|r_j^*u_{ij}|^2}{|\lambda_{ij} - z|^2} \lambda_{ij}v \right) \geq 0.
		\end{aligned}
		\end{equation}
		Using Cauchy-Schwarz inequality, we deduce:
		\begin{equation}\label{}
		\begin{aligned}
			\frac{1}{N} \sum_{j=1}^N \frac{W_{jj}|z|^2\sum_{i=1}^n \frac{|r_j^*u_{ij}|^2}{|\lambda_{ij} - z|^2}}{\left|z + r_j^* \left(\frac{B_{(j)}}{z} - I \right)^{-1}r_j\right|^2}  \leq &\sqrt{\frac{1}{N} \sum_{j=1}^N \frac{W_{jj}\left( \sum_{i=1}^n \frac{|r_j^*u_{ij}|^2}{|\lambda_{ij} - z|^2}\right)^2}{\left|1 + r_j^* \left(B_{(j)} - zI \right)^{-1}r_j\right|^2} } \sqrt{\frac{1}{N} \sum_{j=1}^N \frac{W_{jj}}{\left|z + r_j^* \left(\frac{B_{(j)}}{z} - I \right)^{-1}r_j\right|^2}}  \\
			\frac{1}{N} \sum_{j=1}^N \frac{W_{jj}|z|^2\sum_{i=1}^n \frac{|r_j^*u_{ij}|^2}{|\lambda_{ij} - z|^2}}{\left|z + r_j^* \left(\frac{B_{(j)}}{z} - I \right)^{-1}r_j\right|^2}  \leq & \frac{|z|\sqrt{\bar W}}{v\sqrt{v}} \sqrt{I}.\\
		\end{aligned}
		\end{equation}
		So, combining both previous points, we have:
		\begin{equation}\label{sqrt_ineq}
		\begin{aligned}
			\left| R \right|& \leq \frac{|u|}{v} I + \frac{|z|\sqrt{\bar W}}{v\sqrt{v}} \sqrt{I}.
		\end{aligned}
		\end{equation}
		
		Suppose $u \neq 0$. We denote $K:=\frac{|z|\sqrt{\bar W}}{2v\sqrt{|u|}}$. We have then:
		\begin{equation}\label{}
		\begin{aligned}
			&\sqrt{\frac{v}{|u|}} \left(-  K + \sqrt{ K^2 + \left| R \right|}\right) \leq \sqrt{I}.
		\end{aligned}
		\end{equation}
		Now, let $x \geq 0$. Then:
		\begin{equation}\label{}
		\begin{aligned}
			\left|- \frac{\tilde m_n(z)}{z}x + 1 \right|^2 &= \left(Rx +1\right)^2 + I^2 x^2\\
			\left|- \frac{\tilde m_n(z)}{z}x + 1 \right|^2& \geq (-|R|x +1)^2  + \frac{v^2}{|u|^2} \left(-  K\sqrt{x} + \sqrt{ K^2x + \left| R \right|x}\right)^4.
		\end{aligned}
		\end{equation}
		We denote: $t := \sqrt{ K^2x +\left| R \right|x}\in \R_+$. We have then:
		\begin{equation}\label{}
		\begin{aligned}
			\left|- \frac{\tilde m_n(z)}{z}x + 1 \right|^2 &\geq (t^2 - K^2x - 1)^2 + \frac{v^2}{|u|^2} (t - K\sqrt{x})^4.
		\end{aligned}
		\end{equation}
		We denote $a:=K\sqrt{x}$, $b:= \sqrt{K^2x +1}$ and we split the study of the right part of the previous equation between $[0,(a+b)/2]$ and $[(a+b)/2,+\infty[$. The lower bounds rely mainly on the fact that $b- a \geq \frac{1}{2b}$:
		\begin{itemize}
			\item Let $t \in [0,(a+b)/2]$. Then,
				\begin{equation}\label{}
				\begin{aligned}
					(t^2 - K^2x - 1)^2 + \frac{v^2}{|u|^2} (t - K\sqrt{x})^4 &\geq \left(\left(\frac{a+b}{2}\right)^2 - b^2\right)^2 \\
					& = \frac{(a + 3b)^2(b-a)^2}{16} \\
					& \geq \frac{9(b-a)^4}{16} \\
					(t^2 - K^2x - 1)^2 + \frac{v^2}{|u|^2} (t - K\sqrt{x})^4 &\geq \frac{9}{16^2b^4}.
				\end{aligned}
				\end{equation}
			\item Let $t \in [(a+b)/2, +\infty[$. Then,
				\begin{equation}\label{}
				\begin{aligned}
					(t^2 - K^2x - 1)^2 + \frac{v^2}{|u|^2} (t - K\sqrt{x})^4 &\geq  \frac{v^2}{16|u|^2}(b-a)^4 \\
					(t^2 - K^2x - 1)^2 + \frac{v^2}{|u|^2} (t - K\sqrt{x})^4 &\geq  \frac{v^2}{16^2|u|^2b^4}. \\
				\end{aligned}
				\end{equation}
		\end{itemize}
		Backing up, we have that:
		\begin{equation}\label{}
		\begin{aligned}
			&\left|- \frac{\tilde m_n(z)}{z}x + 1 \right|^2 \geq \frac{1}{16^2\left(K^2x +1 \right)^2}\times \min\left(\frac{v^2}{|u|^2}, 9 \right).
		\end{aligned}
		\end{equation}
		It finally leads to:
		\begin{equation}\label{}
		\begin{aligned}
			&\left \lVert \left( - \frac{\tilde m_n(z)}{z} A + I \right)^{-1} \right \rVert \leq 16 \left( \frac{|z|^2}{4v^2|u|}\bar W \lVert A \rVert +1 \right) \times \max \left(\frac{1}{3}, \frac{|u|}{v} \right),
		\end{aligned}
		\end{equation}
		which proves the first inequality when $u \neq 0$.
		
		Now suppose $u= 0$. From Equation \eqref{sqrt_ineq}, we have:
		\begin{equation}
		\begin{aligned}
			\left| R \right|& \leq \sqrt{\frac{\bar W}{v}} \sqrt{I}.
		\end{aligned}
		\end{equation}
		If $\bar W=0$, then:
		\begin{equation}\label{}
		\begin{aligned}
			&\left|- \frac{\tilde m_n(z)}{z}x + 1 \right|^2 =1 + x^2 I^2 \geq 1 \geq \frac{1}{4}.
		\end{aligned}
		\end{equation}
		Otherwise, with $y = Rx$ and $a = \frac{v^2}{\bar W^4x^2}$:
		\begin{equation}\label{}
		\begin{aligned}
			&\left|- \frac{\tilde m_n(z)}{z}x + 1 \right|^2 \geq (Rx + 1)^2 +x^2 \frac{R^4v^2}{\bar W^4} = ay^4 + (y+1)^2.
		\end{aligned}
		\end{equation}
		Splitting the study between $]-\infty, 1/2]$ and $[1/2, +\infty[$, we have that $ay^4 + (y+1)^2 \geq \min \left(\frac{a}{16}, \frac{1}{4} \right).$ So, we deduce the first inequality when $u = 0$:
		\begin{equation}\label{}
		\begin{aligned}
			&\left \lVert \left( - \frac{\tilde m_n(z)}{z} A + I \right)^{-1} \right \rVert \leq \max\left(2, \frac{4\bar W\lVert A \rVert}{v} \right).
		\end{aligned}
		\end{equation}
		
		Let $j \in \llbracket 1, N \rrbracket$. Using the same method for $\left \lVert \left( - \frac{\tilde m_{(j)}(z)}{z} A + I \right)^{-1} \right \rVert$, we have with $\bar W_{(j)} := \frac{1}{N} \sum_{i \neq j} W_{ii}$ :
		\begin{equation}\label{}
		\begin{aligned}
			\left| \Re\left[- \frac{\tilde m_{(j)}(z)}{z} \right] \right|& \leq \frac{|u|}{v} \Im\left[- \frac{\tilde m_{(j)}(z)}{z} \right] + \frac{|z|\sqrt{\bar W_{(j)}}}{v\sqrt{v}}  \sqrt{ \Im\left[- \frac{\tilde m_{(j)}(z)}{z} \right]}.
		\end{aligned}
		\end{equation}
		As $W_{jj} \geq 0$, $\bar W_{(j)} \leq \bar W$, so we find the same equation as Equation \eqref{sqrt_ineq}:
		\begin{equation}\label{}
		\begin{aligned}
			\left| \Re\left[- \frac{\tilde m_{(j)}(z)}{z} \right] \right|& \leq \frac{|u|}{v} \Im\left[- \frac{\tilde m_{(j)}(z)}{z} \right] + \frac{|z|\sqrt{\bar W}}{v\sqrt{v}}  \sqrt{ \Im\left[- \frac{\tilde m_{(j)}(z)}{z} \right]}.
		\end{aligned}
		\end{equation}
		So, from the previous proof of the first inequality, we deduce immediately the second one, which complete the proof of this lemma:
		\begin{equation}\label{}
		\begin{aligned}
			&u \neq 0 \implies \left \lVert \left( - \frac{\tilde m_{(j)}(z)}{z} A + I \right)^{-1} \right \rVert \leq 16 \left( \frac{|z|^2}{4v^2|u|}\bar W \lVert A \rVert +1 \right) \times \max \left(\frac{1}{3}, \frac{|u|}{v} \right), \\
			&u = 0 \implies \left \lVert \left( - \frac{\tilde m_{(j)}(z)}{z} A + I \right)^{-1} \right \rVert \leq \max\left(2, \frac{4\bar W\lVert A \rVert}{v} \right).
		\end{aligned}
		\end{equation}		
\end{proof}

\begin{crl}\label{tech1_crl1}
	Let $j \in \llbracket 1, N \rrbracket$ and $z = u + iv, v > 0$. Then, for any matrices $A$ and $B$ of same size $n \times n$, $A$ Hermitian non-negative, we have:
	\begin{equation}\label{}
	\begin{aligned}
		&\left | \tr\left[B\left(\left( - \frac{\tilde m_{n}(z)}{z} A + I \right)^{-1} - \left( - \frac{\tilde m_{(j)}(z)}{z} A + I \right)^{-1}\right) \right]\right | \leq 
		\left|\frac{\tilde m_{n}(z)}{z} - \frac{\tilde m_{(j)}(z)}{z} \right|n \lVert B \rVert \lVert A \rVert f(z,\lVert A \rVert)^2.
	\end{aligned}
	\end{equation}
\end{crl}

\begin{proof}
	For any invertible matrices $C_1, C_2$ of the same size as $B$, we have:
	\begin{equation}\label{}
	\begin{aligned}
		|\tr\left[B(C_1^{-1} - C_2^{-1}) \right]| &= \tr\left[BC_1^{-1}(C_2-C_1)C_2^{-1} \right] \\
		|\tr\left[B(C_1^{-1} - C_2^{-1}) \right]| & \leq \lVert B \rVert \times \lVert C_1^{-1} \rVert \times \lVert C_2^{-1} \rVert \times n \lVert C_2-C_1\rVert.
	\end{aligned}
	\end{equation}
	From that point, the result comes immediately from Lemma \ref{tech1}.
\end{proof}

\begin{crl}\label{tech1_crl2}
	Let $j \in \llbracket 1, N \rrbracket$ and $z = u + iv, v > 0$. Then, for any $n \times n$ matrix $A$ Hermitian non-negative, and $r \in \C^n$, $\lVert r \rVert$ denoting its Euclidean norm, we have:
	\begin{equation}\label{}
	\begin{aligned}
		&\left | r^* \left(\left( - \frac{\tilde m_{n}(z)}{z} A + I \right)^{-1} - \left( - \frac{\tilde m_{(j)}(z)}{z} A + I \right)^{-1}\right) r\right | \leq 
		\left|\frac{\tilde m_{n}(z)}{z} - \frac{\tilde m_{(j)}(z)}{z} \right| \lVert r \rVert^2 \lVert A \rVert f(z,\lVert A \rVert)^2.
	\end{aligned}
	\end{equation}
\end{crl}

\begin{proof}
	For any invertible matrices $C_1, C_2$ of size $n \times n$, we have:
	\begin{equation}\label{}
	\begin{aligned}
		r^*(C_1^{-1} - C_2^{-1}) r &= r^*C_1^{-1}(C_2-C_1)C_2^{-1} r \\
		r^*(C_1^{-1} - C_2^{-1}) r & \leq \lVert r \rVert^2 \times \lVert C_1^{-1} \rVert \times \lVert C_2^{-1} \rVert \times \lVert C_2-C_1\rVert.
	\end{aligned}
	\end{equation}
	From that point, the result comes immediately from Lemma \ref{tech1}.
\end{proof}

\begin{lem}\label{tech2}
	Let $j \in \llbracket 1, N \rrbracket$. We denote: $A_{(j)} = \sum_{i \neq j} W_{ii} r_i r_i^*.$ Then,
	\begin{equation}\label{}
	\begin{aligned}
		\left|\frac{\tilde m_{n}(z)}{z} - \frac{\tilde m_{(j)}(z)}{z} \right| \leq \frac{2 d_2}{Nv} + \frac{\lVert A_{(j)} \rVert}{|z|vN}.
	\end{aligned}
	\end{equation}
\end{lem}

\begin{proof}
	Let $j \in \llbracket 1, N \rrbracket$. We have, using Lemma \ref{tech0}:
	\begin{equation}\label{}
	\begin{aligned}
		\left|\frac{\tilde m_{n}(z)}{z} - \frac{\tilde m_{(j)}(z)}{z} \right| & = &&\left|\frac{W_{jj}r_j^*(B_n - zI)^{-1}r_j}{zN} + \sum_{i \neq j} \frac{W_{ii}}{zN}  r_i^*\left[(B_n - zI)^{-1} - (B_{(j)} - zI)^{-1} \right]r_i \right| \\
		\left|\frac{\tilde m_{n}(z)}{z} - \frac{\tilde m_{(j)}(z)}{z} \right|& \leq &&\left|\frac{W_{jj}}{zN}\left(1 - \frac{1}{1 + r_j^*(B_{(j)} - zI)^{-1}r_j} \right) \right| \\
		& &&+ \left|\frac{1}{zN} \frac{r_j^*(B_{(j)} - zI)^{-1} A_{(j)} (B_{(j)} - zI)^{-1}r_j}{1+ r_j^* (B_{(j)} - zI)^{-1}r_j} \right|.
	\end{aligned}
	\end{equation}
	For the first term, using Lemma \ref{ineq}, we have:
	\begin{equation}\label{}
	\begin{aligned}
		\left|\frac{W_{jj}}{zN}\left(1 - \frac{1}{1 + r_j^*(B_{(j)} - zI)^{-1}r_j} \right) \right| \leq \frac{d_2}{N|z|} + \frac{d_2}{Nv} \leq \frac{2d_2}{Nv}.
	\end{aligned}
	\end{equation}
	For the second term, we have:
	\begin{equation}\label{}
	\begin{aligned}
		&\left|\frac{1}{zN} \frac{r_j^*(B_{(j)} - zI)^{-1} A_{(j)} (B_{(j)} - zI)^{-1}r_j}{1+ r_j^* (B_{(j)} - zI)^{-1}r_j} \right| \leq  \frac{\lVert A_{(j)}\rVert}{|z|N} \frac{\lVert (B_{(j)} - zI)^{-1}r_j \rVert^2}{|1+ r_j^* (B_{(j)} - zI)^{-1}r_j| } \\
		&\left|\frac{1}{zN} \frac{r_j^*(B_{(j)} - zI)^{-1} A_{(j)} (B_{(j)} - zI)^{-1}r_j}{1+ r_j^* (B_{(j)} - zI)^{-1}r_j} \right| \leq  \frac{\lVert A_{(j)}\rVert}{|z|vN}.
	\end{aligned}
	\end{equation}
	Using the proof of Lemma 2.6 \cite{Silverstein1995b}, we have additionally that:
	\begin{equation}\label{}
	\begin{aligned}
		&\frac{\lVert (B_{(j)} - zI)^{-1}r_j \rVert^2}{|1+ r_j^* (B_{(j)} - zI)^{-1}r_j| }  \leq \frac{1}{v}.
	\end{aligned}
	\end{equation}
	So,
	\begin{equation}\label{}
	\begin{aligned}
		&\left|\frac{1}{zN} \frac{r_j^*(B_{(j)} - zI)^{-1} A_{(j)} (B_{(j)} - zI)^{-1}r_j}{1+ r_j^* (B_{(j)} - zI)^{-1}r_j} \right| \leq  \frac{\lVert A_{(j)}\rVert}{|z|vN},
	\end{aligned}
	\end{equation}
	which concludes the proof.
\end{proof}

\begin{lem}\label{tech3}
	We have:
	\begin{equation}\label{}
	\begin{aligned}
		&\max_{j \leq N}\lVert q_j \rVert^2 \longrightarrow 1 \text{ a.s.}, \\
		&\max_{j \leq N} \frac{1}{N} \lVert A_{(j)} \rVert \longrightarrow 0 \text{ a.s.}
	\end{aligned}
	\end{equation}
\end{lem}

\begin{proof}
	The first convergence is a direct use of Lemma \ref{lemma_rmt} as remarked in \cite{Silverstein1995b} p.338.
	
	The second convergence comes from the following derivations for $j \in \llbracket 1,N \rrbracket$:
	\begin{equation}\label{}
	\begin{aligned}
		\frac{1}{N^2} \lVert A_{(j)} \rVert^2 &\leq \frac{1}{N^2} \tr \left[\left(\sum_{i \neq j} W_{ii} r_i r_i^* \right) \right] \\
		&\leq \frac{c_n^2}{N^2} \left(\sum_{i \neq j} W_{ii}^4 \lVert T \rVert^2 \lVert q_i \rVert^4 + \sum_{i \neq j}\sum_{i' \neq i, i' \neq j} W_{ii}^2 W_{i'i'}^2| q_{i'}^* T q_i|^2  \right) \\
		&\leq \frac{c_n^2}{N^2} \left(\sum_{i = 1}^N W_{ii}^4 \lVert T \rVert^2 \lVert q_i \rVert^4 + \sum_{i = 1}^N \sum_{i' \neq i} W_{ii}^2 W_{i'i'}^2| q_{i'}^* T q_i|^2  \right) \\
		&\leq \frac{c_n^2 d_2^4}{N^2} \left(\sum_{i = 1}^N h_2^2 \lVert q_i \rVert^4 + \sum_{i = 1}^N \sum_{i' \neq i} | q_{i'}^* T q_i|^2  \right) \\
		\frac{1}{N^2} \lVert A_{(j)} \rVert^2 &\leq \frac{c_n^2 d_2^4h_2}{N} \max_{i \leq N} \lVert q_i \rVert^4 + \frac{c_n^2d_2^4}{N} \max_{i \leq N} \sum_{i' \neq i} | q_{i'}^* T q_i|^2.
	\end{aligned}
	\end{equation}
	The upper bound does not depend on $j$ anymore, so we have:
	\begin{equation}\label{}
	\begin{aligned}
		\max_{j \leq N} \frac{1}{N^2} \lVert A_{(j)} \rVert^2 &\leq \frac{c_n^2 d_2^4h_2^2}{N} \max_{i \leq N} \lVert q_i \rVert^4 + \frac{c_n^2d_2^4}{N} \max_{i \leq N} \sum_{i' \neq i} | q_{i'}^* T q_i|^2.
	\end{aligned}
	\end{equation}
	From the first part of the proof, we have that $ \max_{i \leq N} \lVert q_i \rVert^4 \longrightarrow 1 \text{ a.s.}$, so:
	\begin{equation}\label{}
	\begin{aligned}
		\frac{c_n^2}{N} \max_{i \leq N} \lVert q_i \rVert^4 \longrightarrow 0 \text{ a.s.}
	\end{aligned}
	\end{equation}
	For the second term, we have:
	\begin{equation}\label{}
	\begin{aligned}
		\frac{c_n^2}{N} \max_{i \leq N} \sum_{i' \neq i} | q_{i'}^* T q_i|^2 &= \frac{c_n^2 }{N} \max_{i \leq N} \sum_{i' \neq i} q_{i'}^* (T q_i q_i^* T) q_{i'}.
	\end{aligned}
	\end{equation}
	Let $i \in \llbracket 1,N \rrbracket$. We use Lemma \ref{lemma_rmt} in dimension $n(N-1)$ with $Y = (Z_{ki'})_{k,i'\neq i} \in \R^{n \times (N-1)}$ (in vectorized form) and 
	$C = \begin{pmatrix}
		Tq_1 q_1^*T & & (0) \\
		& \ddots & \\
		(0) & & Tq_{N} q_{N}^*T
	\end{pmatrix}$ of size $n(N-1) \times n(N-1)$, where the index $i$ is removed in the matrix. We have then:
	\begin{equation}\label{}
	\begin{aligned}
		&\E\left[\left|\frac{ 1}{N} \sum_{i' \neq i} q_{i'}^* (T q_i q_i^* T) q_{i'} - \frac{(N-1)}{nN} \lVert T q_i \rVert^2 \right|^6\right] \leq \\
		& \frac{Kn^3(N-1)^3}{N^6n^6} \E[\lVert C \rVert^6] \leq \\
		&\frac{K}{n^3}h_2^{12}\sqrt{B}.
	\end{aligned}
	\end{equation}
	So, for all $\varepsilon > 0$:
	\begin{equation}\label{}
	\begin{aligned}
		&\mathbb{P}\left[\left|\max_{i \leq N} \frac{1}{N} \sum_{i' \neq i} q_{i'}^* (T q_i q_i^* T) q_{i'} - \max_{i \leq N}\frac{N-1}{nN} \lVert T q_i \rVert^2 \right| \geq \varepsilon\right] \leq \\ 
		&N\times \mathbb{P}\left[\left|\frac{1}{N} \sum_{i' \neq 1} q_{i'}^* (T q_1 q_1^* T) q_{i'} - \frac{N-1}{nN} \lVert T q_1 \rVert^2 \right| \geq \varepsilon\right] \leq \\
		&\frac{2KN}{\varepsilon^6n^3}h_2^{12}\sqrt{B},
	\end{aligned}
	\end{equation}
	which is summable. So,
	\begin{equation}\label{}
	\begin{aligned}
		\max_{i \leq N} \frac{1}{N} \sum_{i' \neq i} q_{i'}^* (T q_i q_i^* T) q_{i'} - \max_{i \leq N}\frac{(N-1)}{nN} \lVert T q_i \rVert^2 \longrightarrow 0 \text{ a.s.}
	\end{aligned}
	\end{equation}
	And, from the first part of the proof:
	\begin{equation}\label{}
	\begin{aligned}
		\left| \max_{i \leq N}\frac{c_n^2(N-1)}{nN} \lVert T q_i \rVert^2 \right| \leq \frac{c_n^2 h_2^2}{nN} \max_{i \leq N} \lVert q_i \rVert^2 \longrightarrow 0 \text{ a.s.}
	\end{aligned}
	\end{equation}
	We can conclude the proof:
	\begin{equation}\label{}
	\begin{aligned}
		&\max_{j \leq N} \frac{1}{N} \lVert A_{(j)} \rVert \longrightarrow 0 \text{ a.s.}
	\end{aligned}
	\end{equation}
\end{proof}

\subsubsection{Proof that $\max_{j \leq N} |d_j| \longrightarrow 0$ a.s.} 
\paragraph{Let us prove that $\max_{j \leq N} |d_j^{(1)}| \longrightarrow 0$ a.s.}
Let $j \in \llbracket 1, N \rrbracket$. We recall that:
\begin{equation}\label{}
\begin{aligned}
	d_j^{(1)} &= W_{jj}q_j^* T_n^{1/2} (B_{(j)} - zI)^{-1} \left[\left(-\frac{\tilde{m}_n(z)}{z}  T_n + I\right)^{-1} - \left(-\frac{\tilde{m}_{(j)}(z)}{z}  T_n + I\right)^{-1} \right]T_n^{1/2}q_j.
\end{aligned}
\end{equation}
Then, using Corollary \ref{tech1_crl2} and Lemma \ref{ineq}, we have:
\begin{equation}\label{}
\begin{aligned}
	|d_j^{(1)}| &\leq W_{jj} \lVert (B_{(j)} - zI)^{-1} \rVert \lVert T \rVert^2 \lVert q_j \rVert^2 \left|\frac{\tilde m_{n}(z)}{z} - \frac{\tilde m_{(j)}(z)}{z} \right| f(z,h_2)^2 \\
	|d_j^{(1)}| & \leq \frac{d_2h_2^2}{v}f(z,h_2)^2 \lVert q_j \rVert^2 \left|\frac{\tilde m_{n}(z)}{z} - \frac{\tilde m_{(j)}(z)}{z} \right|.
\end{aligned}
\end{equation}
Using Lemma \ref{tech2}, we have:
\begin{equation}\label{}
\begin{aligned}
	|d_j^{(1)}| & \leq \frac{d_2h_2^2}{v} f(z,h_2)^2  \left(\frac{2 d_2}{Nv} + \frac{\lVert A_{(j)} \rVert}{|z|vN}\right)\lVert q_j \rVert^2.
\end{aligned}
\end{equation}
By Assumption \ref{H4}, $\int x dF^{W_n}(x) \longrightarrow \int x dD(x) < \infty$ a.s., \textit{i.e.} $\bar W \longrightarrow \int x dD(x) < \infty$ a.s., so $\bar W$ is bounded a.s.. As a consequence, $f(z,h_2)^2 = O(1)$ a.s.. Finally, using Lemma \ref{tech3}, we can conclude:
\begin{equation}\label{}
\begin{aligned}
	\max_{j \leq N} |d_j^{(1)}| \longrightarrow 0 \text{ a.s.}
\end{aligned}
\end{equation}

\paragraph{Now let us prove that $\max_{j \leq N} |d_j^{(2)}| \longrightarrow 0$ a.s.}
Let $j \in \llbracket 1, N \rrbracket$. We recall that:
\begin{equation}\label{}
\begin{aligned}
	d_j^{(2)} = & W_{jj}q_j^* T_n^{1/2} (B_{(j)} - zI)^{-1}\left(-\frac{\tilde{m}_{(j)}(z)}{z}  T_n + I\right)^{-1}T_n^{1/2}q_j \\
	& - \frac{W_{jj}}{n}\tr \left( \left(-\frac{ \tilde{m}_{(j)}(z)}{z} T_n + I\right)^{-1} T_n (B_{(j)} - zI)^{-1}\right).
\end{aligned}
\end{equation}
Using Lemma \ref{lemma_rmt}, we have:
\begin{equation}\label{}
\begin{aligned}
	\E\left[\left|d_j^{(2)}\right|^6\right] & \leq \frac{K}{n^3} \left \lVert W_{jj} \left(-\frac{ \tilde{m}_{(j)}(z)}{z} T_n + I\right)^{-1} T_n (B_{(j)} - zI)^{-1} \right \rVert^6.
\end{aligned}
\end{equation}
Using Lemmas \ref{ineq} and \ref{tech1}, we have:
\begin{equation}\label{}
\begin{aligned}
	\E\left[\left|d_j^{(2)}\right|^6\right] & \leq \frac{K d_2^{12}}{v^6n^3} f(z,h_2)^6.
\end{aligned}
\end{equation}
As argued above, $f(z, \log(n)) = O(\log(n))$ a.s., so $\sum_{j=1}^N \E\left[\left|d_j^{(2)}\right|^6\right]$ is a.s. summable. So, we can conclude:
\begin{equation}\label{}
\begin{aligned}
	\max_{j \leq N} |d_j^{(2)}| \longrightarrow 0 \text{ a.s.}
\end{aligned}
\end{equation}

\paragraph{Then, let us prove that $\max_{j \leq N} |d_j^{(3)}| \longrightarrow 0$ a.s.}
Let $j \in \llbracket 1, N \rrbracket$. We recall that:
\begin{equation}\label{}
\begin{aligned}
	d_j^{(3)} &= \frac{W_{jj}}{n}\tr \left(  \left[\left(-\frac{\tilde{m}_{(j)}(z)}{z}  T_n + I\right)^{-1} - \left(-\frac{\tilde{m}_{n}(z)}{z}  T_n + I\right)^{-1} \right] T_n (B_{(j)} - zI)^{-1}\right).
\end{aligned}
\end{equation}
Using Corollary \ref{tech1_crl1} and Lemma \ref{ineq}, we have:
\begin{equation}\label{}
\begin{aligned}
	|d_j^{(3)}| &\leq \frac{d_2h_2^2}{nv} \left|\frac{\tilde m_{n}(z)}{z} - \frac{\tilde m_{(j)}(z)}{z} \right|n f(z, h_2)^2.
\end{aligned}
\end{equation}
Using Lemma \ref{tech2}, we have:
\begin{equation}\label{}
\begin{aligned}
	|d_j^{(3)}| &\leq \frac{d_2h_2^2}{v} \left(\frac{2 d_2}{Nv} + \frac{\lVert A_{(j)} \rVert}{|z| v N}\right) f(z, h_2)^2.
\end{aligned}
\end{equation}
As we argued previously, $f(z, h_2)^2 = O(1)$ a.s., so from Lemma \ref{tech3} we can conclude:
\begin{equation}\label{}
\begin{aligned}
	\max_{j\ \leq N} |d_j^{(3)}|\longrightarrow 0 \text{ a.s.}
\end{aligned}
\end{equation}

\paragraph{Finally, let us prove that $\max_{j \leq N} |d_j^{(4)}| \longrightarrow 0$ a.s.}
Let $j \in \llbracket 1, N \rrbracket$. We recall that:
\begin{equation}\label{}
\begin{aligned}
	d_j^{(4)} &= \frac{W_{jj}}{n}\tr \left(\left(-\frac{\tilde{m}_{n}(z)}{z}  T_n + I\right)^{-1} T_n \left[(B_{(j)} - zI)^{-1} - (B_{n} - zI)^{-1}\right]\right) \\
\end{aligned}
\end{equation}
Using Lemma \ref{tech0}, we have:
\begin{equation}\label{}
\begin{aligned}
	d_j^{(4)} &= \frac{W_{jj}}{n} \frac{r_j^* (B_{(j)} - zI)^{-1} \left(-\frac{\tilde{m}_{n}(z)}{z}  T_n + I\right)^{-1} T_n  (B_{(j)} - zI)^{-1}r_j}{1 + r_j^* (B_{(j)} - zI)^{-1} r_j}
\end{aligned}
\end{equation}
So,
\begin{equation}\label{}
\begin{aligned}
	|d_j^{(4)}| &= \frac{d_2}{n} \left \lVert \left(-\frac{\tilde{m}_{n}(z)}{z}  T_n + I\right)^{-1} T_n \right \rVert \frac{\lVert (B_{(j)} - zI)^{-1}r_j \rVert^2}{|1 + r_j^* (B_{(j)} - zI)^{-1} r_j|}
\end{aligned}
\end{equation}
Using the proof of Lemma 2.6 \cite{Silverstein1995b}, we have that:
\begin{equation}\label{}
\begin{aligned}
	&\frac{\lVert (B_{(j)} - zI)^{-1}r_j \rVert^2}{|1+ r_j^* (B_{(j)} - zI)^{-1}r_j| }  \leq \frac{1}{v}.
\end{aligned}
\end{equation}
So, using Lemma \ref{tech1}, we have:
\begin{equation}\label{}
\begin{aligned}
	|d_j^{(4)}| &= \frac{d_2h_2 f(z, h_2)}{vn}
\end{aligned}
\end{equation}
As argued before, $f(z, h_2) = O(1)$ a.s., so we can conclude:
\begin{equation}\label{}
\begin{aligned}
	\max_{j \leq N} |d_j^{(4)}| \longrightarrow 0 \text{ a.s.}
\end{aligned}
\end{equation}

We can now conclude this section. The last four points prove  that:
\begin{equation}\label{}
\begin{aligned}
	\max_{j \leq N} |d_j| \longrightarrow 0 \text{ a.s.}
\end{aligned}
\end{equation}
And from Lemma \ref{ineq}, we have for $j \in \llbracket 1, N \rrbracket$:
\begin{equation}\label{}
\begin{aligned}
	&\frac{1}{|z(1 + r_j^* (B_{(j)} -zI)^{-1}r_j)|} \leq \frac{1}{v}.
\end{aligned}
\end{equation}
So, 
\begin{equation}\label{}
\begin{aligned}
	&\frac{1}{N} \sum_{j=1}^N \frac{-1}{z(1 + r_j^* (B_{(j)} -zI)^{-1}r_j)} d_j \underset{n \rightarrow \infty}{\longrightarrow} 0 \text{ a.s.}
\end{aligned}
\end{equation}
Using Equation \eqref{eq_cvg}, we can now conclude that:
\begin{equation}\label{cvg}
\begin{aligned}
	&\frac{1}{n} \tr\left((\tilde{m}_n(z) T_n - zI)^{-1}\right) - m_n(z) \underset{n \rightarrow \infty}{\longrightarrow} 0 \text{ a.s.}
\end{aligned}
\end{equation} 

\subsubsection{Proof that a.s., $\left|\tilde m_n(z) \Theta_{n}^{(1)}(z)- (1 + zm_n(z)) \right| \underset{n \rightarrow \infty}{\longrightarrow} 0$} 
For this section, we recall the definition of $\Theta_n^{(1)}$ from Definition \ref{thetadef}:
\begin{equation}\label{}
\begin{aligned}
	\Theta_n^{(1)}(z) =\frac{1}{n} \tr \left((\mathrm{B}_n - zI)^{-1}T_n \right).
\end{aligned}
\end{equation}
We have:
\begin{equation}\label{}
\begin{aligned}
	\left|\tilde m_{n}(z) \Theta_{n}^{(1)}(z) - \left(1 + z m_{n}(z)\right)\right| &= \left|\tilde m_{n}(z) \Theta_{n}^{(1)}(z) - \frac{1}{n} \sum_{j=1}^N \frac{r_j^*(B_{(j)}-zI)^{-1}r_j}{1+r_j^*(B_{(j)}-zI)^{-1}r_j}\right| \\
	&= \left|\frac{1}{N} \sum_{j=1}^N \frac{W_{jj} \left(\Theta_{n}^{(1)}(z) - q_j^*(B_{(j)}-zI)^{-1}q_j\right)}{1+r_j^*(B_{(j)}-zI)^{-1}r_j}\right| \\
	\left|\tilde m_{n}(z) \Theta_{n}^{(1)}(z) - \left(1 + z m_{n}(z)\right)\right| & \leq \frac{|z|}{v}\max_{j \leq N}  d_2 \left|\Theta_{n}^{(1)}(z) - q_j^*T^{1/2}(B_{(j)}-zI)^{-1}T^{1/2}q_j \right|.
\end{aligned}
\end{equation} 
Using Lemma \ref{lemma_rmt}, we have that:
\begin{equation}\label{}
\begin{aligned}
	\E\left[\left|\Theta_{n}^{(1)}(z) - q_j^*T^{1/2}(B_{(j)}-zI)^{-1}T^{1/2}q_j \right|^6\right] \leq \frac{Kh_2^6}{n^3 v^6}.
\end{aligned}
\end{equation} 
So, 
\begin{equation}\label{theta_cvg}
\begin{aligned}
	\max_{j \leq N} \left|\Theta_{n}^{(1)}(z) - q_j^*T^{1/2}(B_{(j)}-zI)^{-1}T^{1/2}q_j \right| \underset{n \rightarrow \infty}{\longrightarrow} 0 \text{ a.s.}
\end{aligned}
\end{equation} 
So we can conclude that:
\begin{equation}\label{theta_cvg2}
\begin{aligned}
	\left|\tilde m_{n}(z) \Theta_{n}^{(1)}(z) - \left(1 + z m_{n}(z)\right)\right| \underset{n \rightarrow \infty}{\longrightarrow} 0 \text{ a.s.}
\end{aligned}
\end{equation} 

From Assumption \ref{H5}, $\bar W \rightarrow \int xdD(x) \in \R_+$ a.s. We focus now on trajectories where $\bar W \rightarrow \int xdD(x) \in \R_+$, Equations \eqref{cvg} and \eqref{theta_cvg} hold, $F^{T_n} \implies H$ and $F^{W_n} \implies D$.

Then, $(\tilde m_n(z))$ is bounded. Indeed, from Lemma \ref{ineq}, $| \tilde m_n(z) | \leq \bar W \frac{|z|}{v}$. Then, as $\Im[\tilde m_n(z)] \leq 0$,it exists a subsequence $\{n_i\}$ of $\N$ and $\tilde m(z) \in \C \backslash \C_+$ such that $\tilde m_{n_i}(z) \underset{i \rightarrow \infty}{\longrightarrow} \tilde m(z) \in \C \backslash \C_+$. 

\subsubsection{Proof that $m_{n_i}(z) \underset{i \rightarrow \infty}{\longrightarrow} \int \frac{1}{\tau \tilde m(z) - z}dH(\tau)$.} 
We want to prove that:
\begin{equation}\label{}
\begin{aligned}
	&m_{n_i}(z) -  \int \frac{1}{\tau \tilde m(z) - z}dH(\tau) \underset{i \rightarrow \infty}{\longrightarrow} 0.
\end{aligned}
\end{equation} 
From Equation \eqref{cvg}, it is equivalent to prove that:
\begin{equation}\label{}
\begin{aligned}
	&\int \frac{1}{\tau \tilde m_{n_i}(z) - z}dF^{T_{n_i}}(\tau) -  \int \frac{1}{\tau \tilde m(z) - z}dH(\tau) \underset{i \rightarrow \infty}{\longrightarrow} 0.
\end{aligned}
\end{equation} 

We prove that $\int \frac{1}{\tau \tilde m_{n_i}(z) - z}dF^{T_{n_i}}(\tau) \underset{i \rightarrow \infty}{\longrightarrow} -\frac{1}{z}$ using the Lebesgue's convergence theorem for weakly converging measures, as detailed in Corollary 5.1 \cite{Feinberg2019}. We denote: $f: \tau \in \R \rightarrow \frac{1}{\tau \tilde m(z) - z}$ and $f_{i}: \tau \in \R \rightarrow \frac{1}{\tau \tilde m_{n_i}(z) - z}$. Regarding the hypotheses of the theorem, we have:
\begin{itemize}
	\item $(f_i)_i$ is a.u.i. w.r.t. $(F^{T_{n_i}})_i$ (see (2.4) \cite{Feinberg2019} for a definition). Indeed, $\forall \tau \in \R_+, \forall i \in \R_+, |f_i(\tau)| \leq \frac{1}{v}$, so $\lim_{K \rightarrow +\infty} \limsup_{i \rightarrow +\infty} \int |f_i(\tau)| 1_{[K,+\infty[}(|f_i(\tau)|) dF^{T_{n_i}} = 0$.
	\item $F^{T_{n_i}} \implies H$ by assumption.
	\item Let $\tau \in \R_+$ and $\varepsilon > 0$. By assumption, $\bar W_{n_i}$ is bounded, and we denote by $\kappa$ one of its finite upper bound. By \eqref{ineq}, we have that $\forall i, |\tilde m_{n_i}(z) | \leq \frac{\kappa}{|z|v}$, so $\tilde m(z) \leq \frac{\kappa}{|z|v}$. We define:
		\begin{equation}\label{}
		\begin{aligned}
			\delta := \min\left(1,\frac{v^2}{\frac{\kappa}{|z|v} + \tau + 1}\right) > 0.
		\end{aligned}
		\end{equation} 
	Then, there exists $i_0 \in \N$ such that $\forall i \geq i_0,  |\tilde m_{n_i}(z) - \tilde m(z) | \leq \delta$. Now, let $\tau' \in ]\tau - \delta, \tau + \delta[ \cap \R_+$ and $i \geq i_0$. Then,
	\begin{equation}\label{}
	\begin{aligned}
		\left|f_i(\tau') - f(\tau) \right| &= \frac{|\tau' \tilde m_{n_i}(z) - \tau \tilde m(z)|}{|\tau' \tilde m_{n_i}(z) - z|\times|\tau \tilde m(z) - z|} \\
		&\leq \frac{1}{v^2}\left(|(\tau' - \tau) \tilde m(z)| + |\tau' (\tilde m_{n_i}(z)-\tilde m(z))| \right) \\
		&\leq \frac{\delta}{v^2}\left(\frac{\kappa}{|z|v} + \tau + 1 \right) \\
		\left|f_i(\tau') - f(\tau) \right| &\leq \varepsilon.
	\end{aligned}
	\end{equation} 
	So $\lim_{i \rightarrow \infty, \tau' \rightarrow \tau} f_i(\tau')$ exists.
\end{itemize}
So, using Corollary 5.1 \cite{Feinberg2019} on the real and imaginary parts, we deduce that $\lim_{i \rightarrow \infty} \int f_i(\tau) dF^{T_{n_i}}(\tau)$ exists and:
\begin{equation}\label{tau}
\begin{aligned}
	&\lim_{i \rightarrow \infty} \int f_i(\tau) dF^{T_{n_i}}(\tau) = \int f(\tau) dH(\tau).
\end{aligned}
\end{equation} 

It immediately leads to:
\begin{equation}
\begin{aligned}
	&m_{n_i}(z) -  \int \frac{1}{\tau \tilde m(z) - z}dH(\tau) \underset{i \rightarrow \infty}{\longrightarrow} 0.
\end{aligned}
\end{equation} 

\subsubsection{Case $\tilde m(z) = 0$.}
Suppose $\tilde m(z) = 0$. Then, $m_{n_i}(z) \underset{i \rightarrow \infty}{\longrightarrow} - \frac{1}{z}$ and by Cauchy-Stieltjes transform property $F^{B_{n_i}} \implies 1_{[0,+\infty[}$. Consequently, $\E\left[\frac{1}{n_i} \tr(B_{n_i})\right] \underset{i \rightarrow +\infty}{\longrightarrow} 0$. As $\E\left[\frac{1}{n_i} \tr(B_{n_i})\right] = \frac{1}{n_i}\tr(T_{n_i}) \times \bar W_{n_i}$. $\bar W_{n_i}$ converges by assumption to $\int \delta dD(\delta)$, and we supposed $D(]0,+\infty[) > 0$ so $\int \delta dD(\delta) > 0$. We deduce that $\frac{1}{n_i}\tr(T_{n_i}) \underset{i \rightarrow +\infty}{\longrightarrow} 0$. So, $H = 1_{[0,+\infty[}$, which is absurd. So $\tilde m(z) \neq 0$.

\subsubsection{Proof that $1 + zm_{n_i}(z) \underset{i \rightarrow \infty}{\longrightarrow} \int \frac{\delta \Theta^{(1)}(z)}{1 + \delta c \Theta^{(1)}(z)}dD(\delta)$ with $\Theta^{(1)}(z) := \frac{1+zm(z)}{\tilde m(z)}$.} 
We remind that we suppose now: $H(]0,+\infty[) > 0$, $D(]0,+\infty[) > 0$ and $\tilde m_{n_i}(z) \underset{i \rightarrow +\infty}{\longrightarrow} \tilde m(z) \neq 0$.

Firstly, we denote: $M_{n} := \log(n) \max_{j \leq N} \left|\Theta_{n}^{(1)}(z) - q_j^*T^{1/2}(B_{(j)}-zI)^{-1}T^{1/2}q_j \right|$. There is $n_0$ large enough so that $\forall n \geq n_0, M_n \leq \frac{v}{2c_n}$, due to Equation \eqref{theta_cvg}. Then, we have for $n_i \geq n_0$:
\begin{equation}\label{cvg_delta}
\begin{aligned}
	&\left| \left(-c_n(1 + zm_{n_i}(z)) +1\right) - \frac{1}{N}\sum_{j=1}^N \frac{1}{1 + W_{jj} c_{n_i} \Theta_{n_i}^{(1)}} \right| = \\
	&\left| \frac{1}{N}\sum_{j=1}^N \frac{1}{1 + r_j^*(B_{(j)}-zI)^{-1}r_j} - \frac{1}{1 + W_{jj} c_{n_i} \Theta_{n_i}^{(1)}} \right| \leq \\
	& \frac{c_{n_i} \bar W M_{n_i}}{v(v - c_{n_i} M_{n_i})} \underset{i \rightarrow \infty}{\longrightarrow} 0.
\end{aligned}
\end{equation} 
Then, we denote $\Theta^{(1)}(z) := \frac{1+zm(z)}{\tilde m(z)} = \int \frac{\tau}{\tau \tilde m(z) - z}dH(\tau)$. From Equation \eqref{theta_cvg2}, $\Theta^{(1)}_{n_i}\underset{i \rightarrow \infty}{\longrightarrow}  \Theta^{(1)}(z)$. We define:
\begin{equation}\label{}
\begin{aligned}
	&g: \delta \in \R \rightarrow \frac{1}{1 + \delta c \Theta^{(1)}(z)}.
\end{aligned}
\end{equation} 
We want to prove that:
\begin{equation}\label{}
\begin{aligned}
	&\frac{1}{N}\sum_{j=1}^N \frac{1}{1 + W_{jj} c_{n_i} \Theta_{n_i}^{(1)}} \underset{i \rightarrow \infty}{\longrightarrow} \int g(\delta) dD(\delta).
\end{aligned}
\end{equation} 
Let $\delta \in \R$. Remark that, as $\Im[\tilde m(z)]  \leq 0$:
\begin{equation}\label{im_theta}
\begin{aligned}
	&\Im[\Theta^{(1)}(z)] = \Im\left[\frac{1+zm(z)}{\tilde m(z)} \right] = \Im \left[\int \frac{\tau}{\tau \tilde m(z) - z}dH(\tau) \right] \geq v\int \frac{\tau}{|\tau \tilde m(z) - z|^2}dH(\tau) > 0.
\end{aligned}
\end{equation} 
So, from \cite{Silverstein1995c} p.338, we have:
\begin{equation}\label{ineq_delta1}
\begin{aligned}
	&|g(\delta)| \leq \frac{|\Theta^{(1)}(z)|}{\Im[\Theta^{(1)}(z)]}.
\end{aligned}
\end{equation} 
And, for $n_i$ large enough so that $\forall n_i \geq n_1, \Im[\Theta^{(1)}_{n_i}] > 0$:
\begin{equation}\label{ineq_delta2}
\begin{aligned}
	&\left| \frac{1}{1 + \delta c_{n_i} \Theta_{n_i}^{(1)}} - g(\delta)\right| \leq \frac{1}{c_{n_i} \Im[\Theta_{n_i}^{(1)}(z)]}\times \frac{|\Theta^{(1)}(z)|}{\Im[\Theta^{(1)}(z)]} \left|c\Theta^{(1)}(z) - c_{n_i} \Theta^{(1)}_{n_i}\right| .
\end{aligned}
\end{equation} 
Using \eqref{cvg_delta}, \eqref{ineq_delta1} and \eqref{ineq_delta2}, we have that:
\begin{equation}\label{}
\begin{aligned}
	&c_{n_i}\left(1 + zm_{n_i}(z)\right) - 1 \underset{i \rightarrow \infty}{\longrightarrow} \int \frac{ 1}{1 + \delta c \Theta^{(1)}(z)}dD(\delta).
\end{aligned}
\end{equation} 
We conclude that:
\begin{equation}\label{delta}
\begin{aligned}
	&1 + zm_{n_i}(z) \underset{i \rightarrow \infty}{\longrightarrow} \int \frac{ \delta \Theta^{(1)}(z)}{1 + \delta c \Theta^{(1)}(z)}dD(\delta).
\end{aligned}
\end{equation} 
Finally, using \eqref{tau} and \eqref{delta}, we have that:
\begin{equation}\label{tilde_m}
\begin{aligned}
	&\tilde m(z) = \int \frac{\delta}{1 + \delta c \int \frac{\tau}{\tau \tilde m(z) - z} dH(\tau)}dD(\delta).
\end{aligned}
\end{equation} 

\subsubsection{Uniqueness}\label{uni} 
In this section, we show that there is at most one solution $m \in \C \backslash \C_+$ so that:
\begin{equation}\label{fun_eq}
\begin{aligned}
	&m = \int \frac{\delta}{1 + \delta c \int \frac{\tau}{\tau m - z} dH(\tau)}dD(\delta).
\end{aligned}
\end{equation}

Let $m \in \C \backslash \C_+$, verifying \eqref{fun_eq}. Then, as $H(]0, +\infty[) > 0$:
\begin{equation}\label{uni_imag}
\begin{aligned}
	\Im[m] &= -\int \frac{\delta^2 c \int \frac{\tau(\tau \Im[-m] + v)}{|\tau m - z|^2}dH(\tau)}{\left| 1 + \delta c \int \frac{\tau}{\tau m - z}dH(\tau) \right|^2} dD(\delta) < 0.
\end{aligned}
\end{equation} 
So, with $\C_- = \{ y \in \C | \Im[y] < 0\}$:
\begin{equation}\label{}
\begin{aligned}
	\left[m \in \C \backslash \C_+ \text{ verifies \eqref{fun_eq} and } H(]0, +\infty[) > 0\right] \iff \left[m \in \C_- \text{ verifies \ref{fun_eq} and } H(]0, +\infty[) > 0\right].
\end{aligned}
\end{equation} 
Let $m_1 = u_1 - i v_1 \in \C_-$ and $m_2 = u_2 - i v_2 \in \C_-$ solving \eqref{fun_eq} at $z = u + iv \in \C_+$. Then, 
\begin{equation}\label{uni_diff}
\begin{aligned}
	m_1 - m_2 &= \int \frac{\delta^2 c \left(\int \frac{\tau}{\tau m_2 - z} dH(\tau) - \int \frac{\tau}{\tau m_1 - z} dH(\tau) \right)}{\left(1 + \delta c \int \frac{\tau}{\tau m_1 - z} dH(\tau)\right)\left(1 + \delta c \int \frac{\tau}{\tau m_2 - z} dH(\tau)\right)} dD(\delta) \\
	m_1 - m_2 &= (m_1 - m_2) \times \int \frac{\delta^2 c \int \frac{\tau^2}{(\tau m_1 - z)(\tau m_2 - z)} dH(\tau) }{\left(1 + \delta c \int \frac{\tau}{\tau m_1 - z} dH(\tau)\right)\left(1 + \delta c \int \frac{\tau}{\tau m_2 - z} dH(\tau)\right)} dD(\delta)
\end{aligned}
\end{equation} 
Using Hölder inequality on the last term and using \eqref{uni_imag} at the end, we have:
\begin{equation}\label{}
\begin{aligned}
	&\left| \int \frac{\delta^2 c \int \frac{\tau^2}{(\tau m_1 - z)(\tau m_2 - z)} dH(\tau) }{\left(1 + \delta c \int \frac{\tau}{\tau m_1 - z} dH(\tau)\right)\left(1 + \delta c \int \frac{\tau}{\tau m_2 - z} dH(\tau)\right)} dD(\delta) \right|^2 \leq \\
	& \left| \int \frac{\delta^2 c \int \frac{\tau^2}{|\tau m_1 - z|^2} dH(\tau) }{\left|1 + \delta c \int \frac{\tau}{\tau m_1 - z} dH(\tau)\right|^2} dD(\delta) \right| \times \left|\int \frac{\delta^2 c \int \frac{\tau^2}{|\tau m_2 - z|^2} dH(\tau) }{\left|1 + \delta c \int \frac{\tau}{\tau m_2 - z} dH(\tau)\right|^2} dD(\delta) \right| < \\
	& \left| \int \frac{\delta^2 c \int \frac{\tau^2 + \tau \frac{v}{v_1}}{|\tau m_1 - z|^2} dH(\tau) }{\left|1 + \delta c \int \frac{\tau}{\tau m_1 - z} dH(\tau)\right|^2} dD(\delta) \right| \times \left|\int \frac{\delta^2 c \int \frac{\tau^2+ \tau \frac{v}{v_2}}{|\tau m_2 - z|^2} dH(\tau) }{\left|1 + \delta c \int \frac{\tau}{\tau m_2 - z} dH(\tau)\right|^2} dD(\delta) \right| = \\
	& \left|\frac{v_1}{v_1} \right| \times \left|\frac{v_2}{v_2} \right| = 1.
\end{aligned}
\end{equation} 
Remark that the inequality is strict because we supposed $H(]0, +\infty[) > 0$. Injecting this inequation in \eqref{uni_diff}, we find:
\begin{equation}\label{}
\begin{aligned}
	|m_1 - m_2| \neq 0 \implies |m_1 - m_2| < |m_1 - m_2|.
\end{aligned}
\end{equation} 
So there is at most one solution $m \in \C \backslash \C_+$ verifying \eqref{fun_eq}.

\subsubsection{Convergences of $\tilde m_n(z)$, $m_n(z)$ and $\Theta^{(1)}_n(z)$.}
Backing up, we proved that almost surely, $\tilde m_n(z)$ is bounded and every convergent subsequence of $\tilde m_n(z)$ converge towards the unique $\tilde m(z) \in \C \backslash \C_+$ verifying \eqref{tilde_m}. So, a.s., $\tilde m_n(z) \underset{n \rightarrow \infty}{\longrightarrow} \tilde m(z)$.

We also proved that, almost surely, if $\tilde m_n(z) \underset{n \rightarrow \infty}{\longrightarrow} \tilde m(z)$ then:
\begin{itemize}
	\item $m_n(z) \underset{n \rightarrow \infty}{\longrightarrow} m(z) = \int \frac{1}{\tau \tilde m(z) - z} dH(\tau)$,
	\item $\Theta^{(1)}_n(z) \underset{n \rightarrow \infty}{\longrightarrow} \Theta^{(1)}(z) = \int \frac{\tau}{\tau \tilde m(z) - z} dH(\tau)$.
\end{itemize}

So, almost surely, $m_n(z) \underset{n \rightarrow \infty}{\longrightarrow} m(z) = \int \frac{1}{\tau \tilde m(z) - z} dH(\tau)$ where $\tilde m(z)$ is the unique solution in $\C \backslash \C_+$ to \eqref{tilde_m}.

\subsubsection{From $\tilde m$ to $X$}
	We have that for all $z \in \C_+$, $\tilde m(z)$ is the unique solution in $\C \backslash \C_+$ of the following equation:
	\begin{equation}\label{}
	\begin{aligned}
		\tilde m(z) &=\int \frac{\delta}{1 +  \delta c \int \frac{\tau}{\tau \tilde m(z) - z}dH(\tau)}dD(\delta).
	\end{aligned}
	\end{equation}
	We then define:
	\begin{equation}\label{}
	\begin{aligned}
		\forall z \in \C_+, X(z) := -\frac{\tilde m(z)}{z}.
	\end{aligned}
	\end{equation}
	Let $z \in \C_+$. It is proved that almost surely:
	\begin{equation}
	\begin{aligned}
		\tilde m_n(z) := \frac{1}{N} \sum_{j=1}^N \frac{W_{jj}}{1 + r_j^*(B_{(j)}-zI)^{-1}r_j} \underset{n \rightarrow +\infty} \longrightarrow \tilde m(z).
	\end{aligned}
	\end{equation}
	So, using that $X(z) := -\frac{\tilde m(z)}{z}$, we have immediately that, almost surely:
	\begin{equation}
	\begin{aligned}
		X_n(z) := -\frac{1}{N} \sum_{j=1}^N \frac{W_{jj}}{z + r_j^*(B_{(j)}/z-I)^{-1}r_j} \underset{n \rightarrow +\infty} \longrightarrow X(z).
	\end{aligned}
	\end{equation}
	Using the proof of Lemma \ref{ineq}, we have that $\Im[r_j^*(B_{(j)}/z-I)^{-1}r_j] \geq 0$, so $\forall n \in \N, X_n(z)\in \C_+$. So, $\Im[X(z)] \geq 0$. 

	Moreover, using the definition $X(z) = -\frac{\tilde m(z)}{z}$, it is easy to check that:
	\begin{equation*}\label{Xeq}
	\begin{aligned}
		X(z) &=-\int \frac{\delta}{z -  \delta c \int \frac{\tau}{\tau X(z) + 1}dH(\tau)}dD(\delta).
	\end{aligned}
	\end{equation*}
	
	It is also immediate that Equation \eqref{Xeq} does not admit any real solution because $z \in \C_+$. As we proved that $\Im[X(z)] \geq 0$, we deduce now that $X(z) \in \C_+$.

	The last point to check is the uniqueness in $\C_+$ of the solution of \eqref{Xeq}. Let $X_1, X_2 \in \C_+$, $X_1 \neq X_2$ solving Equation \eqref{Xeq}. Then, in the spirit of the proof of uniqueness for $\tilde m(z)$, we have with Hölder inequality:
	\begin{equation}
	\begin{aligned}
		|X_1 - X_2|^2 &= |X_1 - X_2|^2 \times \left|\int \frac{\delta^2c\int \frac{\tau^2}{(\tau X_1 + 1)(\tau X_2 + 1)}dH(\tau)}{\left(z - \delta c \int \frac{\tau}{\tau X_1 + 1}dH(\tau)\right)\left(z - \delta c \int \frac{\tau}{\tau X_2 + 1}dH(\tau)\right)}dD(\delta) \right|^2 \\
		&\leq |X_1 - X_2|^2 \times \left|\int \frac{\delta^2c\int \frac{\tau^2}{|\tau X_1 + 1|^2}dH(\tau)}{\left|z - \delta c \int \frac{\tau}{\tau X_1 + 1}dH(\tau)\right|^2}dD(\delta) \right| \times \left|\int \frac{\delta^2c\int \frac{\tau^2}{|\tau X_2 + 1|^2}dH(\tau)}{\left|z - \delta c \int \frac{\tau}{\tau X_2 + 1}dH(\tau)\right|^2}dD(\delta) \right| \\
		&< |X_1 - X_2|^2 \times \left|\frac{\Im[X_1]}{\Im[X_1]}\right| \times \left|\frac{\Im[X_2]}{\Im[X_2]}\right| \\
		|X_1 - X_2|^2 &< |X_1 - X_2|^2
	\end{aligned}
	\end{equation}
	The inequality is strict because $H(]0,\infty[) > 0$ and $D(]0,\infty[) > 0$ according to Assumption \ref{H3} and \ref{H4}. The result is absurd, hence the uniqueness in $\C_+$, which concludes the proof.

	In conclusion of this section, we proved the Lemma \ref{g1}, and the sketch of this proof, and its technical lemmas will be the backbone of the end of the proof of Theorem \ref{FEg}.

\subsection{Case $g = \Id$}
Now that we have proved the result and introduced the tools of the proof in the case $g = \mathrm{1}$, let us continue with the case $g = \Id$.
\begin{lem}
	Assume \ref{H1}-\ref{H5}. For all $z \in \C_+$, $\Theta^{(1)}_n(z) \rightarrow \Theta^{(1)}(z)$ almost surely and:
	\begin{equation*}\label{}
	\begin{aligned}
		 &\Theta^{(1)}(z) = -\frac{1}{z}\int \frac{\tau}{\tau X(z) + 1}dH(\tau),
	\end{aligned}
	\end{equation*}
	where for all $z \in \C_+$, $X(z)$ is the unique solution in $\C_+$ of the following equation:
	\begin{equation*}\label{}
	\begin{aligned}
		X(z) &=-\int \frac{\delta}{z -  \delta c \int \frac{\tau}{\tau X(z) + 1}dH(\tau)}dD(\delta).
	\end{aligned}
	\end{equation*}
\end{lem}

\begin{proof}
	Let $z \in \C_+$. We note $v := \Im[z]$. Using Equation \eqref{theta_cvg2}, we have:
	\begin{equation}\label{}
	\begin{aligned}
		\text{a.s., } \left|\tilde m_n(z) \Theta_{n}^{(1)}(z)- (1 + zm_n(z)) \right| \underset{n \rightarrow \infty}{\longrightarrow} 0.
	\end{aligned}
	\end{equation} 
	We also proved in Equation \ref{tilde_m} that a.s., $\tilde m_n(z) \underset{n \rightarrow \infty}{\longrightarrow} \tilde m(z) = \int \frac{\delta}{1 +  \delta c \int \frac{\tau}{\tau \tilde m(z) - z}dH(\tau)}dD(\delta)$. This equation ensures that $\tilde m(z) \neq 0$. To conclude, a.s., 
	\begin{equation}\label{}
	\begin{aligned}
		\Theta_{n}^{(1)}(z) \underset{n \rightarrow \infty}{\longrightarrow} \frac{1 +zm(z)}{\tilde m(z)} =  -\frac{1 +zm(z)}{zX(z)} = -\frac{1}{z}\int \frac{\tau}{\tau X(z) +1}dH(\tau) =: \Theta^{(1)}(z).
	\end{aligned}
	\end{equation} 
\end{proof}

\subsection{Case $g$ polynomial}
We continue the generalization of the result to polynomial functions, firstly by the following induction property.
\begin{lem}\label{lempoly}
For $q \in \N$ and $z \in \C_+$, $\Theta^{(q+1)}_n(z) \rightarrow \Theta^{(q+1)}(z)$ almost surely and:
\begin{equation}\label{}
\begin{aligned}
	 &\Theta^{(q+1)}(z) = -\frac{1}{zX(z)}\left[z\Theta^{(q)}(z) + \int \tau^q dH(\tau) \right].
\end{aligned}
\end{equation}
\end{lem}

\begin{proof}
	Let $z \in \C_+$ and $q \in \N$. We note $v := \Im[z]$. By induction, suppose a.s. $\Theta_n^{(q)}(z) \rightarrow \Theta^{(q)}(z)$. We have:
	\begin{equation}\label{}
	\begin{aligned}
			&\frac{1}{n}\tr\left(T_n^q + zT_n^q(B_n-zI)^{-1}\right) &= \frac{1}{n} \sum_{j=1}^N \frac{r_j^* T_n^q(B_{(j)}-zI)^{-1}r_j}{1 + r_j^*(B_{(j)}-zI)^{-1}r_j}.
	\end{aligned}
	\end{equation}
	Using Lemma \ref{conc} and elementary matrix calculus, we prove that:
	\begin{itemize}
		\item $\E\left[\left|r_j^*(B_{(j)}-zI)^{-1}r_j - \frac{W_{jj}}{N} \tr\left(T_n(B_{(j)}-zI)^{-1} \right) \right|^6\right] \leq \frac{d_2^6 h_2^6c_n^3}{v^6N^3}$,
		\item $\E\left[\left|r_j^*T_n^q(B_{(j)}-zI)^{-1}r_j - \frac{W_{jj}}{N} \tr\left(T_n^{q+1}(B_{(j)}-zI)^{-1} \right) \right|^6\right] \leq \frac{d_2^6 h_2^{6(q+1)}c_n^3}{v^6N^3}$,
		\item $\forall k \in \N$, $\underset{j \leq N}{\max} \left| \frac{1}{N}\tr\left(T_n^k(B_{(j)}-zI)^{-1} - T_n^k(B_{n}-zI)^{-1} \right) \right| \leq \underset{j \leq N}{\max}  \frac{|z|h_2^kd_2c_n\lVert q_j\rVert^2}{Nv^3} \rightarrow 0$ a.s.
	\end{itemize}
	Using these three points, we prove that:
	\begin{equation}\label{}
	\begin{aligned}
			&\frac{1}{n}\tr\left(T_n^q + zT_n^q(B_n-zI)^{-1}\right) - \frac{1}{n} \sum_{j=1}^N \frac{\frac{W_{jj}}{N} \tr\left(T_n^{q+1}(B_n-zI)^{-1} \right)}{1+\frac{W_{jj}}{N} \tr\left(T_n(B_n-zI)^{-1} \right)} \underset{n \rightarrow +\infty}{\longrightarrow} 0 \text{ a.s.}
	\end{aligned}
	\end{equation}
	Finally, we remark that, by assumption, a.s.:
	\begin{equation}\label{}
	\begin{aligned}
			&\frac{1}{n}\tr\left(T_n^q + zT_n^q(B_n-zI)^{-1}\right) \underset{n \rightarrow +\infty}{\longrightarrow} \int \tau^qdH(\tau) + z\Theta^{(q)}(z).
	\end{aligned}
	\end{equation}
	Moreover, 
	\begin{equation}\label{}
	\begin{aligned}
			&\frac{1}{n} \sum_{j=1}^N \frac{\frac{W_{jj}}{N} \tr\left(T_n^{q+1}(B_n-zI)^{-1} \right)}{1+\frac{W_{jj}}{N} \tr\left(T_n(B_n-zI)^{-1} \right)} = \Theta_n^{q+1}(z) \int \frac{\delta}{1 + \delta c_n \Theta_n^{(1)}(z)}dF^{W_n}(\delta).
	\end{aligned}
	\end{equation}
	As $\left|\frac{\delta}{1 + \delta c_n \Theta_n^{(1)}(z)}\right| \leq \frac{2}{\Im[c\Theta^{(1)}(z)]}$ for $n$ large enough, we have a.s.: 
	\begin{equation}\label{}
	\begin{aligned}
			&\int \frac{\delta}{1 + \delta c_n \Theta_n^{(1)}(z)}dF^{W_n}(\delta) \underset{n \rightarrow +\infty}{\longrightarrow} \int \frac{\delta}{1 + \delta c \Theta^{(1)}(z)}dD(\delta) = -zX(z).
	\end{aligned}
	\end{equation}
	To sum up, we have a.s.:
	\begin{equation}\label{}
	\begin{aligned}
			&\Theta_n^{(q+1)}(z) \underset{n \rightarrow +\infty}{\longrightarrow} -\frac{1}{zX(z)} \left[z\Theta^{(q)}(z) + \int \tau^q dH(\tau) \right].
	\end{aligned}
	\end{equation}
\end{proof}

This immediately leads to Theorem \ref{FEg} for polynomial functions.
\begin{lem}\label{poly}
	For $k \in \N$ and $z \in \C_+$, a.s. $\Theta_n^{(k)}(z) \underset{n \rightarrow +\infty}{\longrightarrow}  \Theta^{(k)}(z)$ and:
	\begin{equation}\label{}
	\begin{aligned}
		&\Theta^{(k)}(z) = -\frac{1}{z}\int \frac{\tau^k}{\tau X(z) + 1}dH(\tau) 
	\end{aligned}
	\end{equation}
	Moreover, let $g$ a real polynomial function. Then, for $z \in \C_+$, a.s. $\Theta_n^{g}(z) \underset{n \rightarrow +\infty}{\longrightarrow}  \Theta^{g}(z)$ and:
	\begin{equation}\label{}
	\begin{aligned}
		&\Theta^{g}(z) = -\frac{1}{z}\int \frac{g(\tau)}{\tau X(z) + 1}dH(\tau) 
	\end{aligned}
	\end{equation}
\end{lem}

\begin{proof}
	The proof is immediate by induction on $k$, using the previous lemmas. The linearity in $g$ proves the second part.
\end{proof}

\subsection{Case $g$ continuous}
\begin{lem}
	Assume \ref{H1}-\ref{H5}. For all $z \in \C_+$, $g$ continuous on $[h_1,h_2]$, $\Theta^{g}_n(z) \rightarrow \Theta^{g}(z)$ almost surely and:
	\begin{equation*}\label{}
	\begin{aligned}
		 &\Theta^{g}(z) = -\frac{1}{z}\int \frac{g(\tau)}{\tau X(z) + 1}dH(\tau),
	\end{aligned}
	\end{equation*}
	where for all $z \in \C_+$, $X(z)$ is the unique solution in $\C_+$ of the following equation:
	\begin{equation*}\label{}
	\begin{aligned}
		X(z) &=-\int \frac{\delta}{z -  \delta c \int \frac{\tau}{\tau X(z) + 1}dH(\tau)}dD(\delta).
	\end{aligned}
	\end{equation*}
\end{lem}

\begin{proof}
	By Weierstrass approximation theorem, there exists a sequence of polynomials $P_i$ that converge to $g$ uniformly on $[h_1,h_2]$. By Lemma \ref{poly}, we have that for all $z \in \C_+$, a.s. $\Theta_n^{P_i}(z) \underset{n \rightarrow +\infty}{\longrightarrow} \Theta^{P_i}(z)$. The, as $P_i$ converges uniformly to $g$ continuous, hence bounded, we have that $P_i$ is uniformly bounded. So, by D.C.T., for all $z \in \C_+$, $ \Theta^{P_i}(z)\underset{i \rightarrow +\infty}{\longrightarrow} \Theta^{g}(z)$. Moreover, for all $z \in \C_+$, $\Theta_n^{P_i}(z)$ converges uniformly to $\Theta_n^{g}(z)$. Indeed, we have from Lemma \ref{ineq}:
	\begin{equation}\label{}
	\begin{aligned}
		\left|\Theta_n^{P_i}(z) - \Theta_n^{g}(z)\right| &= \left|\frac{1}{n}\tr\left((B_n-zI)^{-1}\left(P_i(T_n)-g(T_n)\right)\right)\right|\\
		&\leq \lVert (B_n - zI)^{-1} \rVert \times \lVert P_i(T_n) - g(T_n) \rVert \\
		\left|\Theta_n^{P_i}(z) - \Theta_n^{g}(z)\right| &\leq \frac{ \lVert P_i(T_n) - g(T_n) \rVert}{v} \underset{i \rightarrow +\infty}{\longrightarrow} 0 \text{ uniformly on } n.
	\end{aligned}
	\end{equation}
	So by uniform inversion of limits, we have that:
	\begin{equation}\label{}
	\begin{aligned}
			&\Theta_n^{g}(z) \underset{n \rightarrow +\infty}{\longrightarrow} \Theta^{g}(z) = -\frac{1}{z}\int \frac{g(\tau)}{\tau X(z) + 1}dH(\tau),
	\end{aligned}
	\end{equation}
\end{proof}

\subsection{General case}
By induction on the number of discontinuity points, we prove Theorem \ref{FEg} in the general case. 

\begin{proof}
	Let $k \in \N$. Suppose Theorem \ref{FEg} holds for bounded functions with at most $k$ discontinuity points. Let $g: [h_1,h_2] \rightarrow \R$ bounded with $k+1$ discontinuity points, one of them is denoted by $\nu \in [h_1,h_2]$. We define $\rho: x \in [h_1,h_2] \mapsto g(x) \times (x - \nu)$. By assumption, $\forall z \in \C_+, \Theta_n^\rho \rightarrow \Theta^{\rho}(z) = -\frac{1}{z} \int \frac{\rho(\tau)}{\tau X(z) + 1}dH(\tau)$ a.s. Using the same proof as Lemma \ref{lempoly}, we prove that, for $z \in \C_+$, a.s.:
	\begin{equation}\label{}
	\begin{aligned}
			\Theta_n^\rho(z) = -\frac{1}{zX(z)}\left[\int g(\tau)dH(\tau) + z\Theta_n^g(z)\right] - \nu\Theta_n^g(z) + o(1).
	\end{aligned}
	\end{equation}
	We deduce then that, a.s. for $z \in \C_+$:
	\begin{equation}\label{}
	\begin{aligned}
		\Theta_n^\rho(z) \underset{n \rightarrow +\infty}\longrightarrow -\frac{zX(z)\Theta^\rho(z) + \int g(\tau)dH(\tau)}{z(1 + \nu X(z))} =: \Theta^g(z).
	\end{aligned}
	\end{equation}
	Using the assumption on $\Theta^\rho$, we have for $z \in \C_+$:
	\begin{equation}\label{}
	\begin{aligned}
		\Theta^g(z) &= -\frac{1}{z(1 + \nu X(z))}\left[-X(z) \int \frac{g(\tau)\times(\tau-\nu)}{\tau X(z) + 1}dH(\tau) + \int g(\tau)dH(\tau) \right] \\
		\Theta^g(z) &= -\frac{1}{z}\int \frac{g(\tau)}{\tau X(z) + 1}dH(\tau),
	\end{aligned}
	\end{equation}
	which concludes the proof.
\end{proof}

\section{Appendix: Proof of Theorem \ref{c0}}\label{apxC}
In the whole following section, we assume \ref{H1}-\ref{H5} to hold. We denote $\tilde m(z) := -zX(z)$. 

\begin{lem}
	$\tilde m$ is bounded on every bounded region of $\C_+$ bounded away from $\{0\}$.
\end{lem}

\begin{proof}
	Suppose there exists $(z_n)_n \in \C_+^\N$ bounded and bounded away from $\{0\}$, i.e. there exists $\eta > 0$ such that $|z_n| \geq \eta$, such that $|X(z_n)| \rightarrow +\infty$.
	Moreover, for $z \in \C_+$,
	\begin{equation}\label{}
	\begin{aligned}
		 &\left| \frac{\tau}{\tau X(z) + 1}\right| \geq \frac{h_2}{h_1 |X(z)| - 1}.
	\end{aligned}
	\end{equation}
	So, 
	\begin{equation}\label{}
	\begin{aligned}
		 &\left| \frac{\delta}{z + \delta c \int \frac{\tau}{\tau X(z) + 1} dH(\tau)}\right| \leq \frac{d_2}{|z| - cd_2\frac{h_2}{h_1 |X(z)| - 1}}.
	\end{aligned}
	\end{equation}
	By hypothesis, for $n$ large enough, $|X(z_n)| \geq \max\left(\frac{1}{h_1}\left(1 + \frac{2cd_2h_2}{\eta} \right), 2\frac{d_2}{\eta} \right) + 1$. Then, using the fact that for all $z \in \C_+$:
	\begin{equation}\label{funeq}
	\begin{aligned}
		 &X(z) = - \int \frac{\delta}{z - \delta c \int \frac{\tau}{\tau X(z) + 1}dH(\tau)}dD(\delta),
	\end{aligned}
	\end{equation}
	we have that:
	\begin{equation}\label{}
	\begin{aligned}
		 |X(z_n)| &\leq \int \left|\frac{\delta}{z_n - \delta c \int \frac{\tau}{\tau X(z) + 1}dH(\tau)}\right|dD(\delta) \\
		 &\leq \frac{d_2}{|z_n| - cd_2\frac{h_2}{h_1 |X(z_n)| - 1}} \\
		 &\leq 2\frac{d_n}{|z_n|} \\
		 |X(z_n)|&< |X(z_n)|,
	\end{aligned}
	\end{equation}
	which is absurd, and concludes the proof.
\end{proof}

From \cite{Zhang2007}, $m := \Theta^{(0)}$ is the Cauchy-Stieltjes transform of a probability distribution function - p.d.f. - that we denote $F$.
\begin{lem}
	Let $x_0 \in S_F^c, x_0 \neq 0$. Then, $\lim_{z \in \C_+ \rightarrow x_0} \tilde m(z)$ exists.
\end{lem}

\begin{proof}
	Let $x_0 \in S_F^c, x_0 \neq 0$. By property of the Cauchy-Stieltjes transform of a p.d.f., $m$ is analytic on $\C \backslash S_F$. Thus, $\lim_{z \in \C_+ \rightarrow x_0} m(z) := \check m(z) \in \R$ exists. Then, $\lim_{z \in \C_+ \rightarrow x_0} \Im[\tilde m(z)] = 0$. Now, suppose there exists $z_{1,n} \in \C_+^\N$ such that $z_{1,n} \rightarrow x_0$ and $\tilde m(z_{1,n}) \rightarrow \tilde m_1$, and similarly $z_{2,n} \in \C_+^\N$ such that $z_{2,n} \rightarrow x_0$ and $\tilde m(z_{2,n}) \rightarrow \tilde m_2$. As $m(z_{1,n}) \rightarrow \check m(z) \in \R$ and $m(z_{2,n}) \rightarrow \check m(z) \in \R$, we have that $\Re[m(z_{1,n}) - m(z_{2,n})] \rightarrow 0$. And, 
	\begin{equation}\label{}
	\begin{aligned}
		 \Re[m(z_{1,n}) - m(z_{2,n})] & = \int \frac{\tau \Re[\tilde m(z_{1,n}) - \tilde m(z_{2,n})] - \Re[z_{1,n} - z_{2,n}]}{|\tau \tilde m_{1,n} - z_{1,n}|^2 \times |\tau \tilde m_{2,n} - z_{2,n}|^2}dH(\tau).
	\end{aligned}
	\end{equation}
	So, 
	\begin{equation}\label{}
	\begin{aligned}
		&\left|\Re[\tilde m(z_{1,n}) - \tilde m(z_{2,n})]\right| = \\
		&\left|\frac{\Re[m(z_{1,n}) - m(z_{2,n})] + \Re[z_{1,n} - z_{2,n}] \int \frac{1}{|\tau \tilde m_{1,n} - z_{1,n}|^2 \times |\tau \tilde m_{2,n} - z_{2,n}|^2}dH(\tau)}{\int \frac{\tau}{|\tau \tilde m_{1,n} - z_{1,n}|^2 \times |\tau \tilde m_{2,n} - z_{2,n}|^2}dH(\tau)}\right| \leq \\
		& \left|\frac{\Re[m(z_{1,n}) - m(z_{2,n})]}{\int \frac{\tau}{|\tau \tilde m_{1,n} - z_{1,n}|^2 \times |\tau \tilde m_{2,n} - z_{2,n}|^2}dH(\tau)}\right| + \frac{1}{h_1}\left|\Re[z_{1,n} - z_{2,n}]\right|.
	\end{aligned}
	\end{equation}
	Suppose $\Re[\tilde m_{1}] \neq \Re[\tilde m_{2}]$. Then, from the previous inequation, we deduce that:
	\begin{equation}\label{}
	\begin{aligned}
		\int \frac{\tau}{|\tau \tilde m_{1,n} - z_{1,n}|^2 \times |\tau \tilde m_{2,n} - z_{2,n}|^2}dH(\tau) \rightarrow 0.
	\end{aligned}
	\end{equation}		
	So $|\tilde m_{1,n}| \rightarrow +\infty$ and $|\tilde m_{2,n}| \rightarrow +\infty$, which is absurd. So $\Re[\tilde m_{1}] = \Re[\tilde m_{2}]$. We proved at the beginning of the proof that $\Im[\tilde m_{1}] = \Im[\tilde m_{2}]$. As $\tilde m$ is bounded on every bounded region of $\C_+$ bounded away from $\{0\}$, we deduce that $\lim_{z \in \C_+ \rightarrow x_0} \tilde m(z)$ exists (and is equal to $\tilde m_1$).
\end{proof}

We define:
\begin{equation}\label{}
\begin{aligned}
	 &m_{LD}: z \in \C \backslash S_D \mapsto \int \frac{\delta}{\delta - z}dD(\delta).
\end{aligned}
\end{equation}
Remark that $m_{LD}$ is analytic on $\C \backslash S_D$.

\begin{lem}\label{leminvert}
	Let $x_0 \in S_F, x_0 \neq 0$. We denote $X(z) = -\frac{\tilde m(z)}{z}$. Suppose there exists $z_n \in \C_+ \rightarrow x_0$ such that $X(z_n) \rightarrow X_0 \in \C_+$ and $\hat z_n \in \C_+ \rightarrow x_0$ such that $X(\hat z_n) \rightarrow \hat X_0 \in \C \backslash \C_-$, $\hat X_0 \neq X_0$. Then there exists $\bar z_n \in \C_+ \rightarrow x_0$ such that $X(\bar z_n) \rightarrow \bar X_0 \in \C_+$ and $m_{LD}'(\check Y) \neq 0$ with $\check Y := \frac{x_0}{c\int \frac{\tau}{\tau \check X + 1}dH(\tau)}$.
\end{lem}

\begin{proof}
	Consider the following procedure:
	\begin{itemize}
		\item[(i)] Suppose $\Im[X_0] \neq \Im[\hat X_0]$. Let $n \in \N$. By continuity of $X$ on $\C_+$, there exists $u_n \in \C_+$ in the complex segment $[z_n, \hat z_n]_\C$ such that $\Im[X(u_n)] = \frac{1}{2}\left(\Im[X(z_n)] + \Im[X(\hat z_n)]\right)$. $u_n \rightarrow x_0 \neq 0$ so $X(u_n)$ is bounded, so there exists an extraction $\bar z_n^{(1)}$ of $u_n$ such that $\bar z_n^{(1)} \rightarrow x_0$ and $X(\bar z_n^{(1)}) \rightarrow \bar X^{(1)} \in \C_+$, $\Im[\bar X^{(1)}] = \frac{1}{2}\left(\Im[X_0] + \Im[\hat X_0]\right)$. Do the procedure (i) again with $z_n, \bar z_n^{(1)}$ instead of $z_n, \hat z_n$ in order to construct $\bar z_n^{(2)}$ and $\bar X^{(2)}$.
		\item[(ii)] Suppose $\Im[X_0] = \Im[\hat X_0] > 0$, then $\Re[X_0] \neq \Re[\hat X_0]$. Similarly, there exists $u_n \in \C_+$ in the complex segment $[z_n, \hat z_n]_\C$ such that $\Re[X(u_n)] = \frac{1}{2}\left(\Re[X(z_n)] + \Re[X(\hat z_n)]\right)$. $u_n \rightarrow x_0 \neq 0$ so $X(u_n)$ is bounded, so there exists an extraction $\bar v_n$ of $u_n$ such that $\bar v_n \rightarrow x_0$ and $X(\bar v_n) \rightarrow \bar X \in \C \backslash \C_-$, $\Re[\bar X] = \frac{1}{2}\left(\Re[X_0] + \Re[\hat X_0]\right)$. 
		\begin{itemize}
			\item If $\Im[\bar X] = 0$, do the procedure (i) with $z_n, \bar v_n$ instead of $z_n, \hat z_n$ in order to construct $\bar z_n^{(1)}$ and $\bar X^{(1)} \in \C_+$.
			\item If $\Im[\bar X] = \Im[X_0]$, define $z_n^{(1)} := v_n$ and $\bar X^{(1)} := \bar X \in \C_+$, and do the procedure (ii) again with $z_n, \bar z_n^{(1)}$ instead of $z_n, \hat z_n$ in order to construct $\bar z_n^{(2)},...$ and $\bar X^{(2)} \in \C_+$.
			\item Otherwise, also define $z_n^{(1)} := v_n$ and $\bar X^{(1)} := \bar X  \in \C_+$, and go to procedure (i) with $z_n, \bar z_n^{(1)}$ instead of $z_n, \hat z_n$ in order to construct $\bar z_n^{(2)}$ and $\bar X^{(2)} \in \C_+$.
		\end{itemize}
	\end{itemize}
	Using this procedure, we construct iteratively $\{ \bar X^{(k)}, k\in \N^*\}$, where for all $k \in \N^*, \bar X^{(k)} \in \C_+$ is an adherence point of $X$ as $z \in \C_+ \rightarrow x_0$, and by construction, $\forall k \neq k' \in \N^*, \bar X^{(k)} \neq \bar X^{(k')}$. Suppose that $\forall k \in \N^*, m_{LD}'\left( \frac{x_0}{c\int \frac{\tau}{\tau \bar X^{(k)} + 1}dH(\tau)}\right) = 0$. $D \neq 0$ by hypothesis so $m_{LD}$ is not constant over $\C_+$, and as $m_{LD}$ is analytic on $\C_+$ it implies that $\left\{\int \frac{\tau}{\tau \bar X^{(k)} + 1}dH(\tau), k \in \N^* \right\}$ is finite while $\left\{ X^{(k)}, k \in \N^* \right\}$ is countable. So, as $\phi: X \in \C_+ \mapsto \int \frac{\tau}{\tau X + 1}dH(\tau)$ is analytic, we have that $\phi$ is constant over $\C_+$. But $H \neq 0$ by hypothesis, this is absurd and it concludes the proof.
\end{proof}

\begin{lem}
	Let $x_0 \in S_F, x_0 \neq 0$. Suppose there exists $z_n \in \C_+ \rightarrow x_0$ such that $X(z_n) \rightarrow \check{X} \in \C_+$ and $m_{LD}'(\check Y) \neq 0$ with $\check Y := \frac{x_0}{c\int \frac{\tau}{\tau \check X + 1}dH(\tau)}$. Then $\lim_{z \in \C_+ \rightarrow x_0} X(z)$ exists and $\lim_{z \in \C_+ \rightarrow x_0} X(z) = \check{X}$.
\end{lem}

\begin{proof}
	We have for all $z \in \C_+$, with $X(z) = -\frac{\tilde m(z)}{z}$:
	\begin{equation}\label{}
	\begin{aligned}
		 &m_{LD}\left(\frac{z}{c\int \frac{\tau}{\tau X(z) + 1} dH(\tau)} \right) = c - c \int \frac{1}{\tau X(z) + 1} dH(\tau).
	\end{aligned}
	\end{equation}
	Suppose $z_n \in \C_+ \rightarrow x_0$, $X(z_n) \rightarrow \check{X} \in \C_+$ and, with $\check Y := \frac{x_0}{c\int \frac{\tau}{\tau \check X + 1}dH(\tau)} \in \C_+$ we have $m_{LD}'(\check Y) \neq 0$. So by the holomorphic inverse function theorem, $m_{LD}$ is locally invertible at $\check Y$ in an open set $B$ containing $\check Y$. We denote its local inverse $g_{LD}: m_{LD}(B) \rightarrow B$, which is also analytic.
	
	As $Y_n := \frac{z_n}{c\int \frac{\tau}{\tau X(z_n) + 1} dH(\tau)} \rightarrow \check Y$, for $n$ large enough, $Y_n \in B$. We have then, for $n$ large enough:
	\begin{equation}\label{}
	\begin{aligned}
		 &z_n  = c\int \frac{\tau}{\tau X(z_n) + 1} dH(\tau)g_{LD}\left(c - c \int \frac{1}{\tau X(z_n) + 1} dH(\tau)\right).
	\end{aligned}
	\end{equation}
	So, 
	\begin{equation}\label{}
	\begin{aligned}
		 &x_0  = c\int \frac{\tau}{\tau \check X + 1} dH(\tau)g_{LD}\left(c - c \int \frac{1}{\tau \check X + 1} dH(\tau)\right).
	\end{aligned}
	\end{equation}
	Let $\varepsilon > 0$ and $B' = B(\check X, \varepsilon)$ so that $\left \{ c - c \int \frac{1}{\tau X + 1} dH(\tau), X \in B' \right \} \subset m_{LD}(B)$. We can choose $\varepsilon > 0$ arbitrarily small because $m_{LD}(B)$ is open due to the open mapping theorem, $c - c \int \frac{1}{\tau \check X + 1} dH(\tau) \in m_{LD}(B)$ and $X \mapsto \int \frac{1}{\tau X + 1} dH(\tau)$ is analytic on $\C \backslash \{x \in \R, -x^{-1} \in S_H\}$.
	We define:
	\begin{equation}\label{}
	\begin{aligned}
		 &z_F: X \in B' \mapsto c\int \frac{\tau}{\tau X + 1} dH(\tau)g_{LD}\left(c - c \int \frac{1}{\tau X + 1} dH(\tau)\right).
	\end{aligned}
	\end{equation}
	By the open mapping theorem, as $z_F$ is analytic and non-constant over $B'$, $z_F(B')$ is open and contains $x_0$. Thus, for any $(\bar z_n)_n \in \C_+^\N$ such that $\bar z_n \rightarrow x_0$, we have that $\bar z_n \in z_F(B')$ for $n$ large enough. So, for these $\bar z_n$, there exists $X_n \in B'$ such that $z_F(X_n) = \bar z_n$. So, $X_n = X(\bar z_n) \in B'$. As we can choose $\varepsilon > 0$ arbitrarily small, we have that $ X(\bar z_n) \rightarrow \check X$, which concludes the proof.
\end{proof}

At this stage, we know that, for $x_0 \in S_F^c \backslash \{0\}$, $\lim_{z \in \C_+ \rightarrow x_0} \tilde m(z)$ exists, and for $x_0 \in S_F \backslash \{0\}$, either $\lim_{z \in \C_+ \rightarrow x_0} \Im[\tilde m(z)] = 0$ or $\lim_{z \in \C_+ \rightarrow x_0} \tilde m(z)$ exists.

\begin{lem}\label{lem8}
	Assume $0 \notin ]m_a, m_b[ \subset \R$. If, for all $m_0 \in ]m_a, m_b[$, there exists $(h_n)_n \in (\R_+^*)^\N$ such that $h_n \rightarrow 0$ and $\int \frac{h_n\tau^2dH(\tau)}{(1+\tau m_0)^2 + \tau^2h_n^2} \rightarrow 0$, then $]-m_a^{-1}, -m_b^{-1}[ \subset S_H^c$.
\end{lem}

\begin{proof}
	Let $m_0 \in ]m_a, m_b[$ and $n \in \N$. Let $g(\tau) = \frac{h_n\tau^2}{(1+\tau m_0)^2 + \tau^2h_n^2}$. $g$ is increasing on $]-m_0^{-1} - h_n,-m_0^{-1}]$ and decreasing on $]-m_0^{-1},-m_0^{-1}+h_n]$. Then,
	\begin{equation}\label{}
	\begin{aligned}
		 \int g(\tau) dH(\tau) &\geq \int_{-m_0^{-1} - h_n}^{-m_0^{-1}} \frac{h_n(m_0^{-1}+h_n)^2dH(\tau)}{(h_n m_0)^2 + h_n^2(m_0^{-1}+h_n)^2} \\
		 \int g(\tau) dH(\tau) &\geq \frac{(m_0^{-1} + h_n)^2}{m_0^2 + (m_0^{-1} + h_n)^2} \frac{H(-m_0^{-1}) - H(-m_0^{-1} - h_n)}{h_n}.
	\end{aligned}
	\end{equation}
	Thus, $\frac{H(-m_0^{-1}) - H(-m_0^{-1} - h_n)}{h_n} \rightarrow 0$. Then, the lower-left Dini derivative (see \cite{Hagood2006} for a reference on the subject) is null, \textit{i.e.} $D_- H(-m_0^{-1}) := \liminf_{h \rightarrow 0+} \frac{H(-m_0^{-1}) - H(-m_0^{-1} - h)}{h} = 0$. As it is null, $D_- H$ is continuous on $]-m_a^{-1}, -m_b^{-1}[$. Moreover, $H$ is càdlàg, non-decreasing and $D_- H = 0$ on $]-m_a^{-1}, -m_b^{-1}[$ so $H$ is continuous on $]-m_a^{-1}, -m_b^{-1}[$. The two latter points imply that the three other Dini derivatives - $D^-H, D_+H$ and $D^+H$ - are continuous on $]-m_a^{-1}, -m_b^{-1}[$. As $H$ is monotone on $]-m_a^{-1}, -m_b^{-1}[$, it is almost everywhere differentiable due to Lebesgue's theorem. So, almost everywhere on $]-m_a^{-1}, -m_b^{-1}[$, all four Dini derivatives are equal, and their value is $0$ due to $D_-$. And from continuity, everywhere on $]m_a, m_b[$, all four Dini derivatives are equal to $0$. So $H$ is differentiable on $]-m_a^{-1}, -m_b^{-1}[$ and $H' = 0$. So $H$ is constant on $]-m_a^{-1}, -m_b^{-1}[$, which finally implies that $]-m_a^{-1}, -m_b^{-1}[ \subset S_H^c$.
\end{proof}

\begin{lem}\label{lem36}
	Assume $\lim_{z \in \C_+ \rightarrow x_0} \Im[X(z)] = 0$. Let $z_n \in \C_+ \rightarrow x_0$ and $\hat z_n \in \C_+ \rightarrow x_0$ such that $X(z_n) \rightarrow X_0 \in \R$ and $X(\hat z_n) \rightarrow \hat X_0 \in \R$, $X_0 < \hat X_0$. Then, $\forall \bar X \in ]X_0, \hat X_0[$, there exists $(\bar z_n)_n \in \C_+^\N$ such that $\bar z_n \rightarrow x_0$ and $X(\bar z_n) \rightarrow \bar X$. $(\bar z_n)_n$ can be chosen so that $\Re[X(\bar z_n)] = \bar X$.
\end{lem}

\begin{proof}
	The same proof as Lemma 3.6 \cite{Silverstein1995a} applies here.
\end{proof}

\begin{lem}\label{lem10}
	Let $x_0 \in S_F, x_0 \neq 0$. Suppose $\lim_{z \in \C_+ \rightarrow x_0} \Im[X(z)] = 0$, then $\lim_{z \in \C_+ \rightarrow x_0} X(z)$ exists.
\end{lem}

\begin{proof}
	Suppose $\lim_{z \in \C_+ \rightarrow x_0} \Im[X(z)] = 0$ but $\Re[X(z)]$ does not converge as $z \in \C_+ \rightarrow x_0$. Using Lemma \ref{lem36}, we can find $X_0 \in \R$ and $\hat X_0 \in \R$ adherence points of $X$ at $x_0$, $X_0 < \hat X_0$ such that $0 \notin ]X_0, \hat X_0[$. Let $\bar X \in ]X_0, \hat X_0[$. From Lemma \ref{lem36}, there exists $(\bar z_n)_n \in \C_+^\N$ such that $\bar z_n \rightarrow x_0$ and $X(\bar z_n) \rightarrow \bar X$ and $\Re[X(\bar z_n)] = \bar X$. So, $\Im[X(\bar z_n)] \rightarrow 0$. By the function equation on $X$, we have:
	\begin{equation}\label{}
	\begin{aligned}
		 \Im[X(\bar z_n)] &= \int \frac{\delta \left(\Im[\bar z_n] + \delta c \int \frac{\tau^2 \Im[X(\bar z_n)]}{|\tau X(\bar z_n) + 1|^2}dH(\tau)\right)}{\left|\bar z_n - \delta c \int \frac{\tau}{\tau X(\bar z_n) + 1}dH(\tau) \right|^2}.
	\end{aligned}
	\end{equation}
	$\left|\bar z_n - \delta c \int \frac{\tau}{\tau X(\bar z_n) + 1}dH(\tau) \right|^2$ is bounded, so $\int \frac{\tau^2 \Im[X(\bar z_n)]}{|\tau X(\bar z_n) + 1|^2}dH(\tau) \rightarrow 0$. Using Lemma \ref{lem8} with $m_a := X_0$, $m_b := \hat X_0$ and $h_n := \Im[X(\bar z_n)]$, we deduce that $]-X_0^{-1}, -\hat X_0^{-1}[ \in S_H^c$.
	Without loss of generality, suppose $\left]c - c\int \frac{1}{\tau X_0 +1}dH(\tau), c - c\int \frac{1}{\tau \hat X_0 +1}dH(\tau)\right[ \subset \R^*$. 
	For $\bar z_n \rightarrow x_0$ such that $X(\bar z_n) \rightarrow \bar X \in ]X_0, \bar X_0[$, we have that $m_{LD}\left(\frac{\bar z_n}{c\int \frac{\tau}{\tau X(\bar z_n) + 1} dH(\tau)} \right) = c - c \int \frac{1}{\tau X(\bar z_n) + 1} dH(\tau) \rightarrow c - c\int \frac{1}{\tau \bar X +1}dH(\tau) \in \R^*$. Thus, $\int \frac{\tau}{\tau X(\bar z_n) + 1} dH(\tau) \rightarrow \int \frac{\tau}{\tau \bar X + 1} dH(\tau) \neq 0$ (otherwise $m_{LD}\left(\frac{\bar z_n}{c\int \frac{\tau}{\tau X(\bar z_n) + 1} dH(\tau)} \right) \rightarrow 0$), and $\frac{\bar z_n}{c\int \frac{\tau}{\tau X(\bar z_n) + 1} dH(\tau)} \rightarrow \frac{x_0}{c\int \frac{\tau}{\tau \bar X +1}dH(\tau)} \in \R$.
	Then, $\left] \frac{x_0}{c\int \frac{\tau}{\tau X_0 +1}dH(\tau)},\frac{x_0}{c\int \frac{\tau}{\tau \hat X_0 +1}dH(\tau)} \right[ \subset S_D^c$.
	
	So, for $\bar X \in ]X_0, \hat X_0[$, $m_{LD}$ is locally invertible at $\bar Y = \frac{x_0}{c\int \frac{\tau}{\tau \bar X +1}dH(\tau)}$, on an open $B$ containing this point, we denote its inverse $g_{LD}$. Indeed, $m_{LD}'(\bar Y) = \int \frac{\delta}{(\delta-\bar Y)^2}dD(\delta) > 0$.
	
	We define $B' := \left\{X \in \C^* | -X^{-1} \notin S_H \text{ and }c - c \int \frac{1}{\tau X + 1} dH(\tau) \in m_{LD}(B) \right\}$. Remark that $B'$ is open as $m_{LD}(B)$ open and $X \in \{X \in \C^* | -X^{-1} \notin S_H \} \mapsto c - c \int \frac{1}{\tau X + 1} dH(\tau)$ analytic. We define $z_F: X \in B' \mapsto c\int \frac{\tau}{\tau X + 1} dH(\tau)g_{LD}\left(c - c \int \frac{1}{\tau X + 1} dH(\tau)\right).$ 
	For all $\bar X \in B'$, we have $x_0 = z_F(\bar X)$, so $z_F$ is constant on $B'$, which is absurd. So $\lim_{z \in \C_+ \rightarrow x_0} X(z)$ exists.
\end{proof}

With the previous lemmas, we finally proved that $\lim_{z \in \C_+ \rightarrow x_0} X(z)$ exists for $x_0 \in \R^*$. The point $x_0 = 0$ will be discussed in the end, as we need the following results to prove it. We now prove the existence and formula on $\check \Theta^g$.
\begin{lem}\label{lemImThetag} 
	For $\lambda \in \R^*$ and $g: [h_1,h_2] \rightarrow \R$ bounded with a finite number of discontinuity points:
	\begin{equation}\label{}
	\begin{aligned}
		&\check \Theta^{g}(\lambda) = \frac{1}{\lambda} \int \frac{g(\tau)}{\tau \check X(\lambda) + 1}dH(\tau).
	\end{aligned}
	\end{equation}
\end{lem}

\begin{proof}
	Let $\lambda \in \R^*$. For every $R \subset \C_+$ bounded and bounded away from the imaginary axis, we have that $z \in R \mapsto \int \frac{1}{|1 + X(z)\tau|^2}dH(\tau)$ is bounded from Lemma 3.2 \cite{Couillet2015}. By Fatou's lemma, we have that $\tau \mapsto \frac{1}{|1 + \check X(\lambda)\tau|^2}$ is $H$-mesurable and $\int \frac{1}{|1 + \check X(\lambda)\tau|^2}dH(\tau) < +\infty$. 

	Then, for a neighborhood of $\lambda$ in $\C_+$, there exists a constant $\kappa > 0$ such that $H$-a.s., for all $z$ in this neighborhood, $|\tau X(z) + 1|^2 \geq \kappa |\tau \check X(\lambda) + 1|^2$, using the continuity of the function $z \mapsto \sup_{\tau \in S_H} \frac{|\tau \check X(\lambda) + 1|^2}{|\tau X(z) + 1|^2}$.

	So, for $z \in \C_+$ in this neighborhood of $\lambda$, $X(z)$ is bounded because continuous, say by $K$, we have for $\tau \in S_H$:
	\begin{equation}\label{}
	\begin{aligned}
		&\left|\frac{g(\tau)}{\tau X(z) + 1}\right| \leq \lVert g \rVert_\infty \frac{Kh_2+1}{\kappa|\tau \check X(\lambda) + 1|^2}.
	\end{aligned}
	\end{equation}
	We proved that the right part is integrable, so by D.C.T.:
	\begin{equation}\label{}
	\begin{aligned}
		&\Theta^g(z) \underset{z \in \C_+ \rightarrow \lambda}{\longrightarrow} \check \Theta^{g}(\lambda) = -\frac{1}{\lambda} \int \frac{g(\tau)}{\tau \check X(\lambda) + 1}dH(\tau).
	\end{aligned}
	\end{equation}
\end{proof}

We now extend the convergence to $\lambda=0$ in the case $c < 1$.
\begin{lem}\label{lemcneq1}
	Suppose $c<1$, $\lim_{z \in \C_+ \rightarrow 0} m(z)$ exists and for $g$ bounded with a finite number of discontinuities, $\lim_{z \in \C_+ \rightarrow 0} \Theta^{g}(z)$ exists.
\end{lem}

\begin{proof}
	Using Equation (1.9b) \cite{Bai2004} for $c<$, we have for $\eta = (1-\sqrt{c})^2h_1d_1$ and any $l >0$, $\mathbb{P}(\lambda_{\min}^{B_n} \leq \eta) = o(n^{-l})$. So, $F_n(\eta) \underset{n \rightarrow +\infty}{\longrightarrow} 0$. So, for any $x < \eta$, $F(x) = 0$, in particular $F(0) = 0$ and $F$ continuous in $0$. Consequently, $\lim_{z \in \C_+ \rightarrow 0} m(z) =: \check m(0)$ exists and $\check m(0) \in \R^*$, as $S_F \subset \R^*$. Moreover, $-zX(z)m(z) \underset{z \in \C_+ \rightarrow 0}{\longrightarrow} \int \frac{1}{\tau}dH(\tau) \in \R^*$, so, using Theorem \ref{FEg}, we have by D.C.T. that $\Theta^g(z) \underset{z \in \C_+ \rightarrow 0}{\longrightarrow} \frac{\int \check m(0)g(\tau)dH(\tau)}{\int \tau^{-1} dH(\tau)}$.
\end{proof}

We finish by extending $\check X$ to $\R$ in the case $c > 1$.
\begin{lem}\label{lemcheckX0}
	$X$ - extended with $\check X$ on $\R^*$ if $c > 1$ otherwise on $\R^*$ - is a continuous function.
\end{lem}

\begin{proof}
	We only need to prove that $\lim_{z \in \C_+ \rightarrow 0} X(z)$ exists if $c > 1$, and the lemma is then a direct consequence of Theorem 2.2 \cite{Silverstein1995a}.

	Suppose $c > 1$. Then from the same argument used in the proof of Lemma \ref{lemcneq1} but applied to $\underline{m}$, we have that $\underline{m}(z) \underset{z \in \C_+ \rightarrow 0}{\longrightarrow} \check{\underline{m}}(0) \in \R$. As $c(1+zm(z)) = 1 + z\underline{m}(z)$, we have $\int \frac{1}{\tau X(z) + 1} dH(\tau) = -zm(z) \underset{z \in \C_+ \rightarrow 0}{\longrightarrow} 1 - \frac{1}{c}$. From that property, it is easy to see that $X(z)$ is bounded and bounded away from $0$ in a neighborhood of $0$. So, the adherence points of $X$ in $0$ is a non-empty set included in $\R^*$, denoted $\mathcal{A}$. From Lemma \ref{lem36} we have that $\mathcal{A}$ is a convex set. Let $X_1 < X_2 \in \mathcal{A}$. Using Lemma \ref{lem8} with the subsequences created with Lemma \ref{lem36}, we have that $]-X_1^{-1},-X_2^{-1}[ \in S_H^c$. By D.C.T, we then have that $\forall X \in ]-X_1^{-1},-X_2^{-1}[, \int \frac{1}{\tau X + 1} dH(\tau) = 1 - \frac{1}{c}$. So, $X_1 = X_2$, and $\mathcal{A}$ has only one element. Thus, $\check X(0) = \lim_{z \in \C_+ \rightarrow 0} X(z)$ exists and $\check X(0) \in \R^*$.
\end{proof}

With the last Lemma, we completely proved Theorem \ref{c0}. 

\section{Appendix: Proof of shrinkage formulas and eigenvectors inner product behavior}\label{apxD}
\subsection{Proof of Theorem \ref{csh}}
Theorem \ref{FEg} shows that:
\begin{equation}\label{}
\begin{aligned}
	 \forall z \in \C_+, \Theta_n^{(1)}(z) \rightarrow \Theta^{(1)}(z) = -\frac{1}{z} \int \frac{\tau}{\tau X(z) + 1 }dH(\tau) \text{ a.s.}
\end{aligned}
\end{equation}
Moreover, we have that, from Equation (12) in \cite{Ledoit2009}, for every $x$ where $\Delta_n$ is continuous:
\begin{equation}\label{}
\begin{aligned}
	 \Delta_n(x) = \lim_{\eta \rightarrow 0^+} \frac{1}{\pi} \int_{-\infty}^x \Im\left[\Theta_n^{(1)}(\lambda+i\eta) \right]d\lambda.
\end{aligned}
\end{equation}

The following Lemma generalizes to the weighted setting the Lemma 3.1 \cite{Ledoit2009}.
\begin{lem}\label{lem31}
	Let $g$ a real bounded function on $[h_1, h_2]$ with finitely many points of discontinuity. Consider the function $\Omega_n^g$ defined by:
	\begin{equation}\label{}
	\begin{aligned}
		 \forall x \in \R, \Omega_n^g(x) = \frac{1}{n} \sum_{i=1}^n 1_{[\lambda_i,+\infty[}(x)\sum_{j=1}^n|u_i^*v_j|^2 \times g(\tau_j).
	\end{aligned}
	\end{equation}
	Then, there exists a nonrandom function $\Omega^g$ defined on $\R$ such that $\Omega_n^g(x) \rightarrow \Omega^g(x)$ a.s. at all points of continuity of $\Omega^g$. Furthermore,
	\begin{equation}\label{eq32}
	\begin{aligned}
		 \Omega^g(x) = \lim_{\eta \rightarrow 0^+} \frac{1}{\pi} \int_{-\infty}^x \Im\left[\Theta^g(\lambda+i\eta)\right]d\lambda,
	\end{aligned}
	\end{equation}
	for all $x$ where $\Omega^g$ is continuous.
\end{lem}

\begin{proof}
	The proof is very similar to the proof of Lemma 3.1 \cite{Ledoit2009}. The Cauchy-Stieltjes transform of $\Omega_n^g$ is $\Theta_n^g$. From Theorem \ref{FEg}, $\forall z \in \C_+, \Theta^g_n(z) \rightarrow \Theta^g(z)$ a.s. Therefore, Equation (2.5) \cite{Silverstein1995b} implies that $\lim_{n \rightarrow +\infty} \Omega_n^g(x) =: \Omega^g(x)$ exists for all $x$ where $\Omega^g$ is continuous. Furthermore, the Cauchy-Stieltjes transform of $\Omega^g$ is $\Theta^g$. Then, \eqref{eq32} is simply the Stieltjes-Perron formula, inversion of the Stieltjes transform.
\end{proof}

By Lemma \ref{lem31}, $\Delta(x) := \lim_{n \rightarrow +\infty} \Delta_n(x)$ exists a.s. at all point of continuity of $\Delta$ and is equal to:
\begin{equation}\label{}
\begin{aligned}
	 \Delta(x) = \lim_{\eta \rightarrow 0^+} \frac{1}{\pi} \int_{-\infty}^x \Im\left[\Theta^{(1)}(\lambda+i\eta) \right]d\lambda.
\end{aligned}
\end{equation}

We need the following lemma to invert limit and integral.
\begin{lem}\label{lem32}
	Under the assumptions of Lemma \ref{lem31} and that $c<1$, we have for $(x_1, x_2) \in \R^2$:
	\begin{equation}\label{}
	\begin{aligned}
		 \Omega^g(x_2) - \Omega^g(x_1) = \frac{1}{\pi} \int_{x_1}^{x_2} \lim_{\eta \rightarrow 0^+} \Im[\Theta^g(\lambda+i\eta)]d\lambda.
	\end{aligned}
	\end{equation}
	If $c>1$, the result holds for $x_1x_2 > 0$.
\end{lem}

\begin{proof}
	We know from Theorem \ref{c0} that $\forall x \in \R^*, \lim_{z \in \C_+ \rightarrow x} \Im[\Theta^g(z)] =: \Im[\check{\Theta}^g(x)]$ exists.
	Let us prove it for $x = 0$. From Lemma 3.3 \cite{Ledoit2009} and Marcenko-Pastur theorem for $c<1$, we have that $F$ is constant over the interval $]-\infty, (1-\sqrt{c})^2h_1d_1[$. So $\lim_{z \in \C_+ \rightarrow 0} \Im[m(z)] \rightarrow 0$. So, $\lim_{z \in \C_+ \rightarrow 0} \Im[\Theta^g(z)] \rightarrow 0$.

	$\Theta^g$ is the Cauchy-Stieltjes transform $\Omega^g$. Theorem 2.1 \cite{Silverstein1995c} implies that $\Omega^g$ is differentiable at $x \in \R$ if $c<1$, $x \in \R^*$ if $c > 1$, and its derivative is $\pi^{-1}\Im[\check \Theta^g(\cdot)]$. We get the result by integration.
\end{proof}

We can now complete the proof of Theorem \ref{csh} with the following lemma.
\begin{lem}
	Suppose $c<1$. Let $x \in \R$.
	\begin{equation*}\label{}
	\begin{aligned}
		 \Delta(x) = \int_{-\infty}^x h(\lambda) dF(\lambda),
	\end{aligned}
	\end{equation*}
	where for $\lambda \in \R^* \backslash \{x \in S_F, F'(x) = 0\}$:
	\begin{equation*}\label{}
	\begin{aligned}
		h(\lambda) = \frac{\int \frac{\tau^2}{|\tau \check X(\lambda) + 1|^2}dH(\tau)}{\int \frac{\tau}{|\tau \check X(\lambda) + 1|^2}dH(\tau)}.
	\end{aligned}
	\end{equation*}
\end{lem} 

\begin{proof}
	We have, using Lemma \ref{lem32}:
	\begin{equation}\label{}
	\begin{aligned}
		\Delta(x) &= \lim_{\eta \rightarrow 0^+} \frac{1}{\pi} \int_{-\infty}^x \Im\left[\Theta_n^{(1)}(\lambda+i\eta) \right]d\lambda \\
		\Delta(x) &= \int_{-\infty}^x \frac{1}{\pi}\Im[\check{\Theta}^{(1)}(\lambda)]d\lambda.		
	\end{aligned}
	\end{equation}
	Remark that $\forall \lambda \in \R, \Im[\check{\Theta}^{(1)}(\lambda)] = 0 \iff \Im[\check{m}(\lambda)] = 0$. So, as $F$ has a density $F' := \pi^{-1}\Im[\check m(\cdot)]$ (see \cite{Couillet2015}):
	\begin{equation}\label{}
	\begin{aligned}
		\Delta(x) &= \int_{]-\infty,x[\cap (F')^{-1}(\R_+^*)} \frac{1}{\pi}\Im[\check{\Theta}^{(1)}(\lambda)]d\lambda \\
		 &= \int_{]-\infty,x[\cap (F')^{-1}(\R_+^*)} \frac{\Im[\check \Theta^{(1)}(\lambda)]}{\Im[\check{m}(\lambda)]}\times \frac{1}{\pi}\Im[\check{m}(\lambda)]d\lambda \\	
		\Delta(x) &= \int_{]-\infty,x[\cap (F')^{-1}(\R_+^*)} \frac{\Im[\check \Theta^{(1)}(\lambda)]}{\Im[\check{m}(\lambda)]}dF(\lambda).
	\end{aligned}
	\end{equation}
	Let $\lambda \in (F')^{-1}(\R_+^*)$. We proved that $F'(0) = \pi^{-1} \Im[\check m(0)] = 0$, so $\lambda \neq 0$. Then, using Theorem \ref{c0}:
	\begin{equation}\label{}
	\begin{aligned}
		\frac{\Im[\check \Theta^{(1)}(\lambda)]}{\Im[\check{m}(\lambda)]} &= \frac{\Im\left[\int \frac{\tau}{\tau \check X(\lambda) + 1}dH(\tau)\right]}{\Im\left[\int \frac{1}{\tau \check X(\lambda) + 1}dH(\tau)\right]} \\
		&= \frac{\int \frac{\tau^2}{|\tau \check X(\lambda) + 1|^2}dH(\tau)}{\int \frac{\tau}{|\tau \check X(\lambda) + 1|^2}dH(\tau)} \\
		\frac{\Im[\check \Theta^{(1)}(\lambda)]}{\Im[\check{m}(\lambda)]} &= h(\lambda).
	\end{aligned}
	\end{equation}
	$h(\cdot)$ is well-defined on $\R^* \backslash \{ x\in S_F | F'(x) = 0 \}$, due to Theorem \ref{c0}. We extend $h$ on $\R$ with any constant value  $\lambda \in \{ x\in S_F | F'(x) = 0 \} \cup \{ 0 \}$, for example, for $h(\lambda) := 0$. The value does not matter as $dF\left((F')^{-1}(\{0\})\right) = 0$ and $dF\left(\{0\}\right) = 0$ as $c<1$.
	Finally:
	\begin{equation}\label{}
	\begin{aligned}
		\Delta(x) &= \int_{]-\infty,x[\cap (F')^{-1}(\R_+^*)} h(\lambda) dF(\lambda) \\
		\Delta(x) &= \int_{-\infty}^x h(\lambda) dF(\lambda),
	\end{aligned}
	\end{equation}
	which concludes the proof.
\end{proof}

For the case $c > 1$, we first need the following lemma.
\begin{lem}\label{SHc0}
	Suppose $c > 1$. Then, $-\check X(0)^{-1} \in S_H^c$.
\end{lem}

\begin{proof}
	From Lemma \ref{lemcheckX0}, we know that $\check X(0)$ is well-defined and $\check X(0) \in \R^*$. So, for $x_0 \in S_F^c$, $D_- H\left(-\check X(x_0)^{-1}\right) = 0$. $S_F^c$ is open, and $\check X$ is continuous - due to Theorem 2.2 \cite{Silverstein1995c} and non-constant, so $\{ -\check X(x_0)^{-1}, x_0 \in S_F^c \}$ is open. So, $\forall x_0 \in S_F^c, H'\left(-\check X(x_0)^{-1}\right)$ is defined and $H'\left(-\check X(x_0)^{-1}\right) = 0$. So, $\{ -\check X(x_0)^{-1}, x_0 \in S_F^c \} \subset S_H^c$. In particular, $-\check X(0)^{-1} \in S_H^c$.
\end{proof}

Now, we can finish the proof of the theorem.
\begin{lem}\label{csh0}
	Suppose $c > 1$. Then, 
	\begin{equation}\label{}
	\begin{aligned}
		\Delta(x) &= \int_{-\infty}^x h(\lambda) dF(\lambda) \text{ , with } h(0) = \frac{1}{(c-1)\check X(0)}.
	\end{aligned}
	\end{equation}
\end{lem}

\begin{proof}
	We define, for $z \in \C_+$, the following analytic function:
	\begin{equation}\label{}
	\begin{aligned}
		\mu(z) = \int \frac{\tau}{\tau X(z) + 1} = \frac{1 + zm(z)}{X(z)}.
	\end{aligned}
	\end{equation}
	Then,
	\begin{equation}\label{}
	\begin{aligned}
		\lim_{\varepsilon \rightarrow 0^+} \Delta(\varepsilon) - \Delta(-\varepsilon) = \lim_{\varepsilon \rightarrow 0^+} \lim_{\eta \rightarrow 0^+} \frac{1}{\pi} \int_{-\varepsilon}^{\varepsilon} \Im\left[-\frac{\mu(\xi + i\eta)}{\xi + i\eta}\right] d\xi.
	\end{aligned}
	\end{equation}
	Using Lemma \ref{SHc0}, we have by D.C.T that $\mu(z) \underset{z \in \C_+ \rightarrow 0}{\longrightarrow} \mu(0) = \frac{1}{c\check X(0)}$. So, by Lemma 3.5 \cite{Ledoit2009}, $\lim_{\varepsilon \rightarrow 0^+} \Delta(\varepsilon) - \Delta(-\varepsilon) = \frac{1}{c\check X(0)}$.

	By Lemma 3.3 \cite{Ledoit2009}, $F(\lambda) = \left(1 - \frac{1}{c}\right)\mathrm{1}_{[0,+\infty[}(\lambda)$ in a neighborhood $N$ of $0$. So, for $x \in N$:
	\begin{equation}\label{}
	\begin{aligned}
		&\Delta(x) = \int_{-\infty}^x \frac{1}{c \check X(0)}d\mathrm{1}_{[0,+\infty[}(\lambda)\\
		&\Delta(x) = \int_{-\infty}^x \frac{1}{(c-1) \check X(0)}dF(\lambda).
	\end{aligned}
	\end{equation}
	Outside $N$, the proof with $c < 1$ applies, hence the result.
\end{proof}

\subsection{Proof of Theorem \ref{psh}}
	\begin{proof}
		We begin by the case $c < 1$. By Lemma \ref{lem31}, $\Psi(x) := \lim_{n \rightarrow +\infty} \Psi_n(x)$ exists a.s. at all points of continuity of $\Psi$ and is equal to:
	\begin{equation}\label{}
	\begin{aligned}
		\Psi(x) = \lim_{\eta \rightarrow 0^+} \frac{1}{\pi} \int_{-\infty}^x \Im\left[\Theta^{(-1)}(\lambda+i\eta) \right]d\lambda.
	\end{aligned}
	\end{equation}

	We have, using Lemma \ref{lem32}:
	\begin{equation}\label{}
	\begin{aligned}
		\Psi(x) &= \lim_{\eta \rightarrow 0^+} \frac{1}{\pi} \int_{-\infty}^x \Im\left[\Theta_n^{(-1)}(\lambda+i\eta) \right]d\lambda \\
		\Psi(x) &= \int_{-\infty}^x \frac{1}{\pi}\Im[\check{\Theta}^{(-1)}(\lambda)]d\lambda.		
	\end{aligned}
	\end{equation}
	Remark that $\forall \lambda \in \R, \Im[\check{\Theta}^{(-1)}(\lambda)] = 0 \iff \Im[\check{m}(\lambda)] = 0$. So, as $F$ has a density $F' := \pi^{-1}\Im[\check m(\cdot)]$ (see \cite{Couillet2015}):
	\begin{equation}\label{}
	\begin{aligned}
		\Psi(x) &= \int_{]-\infty,x[\cap (F')^{-1}(\R_+^*)} \frac{1}{\pi}\Im[\check{\Theta}^{(-1)}(\lambda)]d\lambda \\
		 &= \int_{]-\infty,x[\cap (F')^{-1}(\R_+^*)} \frac{\Im[\check \Theta^{(-1)}(\lambda)]}{\Im[\check{m}(\lambda)]}\times \frac{1}{\pi}\Im[\check{m}(\lambda)]d\lambda \\	
		\Psi(x) &= \int_{]-\infty,x[\cap (F')^{-1}(\R_+^*)} \frac{\Im[\check \Theta^{(-1)}(\lambda)]}{\Im[\check{m}(\lambda)]}dF(\lambda).
	\end{aligned}
	\end{equation}
	Let $\lambda \in (F')^{-1}(\R_+^*)$. We proved that $F'(0) = \pi^{-1} \Im[\check m(0)] = 0$, so $\lambda \neq 0$. Then, using Theorem \ref{c0}:
	\begin{equation}\label{}
	\begin{aligned}
		\frac{\Im[\check \Theta^{(-1)}(\lambda)]}{\Im[\check{m}(\lambda)]} &= \frac{\Im\left[\int \frac{\tau^{-1}}{\tau \check X(\lambda) + 1}dH(\tau)\right]}{\Im\left[\int \frac{1}{\tau \check X(\lambda) + 1}dH(\tau)\right]} \\
		&= \frac{\int \frac{1}{|\tau \check X(\lambda) + 1|^2}dH(\tau)}{\int \frac{\tau}{|\tau \check X(\lambda) + 1|^2}dH(\tau)} \\
		\frac{\Im[\check \Theta^{(-1)}(\lambda)]}{\Im[\check{m}(\lambda)]} &= t(\lambda).
	\end{aligned}
	\end{equation}
	$t(\cdot)$ is well-defined on $\R^* \backslash \{ x\in S_F | F'(x) = 0 \}$, due to Theorem \ref{c0}. We extend $h$ on $\R$ with any constant value  $\lambda \in \{ x\in S_F | F'(x) = 0 \} \cup \{ 0 \}$, for example, for $t(\lambda) := 0$. The value does not matter as $dF\left((F')^{-1}(\{0\})\right) = 0$ and $dF\left(\{0\}\right) = 0$ as $c<1$.
	Finally:
	\begin{equation}\label{}
	\begin{aligned}
		\Psi(x) &= \int_{]-\infty,x[\cap (F')^{-1}(\R_+^*)} t(\lambda) dF(\lambda) \\
		\Psi(x) &= \int_{-\infty}^x t(\lambda) dF(\lambda),
	\end{aligned}
	\end{equation}
	which concludes the proof.
\end{proof}

Now, we can finish the proof of the theorem with the case $c>1$.
\begin{lem}\label{psh0}
	Suppose $c > 1$. Then, 
	\begin{equation}\label{}
	\begin{aligned}
		\Psi(x) &= \int_{-\infty}^x t(\lambda) dF(\lambda) \text{ , with } t(0) = -\check X(0) + \frac{c}{c-1} \int \frac{1}{\tau}dH(\tau).
	\end{aligned}
	\end{equation}
\end{lem}

\begin{proof}
	We define, for $z \in \C_+$, the following analytic function:
	\begin{equation}\label{}
	\begin{aligned}
		\mu(z) = zm(z)X(z) + \int \frac{1}{\tau}dH(\tau).
	\end{aligned}
	\end{equation}
	Then,
	\begin{equation}\label{}
	\begin{aligned}
		\lim_{\varepsilon \rightarrow 0^+} \Psi(\varepsilon) - \Psi(-\varepsilon) = \lim_{\varepsilon \rightarrow 0^+} \lim_{\eta \rightarrow 0^+} \frac{1}{\pi} \int_{-\varepsilon}^{\varepsilon} \Im\left[-\frac{\mu(\xi + i\eta)}{\xi + i\eta}\right] d\xi.
	\end{aligned}
	\end{equation}
	Using Lemma \ref{SHc0}, we have by D.C.T that $\mu(z) \underset{z \in \C_+ \rightarrow 0}{\longrightarrow} \mu(0) = -\check X(0)\frac{c}{c-1} + \int \frac{1}{\tau}dH(\tau)$. So, by Lemma 3.5 \cite{Ledoit2009}, $\lim_{\varepsilon \rightarrow 0^+} \Psi(\varepsilon) - \Psi(-\varepsilon) = -\check X(0)\frac{c}{c-1} + \int \frac{1}{\tau}dH(\tau)$.

	By Lemma 3.3 \cite{Ledoit2009}, $F(\lambda) = \left(1 - \frac{1}{c}\right)\mathrm{1}_{[0,+\infty[}(\lambda)$ in a neighborhood $N$ of $0$. So, for $x \in N$:
	\begin{equation}\label{}
	\begin{aligned}
		&\Psi(x) = \int_{-\infty}^x -\check X(0)\frac{c}{c-1} + \int \frac{1}{\tau}dH(\tau)d\mathrm{1}_{[0,+\infty[}(\lambda)\\
		&\Psi(x) = \int_{-\infty}^x -\check X(0) + \frac{c}{c-1} \int \frac{1}{\tau}dH(\tau)dF(\lambda).
	\end{aligned}
	\end{equation}
	Outside $N$, the proof with $c < 1$ applies, hence the result.
\end{proof}

\subsection{Proof of Theorem \ref{phithm}}
This proof follows the same ideas we developed in the previous proofs for shrinkage.

Assume \ref{H1}-\ref{H5}. Let $\tau \in \R$. From Lemma \ref{lem31} with $g = \mathrm{1}_{]-\infty,\tau]}$, we have that $\forall z \in \C_+$, almost surely, $\Theta_n^g \rightarrow \Theta^g(z)$ from Theorem \ref{FEg}. We also have: $\Phi_n(\lambda, \tau) = \lim_{\eta \rightarrow 0^+} \frac{1}{\pi} \int_{-\infty}^\lambda \Im[\Theta_n^g(\ell +i\eta)]d\ell $. From Lemma \ref{lem31}, we have that $\Phi(\lambda, \tau) := \lim_{n \rightarrow +\infty} \Phi_n(\lambda, \tau)$ exists and:
\begin{equation}\label{}
\begin{aligned}
	&\Phi(\lambda,\tau) = \lim_{\eta \rightarrow 0^+} \frac{1}{\pi} \int_{-\infty}^\lambda \Im[\Theta^g(\ell +i\eta)]d\ell ,
\end{aligned}
\end{equation}
for all $\lambda, \tau$ where $\Phi$ is continuous.

Firstly, suppose $c<1$. Then $F$ has a density over $\R$, and using Lemma \ref{lem32} and Theorem \ref{c0}:
\begin{equation}\label{}
\begin{aligned}
	\Phi(\lambda,\tau) &= \frac{1}{\pi} \int_{-\infty}^\lambda \lim_{\eta \rightarrow 0^+}\Im[\Theta^g(\ell +i\eta)]d\ell \\
	& =\frac{1}{\pi} \int_{-\infty}^\lambda \Im[\check \Theta^g(\ell +i\eta)]d\ell \\
	\Phi(\lambda,\tau) &= \int_{-\infty}^\lambda \int_{-\infty}^\tau \phi(\ell ,t) dH(t)dF(\ell ),
\end{aligned}
\end{equation}
using for the last line that $F'(\ell ) = \frac{1}{\pi}\Im[\check m(\ell )]$ and Theorem \ref{c0} for the characterization of $\Im[\check \Theta^g]$.

In the case $c>1$, following the same procedure as for shrinkage formulas, we define the following analytic function:
\begin{equation}\label{}
\begin{aligned}
	&\forall z \in \C_+, \mu(z) = \int_{-\infty}^\tau \frac{1}{t X(z) + 1}dH(t).
\end{aligned}
\end{equation}
We then have due to Lemma 3.5 \cite{Ledoit2009} and Lemma \ref{SHc0}:
\begin{equation}\label{}
\begin{aligned}
	&\lim_{\varepsilon \rightarrow 0^+} \Phi(\varepsilon,\tau) - \Phi(-\varepsilon,\tau) = \mu(0) = \int_{-\infty}^\tau \frac{1}{t \check X(0) + 1}dH(t).
\end{aligned}
\end{equation}
So, from the same argument used in shrinkage formulas, as in a neighborhood of $0$ we have $F(\lambda) = \left(1 - \frac{1}{c}\right)1_{[0,\infty[}(\lambda)$, we have :
\begin{equation}\label{}
\begin{aligned}
	&\Phi(\lambda,\tau) = \int_{-\infty}^\lambda \int_{-\infty}^\tau \phi(\ell ,t) dH(t)dF(\ell ), \\
	&\text{with } \phi(0,t) = \frac{c}{(c-1)(t \check X(0) + 1)},
\end{aligned}
\end{equation}
which concludes the proof.

\subsection{Proof of Theorem \ref{exp}}
Let $\alpha \in \R_+$, $\beta = \frac{\alpha}{1-e^{-\alpha}}$, and $z \in \C_+$. We remark that, noting $m := \Theta^{(0)}$:
	\begin{equation}\label{}
	\begin{aligned}
		1 + zm(z) = -zX(z) \Theta^{(1)}(z) \text{ and } -zX(z) = \int \frac{\delta}{1 + \delta c\Theta^{(1)}(z)}dD(\delta).
	\end{aligned}
	\end{equation}
	So, we deduce that:
	\begin{equation}\label{}
	\begin{aligned}
		1 + zm(z) = \Theta^{(1)}(z) \int \frac{\delta}{1 + \delta c \Theta^{(1)}(z)}dD(\delta).
	\end{aligned}
	\end{equation}
	In the case of $\alpha$-exponential weight law, we have, assuming complex-valued logarithm:
	\begin{equation}\label{}
	\begin{aligned}
		1 + zm(z) = \frac{1}{\alpha c} \left(\log \left[1 + \beta c \Theta^{(1)}(z)\right] - \log \left[1 + \beta e^{-1}c \Theta^{(1)}(z)\right] \right).
	\end{aligned}
	\end{equation}
	So,
	\begin{equation}\label{}
	\begin{aligned}
		\exp\left(\alpha c(1 + zm(z))\right) = \frac{1 + \beta c \Theta^{(1)}(z)}{1 + \beta e^{-1}c \Theta^{(1)}(z)}.
	\end{aligned}
	\end{equation}
	And,
	\begin{equation}\label{}
	\begin{aligned}
		\Theta^{(1)}(z) = \frac{e^{\alpha c(1 + zm(z))} - 1}{\beta c\left(1 - e^{-\alpha + \alpha c(1 + zm(z))}\right)}
	\end{aligned}
	\end{equation}
	Due to Theorem \ref{c0} and Theorem 2.2 \cite{Silverstein1995c}, $\check \Theta^{(1)}$ is continuous on $\R^*$. So, for all $x \in \R^*$, $c(1 + \lambda \check m(\lambda)) \neq 1$. So, by continuity, we have for all $\lambda \in \R^*$:
	\begin{equation}\label{}
	\begin{aligned}
		\check \Theta^{(1)}(\lambda) = \frac{e^{\alpha c(1 + \lambda\check m(\lambda))} - 1}{\beta c\left(1 - e^{-\alpha + \alpha c(1 + \lambda m(\lambda))}\right)},
	\end{aligned}
	\end{equation}
	which concludes the proof.

	As final remark, we have for $\lambda \in (F')^{-1}(\R_+^*)$, through the method of Lemma \ref{lempoly}:
	\begin{equation}\label{}
	\begin{aligned}
		\check \Theta^{(-1)}(\lambda) = - \check X(\lambda) \check m(\lambda) - \frac{1}{\lambda} \int \frac{1}{\tau}dH(\tau) = \frac{\check m(\lambda)(1+\lambda \check m(\lambda))}{\lambda \check \Theta^{(1)}(\lambda)} - \frac{1}{\lambda} \int \frac{1}{\tau}dH(\tau),
	\end{aligned}
	\end{equation}
	so, 
	\begin{equation}\label{}
	\begin{aligned}
		t(\lambda) = \frac{\Im\left[\frac{\check{m}(\lambda)(1+\lambda \check{m}(\lambda))}{\lambda \check{\Theta}^{(1)}(\lambda)}\right]}{\Im[\check{m}(\lambda)]}.
	\end{aligned}
	\end{equation}

	As for all $\lambda \in (F')^{-1}(\R_+^*)$, $h(\lambda) = \Im[\check{\Theta}^{(1)}(\lambda)]/\Im[\check{m}(\lambda)]$ and $t(\lambda) = \Im\left[\frac{\check{m}(\lambda)(1+\lambda \check{m}(\lambda))}{\lambda \check{\Theta}^{(1)}(\lambda)}\right]/\Im[\check{m}(\lambda)]$, Theorem \ref{exp} states that the statistical challenge for exponentially-weighted shrinkage is the same as it is for equally weighted sample covariance: estimate $\check{m}$, which can be done through kernel estimation for example in \cite{Jing2010}.

\section{Appendix: Proof of the numerical procedures}\label{apxE}
\subsection{Proof of Theorem \ref{checkX}}
Let us proof the equivalence in two steps. Assume \ref{H1}-\ref{H5} and let $\lambda \in \R^* \backslash \{x \in S_F, F'(x) = 0\}$.

We prove firstly that $f_\lambda(\check X(\lambda)) = 0$. From Lemma 3.2 \cite{Couillet2015}, we have that for all $R _subset \C_+$ bounded and bounded away from the imaginary axis, $z \in R \mapsto \int \frac{1}{|1 + \delta c\Theta^{(1)}(z)|^2}dD(\delta)$ is bounded. So, by Fatou's lemma, for $\lambda \in \R^*$, $\int \frac{1}{|1 + \delta c\check \Theta^{(1)}(\lambda)|^2}dD(\delta) < +\infty$. As $\forall z \in \C_+, X(z) = \frac{1}{z} \int \frac{\delta}{1+\delta c \Theta^{(1)}(z)}dD(\delta)$, we have by D.C.T. that, for $\lambda \in \R^*$, $f_z(X(z)) \underset{ z\in \C_+ \rightarrow \lambda}{\longrightarrow} f_\lambda(\check X(\lambda)) $. Hence, $f_\lambda(\check X(\lambda)) = 0$.

For the direct way, suppose that $\check X(\lambda) \in \C_+$. Then, by D.C.T, we have immediately that $f_{\lambda}\left(\check X(\lambda)\right) = 0$. So $f_\lambda(X) = 0$ has indeed at least one solution in $\C_+$.

For the other way, suppose now that there exists $X_0 \in \C_+$ such that $f_\lambda(X_0) = 0$. We are using the Theorem 9.1.1 \cite{Moller2015} to conclude. In order to apply it, we have that:
\begin{itemize}
	\item $f_\lambda(X_0) = 0$,
	\item $(x,X) \mapsto f_x(X)$ is analytic on an open of $\C^2$ containing $(\lambda, X_0)$,
	\item $f_\lambda(\cdot)$ is not identically $0$ in a neighborhood of $X_0$ (because of the imaginary part for example).
\end{itemize}
Then, there exists $\delta > 0$, $\varepsilon > 0$, $m \in \N^*$ such that for all $z \in \C$, $|z - x_0| < \varepsilon$, $f_z(\cdot)$ has exactly $m$ zeros in $\{ X \in \C, |X - X_0| < \delta\}$ and $X_{k_j}(z) = X_0 + \sum_{n=1}^\infty a_{k_n}\left( \left((z-\lambda)^{1/p_k} \right)_j \right)^n$. We refer to Section 9 \cite{Moller2015} for the notation.

So, for $|z-\lambda|$ small enough, $X_{k_j}(z) \in \C_+$ as $X_0 \in \C_+$ and $X_{k_j}$ continuous. Moreover, $f_z(X_{k_j}(z)) = 0$ for all $z\in \C_+$ and $|z-\lambda|$ small enough. By uniqueness of the solution of $f_z(\cdot) = 0$ when $z \in \C_+$, we deduce that $\forall k_j \in \llbracket 1,m\rrbracket, X_{k_j}(z) = X(z)$.

In conclusion, that $\forall k_j \in \llbracket 1,m\rrbracket, X_{k_j}(z) \underset{z \rightarrow x_0}{\longrightarrow} X_0$ by continuity of $X_{k_j}$, and $\forall k_j \in \llbracket 1,m\rrbracket, X_{k_j}(z) \underset{z \rightarrow x_0}{\longrightarrow} \check X(\lambda)$ by continuity of $X(\cdot)$ on $\C^*$. So, $\check X(\lambda) = X_0 \in \C_+$, which concludes the proof.

\bibliography{mybib}
\end{document}